\documentclass[11pt]{article}%
\usepackage{amssymb}
\usepackage{graphicx}
\usepackage{amsmath}
\usepackage{amsfonts}%
\providecommand{\U}[1]{\protect\rule{.1in}{.1in}}
\newtheorem{theorem}{Theorem}[section]

\newtheorem{corollary}[theorem]{Corollary}

\newtheorem{definition}[theorem]{Definition}

\newtheorem{lemma}[theorem]{Lemma}
\newtheorem{notation}[theorem]{Notation}

\newtheorem{proposition}[theorem]{Proposition}
\newtheorem{remark}[theorem]{Remark}

\newenvironment{proof}[1][Proof]{\textbf{#1.} }{\hfill\rule{0.5em}{0.5em}}
{\catcode`\@=11\global\let\AddToReset=\@addtoreset
\AddToReset{equation}{section}

\AddToReset{theorem}{section}

\begin{document}

\title{Self-similar solutions of the $p$-Laplace heat equation: the case $p>2.$}
\author{Marie Fran\c{c}oise Bidaut-V\'{e}ron\thanks{Laboratoire de Math\'{e}matiques
et Physique Th\'{e}orique, CNRS UMR 6083, Facult\'{e} des Sciences, Parc
Grandmont, 37200 Tours, France. e-mail:veronmf@univ-tours.fr}}
\maketitle
\date{.}

\begin{abstract}
We study the self-similar solutions of the equation
\[
u_{t}-div(\left\vert \nabla u\right\vert ^{p-2}\nabla u)=0,
\]
in $\mathbb{R}^{N},$ when $p>2.$ We make a complete study of the existence and
possible uniqueness of solutions of the form%
\[
u(x,t)=(\pm t)^{-\alpha/\beta}w((\pm t)^{-1/\beta}\left\vert x\right\vert )
\]
of any sign, regular or singular at $x=0.$ Among them we find solutions with
an expanding compact support or a shrinking hole (for $t>0),$ or a spreading
compact support or a focussing hole (for $t<0).$ When $t<0,$ we show the
existence of positive solutions oscillating around the particular solution
$U(x,t)=C_{N,p}(\left\vert x\right\vert ^{p}/(-t))^{1/(p-2)}.$

\end{abstract}

.\pagebreak

\section{ Introduction and main results\label{intro}}

Here we consider the self-similar solutions of the degenerate heat equation
involving the $p$-Laplace operator
\[
\qquad\qquad\qquad\qquad\qquad\qquad\qquad u_{t}-div(\left\vert \nabla
u\right\vert ^{p-2}\nabla u)=0.\qquad\qquad\qquad\qquad\qquad\qquad
\qquad\qquad\text{(\textbf{E}}_{u}\text{)}%
\]
in $\mathbb{R}^{N},$ with $p>2.$ This study is the continuation of the work
started in \cite{Bi1}, relative to the case $p<2.$ It can be read
independently. We set
\begin{equation}
\gamma=\frac{p}{p-2},\qquad\eta=\frac{N-p}{p-1}, \label{del}%
\end{equation}
thus $\gamma>1,$ $\eta<N,$
\begin{equation}
\frac{N+\gamma}{p-1}=\eta+\gamma=\frac{N-\eta}{p-2}. \label{rel}%
\end{equation}
$\medskip$

If $u$ is a solution, then for any $\alpha,\beta\in\mathbb{R},$ $u_{\lambda
}(x,t)=\lambda^{\alpha}u(\lambda x,\lambda^{\beta}t)$ is a solution of
(\textbf{E}$_{u}$) if and only if
\begin{equation}
\beta=\alpha(p-2)+p=(p-2)(\alpha+\gamma); \label{bet}%
\end{equation}
notice that $\beta>0\Longleftrightarrow\alpha>-\gamma.$ Given $\alpha
\in\mathbb{R}$ such that $\alpha\neq-\gamma,$ we search self-similar
solutions, radially symmetric in $x,$ of the form:%
\begin{equation}
u=u(x,t)=(\varepsilon\beta t)^{-\alpha/\beta}w(r),\qquad r=(\varepsilon\beta
t)^{-1/\beta}\left\vert x\right\vert , \label{chg}%
\end{equation}
where $\varepsilon=\pm1.$ By translation, for any real $T,$ we obtain
solutions defined for any $t>T$ when $\varepsilon\beta>0,$ or $t<T$ when
$\varepsilon\beta<0.$ We are lead to the equation%
\[
\qquad\qquad\qquad\left(  \left\vert w^{\prime}\right\vert ^{p-2}w^{\prime
}\right)  ^{\prime}+\frac{N-1}{r}\left\vert w^{\prime}\right\vert
^{p-2}w^{\prime}+\varepsilon(rw^{\prime}+\alpha w)=0\qquad\text{in }\left(
0,\infty\right)  .\qquad\qquad\qquad\text{(\textbf{E}}_{w}\text{)}%
\]
\textit{ }

Our purpose is to give a complete description of all the solutions, with
constant or changing sign. Equation (\textbf{E}$_{w}$) is very interesting,
because it is singular at any zero of $w^{\prime},$ since $p>2,$ implying a
nonuniqueness phenomena.\medskip

For example, concerning the constant sign solutions near the origin, it can
happen that
\[
\lim_{r\rightarrow0}w=a\neq0,\quad\quad\lim_{r\rightarrow0}w^{\prime}=0,
\]
we will say that $w$ is \textit{regular}, or
\[
\lim_{r\rightarrow0}w=\lim_{r\rightarrow0}w^{\prime}=0,
\]
we say that $w$ is \textit{flat}. Or different kinds of singularities may
occur, either at the level of $w:$%
\[
\lim_{r\rightarrow0}w=\infty,
\]
\medskip or at the level of the gradient:
\begin{align*}
\lim_{r\rightarrow0}w  &  =a\in\mathbb{R},\quad\quad\lim_{r\rightarrow
0}w^{\prime}=\pm\infty,\quad\quad\quad\text{when }p>N>1,\\
\lim_{r\rightarrow0}w  &  =a\in\mathbb{R},\quad\quad\lim_{r\rightarrow
0}w^{\prime}=b\neq0\quad\quad\text{when }p>N=1.
\end{align*}

We first show that any local solution $w$ of (\textbf{E}$_{w}$) can be defined
on $\left(  0,\infty\right)  ,$ thus any solution $u$ of equation
(\textbf{E}$_{u}$) associated to $w$ by (\ref{chg}) is defined on
$\mathbb{R}^{N}\backslash\left\{  0\right\}  \times\left(  0,\pm\infty\right)
.$ Then we prove the existence of regular solutions, flat ones, and of all
singular solutions mentioned above.\medskip

Moreover, for $\varepsilon=1$, there exist solutions $w$ with a compact
support $\left(  0,\bar{r}\right)  $; then $u\equiv0$ on the set
\[
D=\left\{  \left(  x,t\right)  :x\in\mathbb{R}^{N},\quad\beta t>0,\quad
\left\vert x\right\vert >(\beta t)^{1/\beta}\bar{r}\right\}  .
\]
For $\varepsilon=-1,$ there exist solutions with a hole:
$w(r)=0\Longleftrightarrow r\in\left(  0,\bar{r}\right)  $. Then $u\equiv0$ on
the set%
\[
H=\left\{  \left(  x,t\right)  :x\in\mathbb{R}^{N},\quad\beta t<0,\quad
\left\vert x\right\vert <(-\beta t)^{1/\beta}\bar{r}\right\}  .
\]
The free boundary is of parabolic type for $\beta>0,$ of hyperbolic type for
$\beta<0.$ This leads to four types of solutions, and we prove their
existence:\medskip

$\bullet$ If $t>0$, with $\varepsilon=1,\beta>0$, we say that $u$ has an
\textit{expanding support}; the support increases from $\left\{  0\right\}  $
as $t$ increases from $0$. \medskip

$\bullet$ If $t>0$, with $\varepsilon=-1,\beta<0$, we say that $u$ has a
\textit{shrinking hole}: the hole decreases from infinity as $t$ increases
from $0$;\medskip

$\bullet$ If $t<0$, with $\varepsilon=1,\beta<0$, we say that $u$ has a
\textit{spreading support}: the support increases to be infinite as $t$
increases to $0.$\medskip\ 

$\bullet$ If $t<0$, with $\varepsilon=-1,\beta>0,$ we say that $u$ has a
\textit{focussing hole}: the hole disappears as $t$ increases to $0.$ \medskip

Up to our knowledge, some of them seem completely new, as for example the
solutions with a shrinking hole or a spreading support. In particular we find
again and improve some results of \cite{GiVa} concerning the existence of
focussing type solutions. \medskip

Finally for $t<0$ we also show the existence of positive solutions turning
around the fundamental solution $U$ given at (\ref{uinf}) with a kind of
periodicity, and also the existence of changing sign solutions doubly
oscillating in $\left\vert x\right\vert $ near $0$ and infinity.\medskip

As in \cite{Bi1} we reduce the problem to dynamical systems.\medskip

When $\varepsilon=-1,$ a critical negative value of $\alpha$ is involved:
\begin{equation}
\alpha^{\ast}=-\gamma+\frac{\gamma(N+\gamma)}{(p-1)(N+2\gamma)}. \label{star}%
\end{equation}

\subsection{Explicit solutions}

Obviously if $w$ is a solution of (\textbf{E}$_{w}$), $-w$ is also a solution.
Some particular solutions are well-known. \medskip

\noindent\ \textbf{The solution }$U$. For any $\alpha$ such that
$\varepsilon(\alpha+\gamma)<0,$ that means $\varepsilon\beta<0,$ there exist
flat solutions of (\textbf{E}$_{w}$), given by
\begin{equation}
w(r)=\pm\ell r^{\gamma}, \label{rdel}%
\end{equation}
where%
\begin{equation}
\ell=\left(  \frac{\left\vert \alpha+\gamma\right\vert }{\gamma^{p-1}%
(\gamma+N)}\right)  ^{1/(p-2)}>0. \label{defl}%
\end{equation}
They correspond to a unique solution of (\textbf{E}$_{u})$ called $U$, defined
for $t<0,$ such that $U(0,t)=0,$ flat, blowing up at $t=0$ for fixed
$x\neq0:$
\begin{equation}
U(x,t)=C\left(  \frac{\left\vert x\right\vert ^{p}}{-t}\right)  ^{1/(p-2)}%
,\qquad C=((p-2)\gamma^{p-1}(\gamma+N))^{1/(2-p)}. \label{uinf}%
\end{equation}
\medskip

\noindent\ \textbf{The case }$\alpha=N.$ Then $\beta=\beta_{N}=N(p-2)+p>0,$
and the equation has a first integral%
\begin{equation}
w+\varepsilon r^{-1}\left\vert w^{\prime}\right\vert ^{p-2}w^{\prime}=Cr^{-N}.
\label{fat}%
\end{equation}
All the solutions corresponding to $C=0$ are given by
\begin{align}
w  &  =w_{K,\varepsilon}(r)=\pm\left(  K-\varepsilon\gamma^{-1}r^{p^{\prime}%
}\right)  _{+}^{(p-1)/(p-2)},\qquad K\in\mathbb{R},\nonumber\\
u  &  =\pm u_{K,\varepsilon}(x,t)=\pm(\varepsilon\beta_{N}t)^{-N/\beta_{N}%
}\left(  K-\varepsilon\gamma^{-1}(\varepsilon\beta_{N}t)^{-p^{\prime}%
/\beta_{N}}\left\vert x\right\vert ^{p^{\prime}}\right)  _{+}^{(p-1)/(p-2)}.
\label{fut}%
\end{align}
For $\varepsilon=1,$ $K>0,$ they are defined for $t>0,$ called
\textit{Barenblatt solutions,} regular with a compact support. Given $c>0,$
the function $u_{K,1},$ defined on $\mathbb{R}^{N}\times\left(  0,\infty
\right)  ,$ is the unique solution of equation (\textbf{E}$_{u}$) with initial
data $u(0)=c\delta_{0},$ where $\delta_{0}$ is the Dirac mass at $0,$ and $K$
being linked by $%
{\displaystyle\int\limits_{\mathbb{R}^{N}}}
u_{K}(x,t)dt=c$. The $u_{K,1}$ are the only nonnegative solutions defined on
$\mathbb{R}^{N}\times\left(  0,\infty\right)  ,$ such that $u(x,0)=0$ for any
$x\neq0.$ For $\varepsilon=-1,$ the $u_{K,-1}$ are defined for $t<0;$ for
$K>0,$ $w$ does not vanish on $\left(  0,\infty\right)  ;$ for $K<0,$ $w$ is
flat with a hole near $0$. For $K=0,$ we find again the function $w$ given at
(\ref{rdel}).\medskip

\noindent\ \textbf{The case }$\alpha=\eta\neq0$. We exhibit a family of
solutions of (\textbf{E}$_{w}):$
\begin{equation}
w(r)=Cr^{-\eta},\qquad u(t,x)=C\left\vert x\right\vert ^{-\eta},\qquad C\neq0.
\label{aeta}%
\end{equation}
The solutions $u,$ independent of $t,$ are $p$-harmonic in $\mathbb{R}^{N}$;
they are fundamental solutions when $p<N$. When $p>N,$ $w$ satisfies
$\lim_{r\rightarrow0}w=0,$ and $\lim_{r\rightarrow0}w^{\prime}=\infty$ for
$N>1,$ $\lim_{r\rightarrow0}w^{\prime}=b$ for $N=1.$\medskip

\noindent\textbf{The case }$\alpha=-p^{\prime}.$ Equation (\textbf{E}$_{w}$)
admits regular solutions of the form%
\begin{equation}
w(r)=\pm K\left(  N(Kp^{\prime})^{p-2}+\varepsilon r^{p^{\prime}}\right)
,\qquad u(x,t)=\pm K\left(  N(Kp^{\prime})^{p-2}t+\left\vert x\right\vert
^{p^{\prime}}\right)  ,\qquad K>0. \label{moi}%
\end{equation}
Here $\beta>0;$ in the two cases $\varepsilon=1,t>0$ and $\varepsilon=-1,t<0,$
$u$ is defined for any $t\in\mathbb{R}$ and of the form $\psi(t)+\Phi
(\left\vert x\right\vert )$ with $\Phi$ nonconstant, and $u(.,t)$ has a
constant sign for $t>0$ and changing sign for $t<0.$\ \medskip

\noindent\textbf{The case }$\alpha=0.$ Equation (\textbf{E}$_{w}$) can be
explicitely solved: either $w^{\prime}\equiv0,$ thus $w\equiv a\in\mathbb{R},$
$u$ is a constant solution of (\textbf{E}$_{u}),$ or there exists
$K\in\mathbb{R}$ such that
\begin{equation}
\left\vert w^{\prime}\right\vert =r^{-(\eta+1)}\left(  K-\frac{\varepsilon
}{\gamma+N}r^{N-\eta}\right)  _{+}^{1/(p-2)}; \label{azer}%
\end{equation}
and $w$ follows by integration, up to a constant, and then
$u(x,t)=w(\left\vert x\right\vert /(\varepsilon pt)^{1/p}).$ If $\varepsilon
=1,$ then $t>0,$ $K>0$ and $w^{\prime}$ has a compact support; up to a
constant, $u$ has a compact support. If $\varepsilon=-1,$ then $t<0;$ for
$K>0,$ $w$ is strictly monotone; for $K<0,$ $w$ is flat, constant near $0;$
for $K=0,$ we find again (\ref{rdel}). For $\varepsilon=\pm1,K>0,$ observe
that $\lim_{r\rightarrow0}w=\pm\infty$ if $p\leqq N;$ and $\lim_{r\rightarrow
0}w=a\in\mathbb{R},$ $\lim_{r\rightarrow0}w^{\prime}=\pm\infty$ if $p>N>1;$
and $\lim_{r\rightarrow0}w=a\in\mathbb{R},$ $\lim_{r\rightarrow0}w^{\prime}=K$
if $p>N=1.$ In particular we find solutions such that $w=cr^{\left\vert
\eta\right\vert }(1+o(1))$ near $0,$ with $c>0.$ \medskip\ 

\textbf{(v) Case }$N=1$\textbf{ and }$\alpha=-(p-1)/(p-2)<0.$ Here $\beta=1,$
and we find the solutions
\begin{equation}
w(r)=\pm\left(  Kr+\varepsilon\left\vert \alpha\right\vert ^{p-1}\left\vert
K\right\vert ^{p}\right)  _{+}^{(p-1)/(p-2)},\quad u(x,t)=\pm\left(
K\left\vert x\right\vert +\left\vert \alpha\right\vert ^{p-1}\left\vert
K\right\vert ^{p}t\right)  _{+}^{(p-1)/(p-2)}, \label{exp}%
\end{equation}
If $\varepsilon=1,t>0,$ then $w$ has a singularity at the level of the
gradient, and either $K>0,$ $w>0,$ or $K<0$ and $w$ has a compact support$.$
If $\varepsilon=-1,t<0$ then $K>0$, $w$ has a hole.

\subsection{Main results}

In the next sections we provide an exhaustive study of equation (\textbf{E}%
$_{w}$). Here we give the main results relative to the function $u.$ Let us
show how to return from $w$ to $u$. Suppose that the behaviour of $w$ is given
by
\[
\lim_{r\rightarrow0}r^{\lambda}w(r)=c\neq0,\qquad\lim_{r\rightarrow\infty
}r^{\mu}w(r)=c^{\prime}\neq0,\qquad\text{where }\lambda,\mu\in\mathbb{R}.
\]

\noindent(i) Then for fixed $t\neq0,$ the function $u$ has a behaviour in
$\left\vert x\right\vert ^{-\lambda}$ near $x=0,$ and a behaviour in
$\left\vert x\right\vert ^{-\mu}$ for large $\left\vert x\right\vert
.\medskip$

If $\lambda=0,$ then $u$ is defined on $\mathbb{R}^{N}\times\left(
0,\pm\infty\right)  .$ Either $w$ is regular, then $u(.,t)\in C^{1}\left(
\mathbb{R}^{N}\times\left(  0,\infty\right)  \right)  ;$ we will say that $u$
is \textbf{regular}; nevertheless the regular solutions $u$ presents a
singularity at time $t=0$ if and only if $\alpha<-\gamma$ or $\alpha>0.$ Or a
singularity can appear for $u$ at the level of the gradient.$\medskip$

If $\lambda<0,$ thus $u$ is defined on $\mathbb{R}^{N}\times\left(
0,\pm\infty\right)  $ and $u(0,t)=0$; either $w$ is flat, we also say that $u$
is \textbf{flat}, or a singularity appears at the level of the
gradient.$\medskip$

If $0<\lambda<N,$ then $u(.,t)\in L_{loc}^{1}\left(  \mathbb{R}^{N}\right)  $
for $t\neq0,$ we say that $x=0$ is a \textbf{weak}\textit{ singularity. }We
will show that there exist no stronger singularity.$\medskip$

If $\lambda<N<\mu;$ then $u(.,t)\in L^{1}\left(  \mathbb{R}^{N}\right)
.\medskip$

\noindent(ii) For fixed $x\neq0,$ the behaviour of $u$ near $t=0,$ depends on
the sign of $\beta$:
\begin{align*}
\lim_{t\rightarrow0}\left\vert x\right\vert ^{\mu}\left\vert t\right\vert
^{(\alpha-\mu)/\beta}u(x,t)  &  =C\neq0\quad\text{if}\quad\alpha>-\gamma,\\
\lim_{t\rightarrow0}\left\vert x\right\vert ^{\lambda}\left\vert t\right\vert
^{(\alpha-\lambda)/\beta}u(x,t)  &  =C\neq0\quad\text{if}\quad\alpha<-\gamma.
\end{align*}
If $\mu<0,\alpha>-\gamma$ or $\lambda<0,\alpha<-\gamma$, then $\lim
_{t\rightarrow0}u(x,t)=0.$

\subsubsection{Solutions defined for $t>0$}

Here we look for solutions $u$ of (\textbf{E}$_{u})$ of the form (\ref{chg})
defined on $\mathbb{R}^{N}\backslash\left\{  0\right\}  \times\left(
0,\infty\right)  .$ That means $\varepsilon\beta>0$ or equivalently
$\varepsilon=1,$ $-\gamma<\alpha$ (see Section \ref{one}) or $\varepsilon
=-1,\alpha<-\gamma$ see (Section \ref{two}). We begin by the case
$\varepsilon=1,$ treated at Theorem \ref{pin}.

\begin{theorem}
\label{T1}Assume $\varepsilon=1,$ and $-\gamma<\alpha.\medskip$

\noindent(1) Let $\alpha<N.\medskip$

All regular solutions on $\mathbb{R}^{N}\backslash\left\{  0\right\}
\times\left(  0,\infty\right)  $ have a strict constant sign, in $\left\vert
x\right\vert ^{-\alpha}$ near $\infty$ for fixed $t,$ with initial data
$L\left\vert x\right\vert ^{-\alpha}(L\neq0)$ in $\mathbb{R}^{N};$ thus
$u(.,t)\not \in L^{1}\left(  \mathbb{R}^{N}\right)  $, and $u$ is unbounded
when $\alpha<0.$

There exist nonnegative solutions such that near $x=0,$%
\begin{equation}
\left.
\begin{array}
[c]{c}%
\text{for }p<N,\text{\quad\ }u\text{ has a weak singularity in }\left\vert
x\right\vert ^{-\eta},\qquad\qquad\qquad\qquad\qquad\qquad\quad\\
\text{for }p=N,\text{ \quad}u\text{ has a weak singularity in }\ln\left\vert
x\right\vert ,\qquad\qquad\qquad\qquad\qquad\qquad\quad\\
\text{for }p>N,\text{ \quad}u\in C^{0}(\mathbb{R}^{N}\times\left(
0,\infty\right)  ,\quad u(0,t)=a>0,\text{ with a singular gradient},
\end{array}
\right\}  \label{rr}%
\end{equation}
and $u$ has an \textbf{expanding} \textbf{compact support} for any $t>0$, with
initial data $L\left\vert x\right\vert ^{-\alpha}in$ $\mathbb{R}^{N}%
\backslash\left\{  0\right\}  .$

There exist positive solutions with the same behaviour as $x\rightarrow0,$ in
$\left\vert x\right\vert ^{-\alpha}$ near $\infty$ for fixed $t;$ and also
solutions such that $u$ has one zero for fixed $t\neq0,$ and the same behaviour.

If $p>N,$ there exist positive solutions satisfying (\ref{rr}), and also
positive solutions such that%
\begin{equation}
u\in C^{0}(\mathbb{R}^{N}\times\left(  0,\infty\right)  ,\quad u(0,t)=0,\text{
in }\left\vert x\right\vert ^{\left\vert \eta\right\vert }\text{ near
}0,\text{ with a singular gradient,} \label{ro}%
\end{equation}
in $\left\vert x\right\vert ^{-\alpha}$ near $\infty$ for fixed $t,$ with and
initial data $L\left\vert x\right\vert ^{-\alpha}in$ $\mathbb{R}^{N}%
\backslash\left\{  0\right\}  .\medskip$

\noindent(2) Let $\alpha=N.\medskip$

All \textbf{regular (Barenblatt}) solutions are nonnegative, have a
\textbf{compact support} for any $t>0$. If $p\leqq N,$ all the other solutions
have one zero for fixed $t$, satisfy (\ref{rr}) or (\ref{ro}) and have the
same behaviour at $\infty.\medskip$

\noindent(3) Let $N<\alpha.\medskip$

All regular solutions $u$ have a finite number $m\geqq1$ of simple zeros for
fixed $t$, and $u(.,t)\in L^{1}\left(  \mathbb{R}^{N}\right)  .$ Either they
are in $\left\vert x\right\vert ^{-\alpha}$ near $\infty$ for fixed $t,$ then
there exist solutions with $m$ zeros, compact support, satisfying (\ref{rr});
or they have a compact support.$\ $All the solutions have $m$ or $m+1$ zeros.
There exist solutions satisfying (\ref{rr}) with $m+1$ zeros, and in
$\left\vert x\right\vert ^{-\alpha}$ near $\infty.$ If $p>N,$ there exist
solutions satisfying (\ref{rr}) with $m$ zeros; there exist also solutions
with $m$ zeros, $u(0,t)=0,$ and a singular gradient, in $\left\vert
x\right\vert ^{-\alpha}$ near $\infty.$
\end{theorem}

Next we come to the case $\varepsilon=-1,$ which is the subject of Theorem
\ref{osc}.

\begin{theorem}
\label{T2}Assume $\varepsilon=-1$ and $\alpha<-\gamma.$\medskip

All the solutions $u$ on $\mathbb{R}^{N}\backslash\left\{  0\right\}
\times\left(  0,\infty\right)  ,$ in particular the regular ones, are
\textbf{oscillating around} $0$ for fixed $t>0$ and large $\left\vert
x\right\vert ,$ and $r^{-\gamma}w$ is asymptotically periodic in $\ln r$.
Moreover there exist\medskip

solutions such that $r^{-\gamma}w$ is \textbf{periodic} in $\ln r,$ in
particular $C_{1}t^{-\left\vert \alpha/\beta\right\vert }\leqq\left\vert
u\right\vert \leqq C_{2}t^{-\left\vert \alpha/\beta\right\vert }$ for some
$C_{1},C_{2}>0;$

solutions $u\in C^{1}(\mathbb{R}^{N}\times\left[  0,\infty\right)  ),$
$u(x,0)\equiv0,$ with a \textbf{shrinking} \textbf{hole}$;$

\textbf{flat} solutions $u\in C^{1}(\mathbb{R}^{N}\times\left[  0,\infty
\right)  ),$ in $\left\vert x\right\vert ^{\left\vert \alpha\right\vert }$
near $0,$ with initial data $L\left\vert x\right\vert ^{\left\vert
\alpha\right\vert }(L\neq0);$

solutions satisfying (\ref{rr}) near $x=0,$ and if $p>N,$ solutions satisfying
(\ref{ro}) near $0.$
\end{theorem}

\subsubsection{Solutions defined for $t<0$}

We look for solutions $u$ of (\textbf{E}$_{u})$ of the form (\ref{chg})
defined on $\mathbb{R}^{N}\backslash\left\{  0\right\}  \times\left(
-\infty,0\right)  .$ That means $\varepsilon\beta<0$ or equivalently
$\varepsilon=1,$ $\alpha<-\gamma$ (see Section \ref{three}, Theorem \ref{mel})
or $\varepsilon=-1,\alpha>-\gamma$ (see Section \ref{four}). In the case
$\varepsilon=1,$ we get the following:

\begin{theorem}
\label{T3}Assume $\varepsilon=1,$ and $\alpha<-\gamma$.\medskip

The function $U(x,t)=C\left(  \frac{\left\vert x\right\vert ^{p}}{-t}\right)
^{1/(p-2)}$ is a positive \textbf{flat} solution on $\mathbb{R}^{N}%
\backslash\left\{  0\right\}  \times\left(  -\infty,0\right)  $.

All regular solutions have a constant sign, are unbounded in $\left\vert
x\right\vert ^{\gamma}$ near $\infty$ for fixed $t,$ and blow up at $t=0$ like
$\left(  -t\right)  ^{-\left\vert \alpha\right\vert /\left\vert \beta
\right\vert }$ for fixed $x\neq0.$

There exist \textbf{flat} \textbf{positive} solutions $u\in C^{1}%
(\mathbb{R}^{N}\times\left(  -\infty,0\right]  ),$ in $\left\vert x\right\vert
^{\gamma}$ near $\infty$ for fixed $t,$ with \textbf{final data} $L\left\vert
x\right\vert ^{\left\vert \alpha\right\vert }$ $(L>0)$.

There exist \textbf{nonnegative} solutions satisfying (\ref{rr}) near $0,$
with a \textbf{spreading} \textbf{compact support}, blowing up near $t=0$
(like $\left\vert t\right\vert ^{-(\eta+\left\vert \alpha\right\vert
)/\left\vert \beta\right\vert }$ for $p<N,$ or $\left\vert t\right\vert
^{-\left\vert \alpha\right\vert /\left\vert \beta\right\vert }\ln\left\vert
t\right\vert $ for $p=N$, or $\left(  -t\right)  ^{-\left\vert \alpha
\right\vert /\left\vert \beta\right\vert }$for $p>N).$

There exist positive solutions with the same behaviour near $0,$ in
$\left\vert x\right\vert ^{\gamma}$ near $\infty,$ blowing up as above at
$t=0$, and solutions with one zero for fixed $t$, and the same behaviour. If
$p>N,$ there exist positive solutions satisfying (\ref{rr}) (resp. (\ref{ro}))
near 0, in $\left\vert x\right\vert ^{\gamma}$ near $\infty$ for fixed $t,$
blowing up at $t=0$ like $\left\vert t\right\vert ^{-\left\vert \alpha
\right\vert /\left\vert \beta\right\vert }$ ( resp. $\left\vert t\right\vert
^{(\left\vert \eta\right\vert -\left\vert \alpha\right\vert )/\left\vert
\beta\right\vert }$) for fixed $x$.

Up to a symmetry, all the solutions are described.
\end{theorem}

The most interesting case is $\varepsilon=-1,-\gamma<\alpha.$ For simplicity
we will assume that\textbf{ }$p<N.$ The case $p\geqq N$ is much more delicate,
and the complete results can be read in terms of $w$ at Theorems \ref{int},
\ref{pom}, \ref{clin}, \ref{sou}, \ref{orb} and \ref{ent}. We discuss
according to the position of $\alpha$ with respect to $-p^{\prime}$ and
$\alpha^{\ast}$ defined at (\ref{star}). Notice that $\alpha^{\ast}%
<-p^{\prime}.$

\begin{theorem}
\label{T4}Assume $\varepsilon=-1,$ and $-p^{\prime}\leqq\alpha\neq0$. The
function $U$ is still a flat solution on $\mathbb{R}^{N}\backslash\left\{
0\right\}  \times\left(  -\infty,0\right)  .\medskip$

(1) Let $0<\alpha.\medskip$

All regular solutions have a strict constant sign, in $\left\vert x\right\vert
^{\gamma}$ near $\infty$ for fixed $t,$ blowing up at $t=0$ like $\left(
-t\right)  ^{-1/(p-2)}$ for fixed $x\neq0.$

There exist nonnegative solutions with a \textbf{focussing} \textbf{hole}:
$u(x,t)\equiv0$ for $\left\vert x\right\vert \leqq C\left\vert t\right\vert
^{1/\beta},$ $t>0,$ in $\left\vert x\right\vert ^{\gamma}$ near $\infty$ for
fixed $t,$ blowing up at $t=0$ like $\left(  -t\right)  ^{-1/(p-2)}$ for fixed
$x\neq0.$

There exist positive solutions $u$ with a (\textbf{weak}) \textbf{singularity}
in $\left\vert x\right\vert ^{-\eta}$at $x=0,$ in $\left\vert x\right\vert
^{-\alpha}$ near $\infty$ for fixed $t,$ with \textbf{ }$u(.,t)\in
L^{1}\left(  \mathbb{R}^{N}\right)  $ if $\alpha>N,$ with final data
$L\left\vert x\right\vert ^{-\alpha}$ $(L>0)$ in $\mathbb{R}^{N}%
\backslash\left\{  0\right\}  .$

There exist positive solutions $u$ in $\left\vert x\right\vert ^{-\eta}$at
$x=0,$ in $\left\vert x\right\vert ^{\gamma}$ near $\infty$ for fixed $t,$
blowing up at $t=0$ like $\left(  -t\right)  ^{-1/(p-2)}$ for fixed $x\neq0;$
solutions with one zero and the same behaviour.$\medskip$

\noindent(2) Let $-p^{\prime}<\alpha<0.\medskip$

All regular solutions have \textbf{one zero }for fixed $t$, and the same
behaviour. There exist solutions with one zero, in $\left\vert x\right\vert
^{-\eta}$at $x=0,$ in $\left\vert x\right\vert ^{\left\vert \alpha\right\vert
}$ near $\infty$ for fixed $t,$ with final data $L\left\vert x\right\vert
^{-\alpha}$ $(L>0)$ $in$ $\mathbb{R}^{N}\backslash\left\{  0\right\}  .$ There
exist solutions with one zero, $u$ in $\left\vert x\right\vert ^{-\eta}$at
$x=0,$ in $\left\vert x\right\vert ^{\gamma}$ near $\infty$ for fixed $t,$
blowing up at $t=0$ like $\left(  -t\right)  ^{-1/(p-2)}$ for fixed $x\neq0;$
solutions with two zeros and the same behaviour.$\medskip$

\noindent3) Let $\alpha=-p^{\prime}.\medskip$

All regular solutions have \textbf{one zero }and are in $\left\vert
x\right\vert ^{\left\vert \alpha\right\vert }$ near $\infty$ for fixed $t,$
and with \textbf{final data} $L\left\vert x\right\vert ^{\left\vert
\alpha\right\vert }$ $(L>0)$. The other solutions have one or two zeros, are
in $\left\vert x\right\vert ^{-\eta}$at $x=0,$ in $\left\vert x\right\vert
^{\gamma}$ near $\infty$ for fixed $t.$

In any case, up to a symmetry, all the solutions are described.
\end{theorem}

\begin{theorem}
\label{T5}Assume $\varepsilon=-1,-\gamma<\alpha<-p^{\prime}.$ Then $U$ is
still a flat solution on $\mathbb{R}^{N}\backslash\left\{  0\right\}
\times\left(  -\infty,0\right)  .\medskip$

(1)Let $\alpha\leqq\alpha^{\ast}.\medskip$

Then there exist \textbf{positive flat solutions}, in $\left\vert x\right\vert
^{\gamma}$ near 0, in $\left\vert x\right\vert ^{\left\vert \alpha\right\vert
}$ near $\infty$ for fixed $t,$ with \textbf{final data} $L\left\vert
x\right\vert ^{-\alpha}$ $(L>0)$ in $\mathbb{R}^{N}.$

All the other solutions, among them the \textbf{regular ones}, have an
\textbf{infinity of zeros}: $u(t,.)$ is oscillating around $0$ for large
$\left\vert x\right\vert .$ There exist solutions with a focussing hole, and
solutions with a singularity in $\left\vert x\right\vert ^{-\eta}$at $x=0.$
There exist solutions \textbf{oscillating also for small }$\left\vert
x\right\vert ,$ such that $r^{-\gamma}w$ is periodic in $\ln r.\medskip$

(2) There exist a \textbf{critical unique value} $\alpha_{c}\in\left(
\max(\alpha^{\ast},-p^{\prime}\right)  $ such that for $\alpha=\alpha_{c},$
there exists nonnegative solutions with a \textbf{focussing hole} near 0, in
$\left\vert x\right\vert ^{\left\vert \alpha\right\vert }$ near $\infty$ for
fixed $t,$ with \textbf{final data} $L\left\vert x\right\vert ^{-\alpha}$
$(L>0)$ in $\mathbb{R}^{N}.$ And $\alpha_{c}>-(p-1)/(p-2).$

\noindent There exist positive flat solutions, such that $\left\vert
x\right\vert ^{-\gamma}u$ is bounded on $\mathbb{R}^{N}$ for fixed $t,$
blowing up at $t=0$ like $\left(  -t\right)  ^{-1/(p-2)}$ for fixed $x\neq0.$
The regular solutions are oscillating around $0$ as above. There exist
solutions \textbf{oscillating around} $0,$ such that $r^{-\gamma}w$ is
\textbf{periodic} in $\ln r.$ There are solutions with a weak singularity in
$\left\vert x\right\vert ^{-\eta}$at $x=0$, and oscillating around 0 for large
$\left\vert x\right\vert .\medskip$

(3) Let $\alpha^{\ast}<\alpha<\alpha_{c}.\ \medskip$

\noindent The regular solutions are as above. There exist solutions of the
same types as above. Moreover there exist \textbf{positive} solutions, such
that $r^{-\gamma}w$ is \textbf{periodic} in $\ln r,$ thus there exist
$C_{1},C_{2}>0$ such that
\[
C_{1}\left(  \frac{\left\vert x\right\vert ^{p}}{\left\vert t\right\vert
}\right)  ^{1/(p-2)}\leqq u\leqq C_{2}\left(  \frac{\left\vert x\right\vert
^{p}}{\left\vert t\right\vert }\right)  ^{1/(p-2)}%
\]
There exist \textbf{positive} solutions, such that $r^{-\gamma}w$ is
asymptotically periodic in $\ln r$ near $0$ and in $\left\vert x\right\vert
^{\gamma}$ near $\infty$ for fixed $t;$ and also, solutions with a hole, and
oscillating around 0 for large $\left\vert x\right\vert $. There exist
solutions positive near $0$, oscillating near $\infty,$ and $r^{-\gamma}w$ is
\textbf{doubly asymptotically periodic} in $\ln r.\medskip$

4) Let $\alpha_{c}<\alpha<-p^{\prime}.\medskip$

\noindent There exist nonnegative solutions with a focussing hole near 0, in
$\left\vert x\right\vert ^{\gamma}$ near $\infty$ for fixed $t,$ blowing up at
$t=0$ like $\left(  -t\right)  ^{-1/(p-2)}$ for fixed $x\neq0.$ Either the
regular solutions have an \textbf{infinity} of zeros for fixed $t$, then the
same is true for all the other solutions. Or they have a \textbf{finite}
number $m\geqq2$ of zeros, and can be in $\left\vert x\right\vert ^{\gamma}$
or $\left\vert x\right\vert ^{\left\vert \alpha\right\vert }$ near $\infty$
(in that case they have a final data $L\left\vert x\right\vert ^{\left\vert
\alpha\right\vert }$); all the other solutions have $m$ or $m+1$ zeros.
\end{theorem}

In the case $\alpha=\alpha_{c},$ we find again the existence and uniqueness of
the focussing solutions introduced in \cite{GiVa}.

\section{Different formulations of the problem\label{form}}

In all the sequel we assume
\[
\alpha\neq0,
\]
recalling that the solutions $w$ are given explicitely by (\ref{azer}) when
$\alpha=0.$ Defining
\begin{equation}
J_{N}(r)=r^{N}\left(  w+\varepsilon r^{-1}\left\vert w^{\prime}\right\vert
^{p-2}w^{\prime}\right)  ,\qquad J_{\alpha}(r)=r^{\alpha-N}J_{N}(r),
\label{gg}%
\end{equation}
equation (\textbf{E}$_{w}$) can be written in an equivalent way under the
forms
\begin{equation}
J_{N}^{\prime}(r)=r^{N-1}(N-\alpha)w,\qquad J_{\alpha}^{\prime}%
(r)=-\varepsilon(N-\alpha)r^{\alpha-2}\left\vert w^{\prime}\right\vert
^{p-2}w^{\prime}. \label{jpn}%
\end{equation}
If $\alpha=N,$ then $J_{N}$ is constant, so we find again (\ref{fat}).
\medskip

\noindent We mainly use logarithmic substitutions; given $d\in\mathbb{R},$
setting
\begin{equation}
w(r)=r^{-d}y_{d}(\tau),\qquad Y_{d}=-r^{(d+1)(p-1)}\left\vert w^{\prime
}\right\vert ^{p-2}w^{\prime},\qquad\tau=\ln r, \label{cge}%
\end{equation}
we obtain the equivalent system:%
\begin{equation}
\left.
\begin{array}
[c]{c}%
y_{d}^{\prime}=dy_{d}-\left\vert Y_{d}\right\vert ^{(2-p)/(p-1)}Y_{d}%
,\qquad\qquad\qquad\qquad\qquad\qquad\qquad\\
\\
Y_{d}^{\prime}=(p-1)(d-\eta)Y_{d}+\varepsilon e^{(p+(p-2)d)\tau}(\alpha
y_{d}-\left\vert Y_{d}\right\vert ^{(2-p)/(p-1)}Y_{d}).
\end{array}
\right\}  \label{sysd}%
\end{equation}
At any point $\tau$ where $w^{\prime}(\tau)\neq0,$ the functions $y_{d},Y_{d}$
satisfy the equations
\begin{equation}
y_{d}^{\prime\prime}+(\eta-2d)y_{d}^{\prime}-d(\eta-d)y_{d}+\frac{\varepsilon
}{p-1}e^{((p-2)d+p)\tau}\left\vert dy_{d}-y_{d}^{\prime}\right\vert
^{2-p}\left(  y_{d}^{\prime}+(\alpha-d)y_{d}\right)  =0, \label{phi}%
\end{equation}%
\begin{equation}%
\begin{array}
[c]{c}%
Y_{d}^{\prime\prime}+(p-1)(\eta-2d-p^{\prime})Y_{d}^{\prime}+\varepsilon
e^{((p-2)d+p)\tau}\left\vert Y_{d}\right\vert ^{(2-p)/(p-1)}\left(
Y_{d}^{\prime}/(p-1)+(\alpha-d)Y_{d}\right)  \qquad\qquad\qquad\qquad\\
\qquad\qquad\qquad\qquad\qquad\qquad\qquad\qquad\qquad-(p-1)^{2}%
(\eta-d)(p^{\prime}+d)Y_{d}=0,
\end{array}
\label{phu}%
\end{equation}
The main case is $d=-\gamma$: setting $y=y_{-\gamma},$
\begin{equation}
w(r)=r^{\gamma}y(\tau),\qquad Y=-r^{(-\gamma+1)(p-1)}\left\vert w^{\prime
}\right\vert ^{p-2}w^{\prime},\qquad\tau=\ln r, \label{cha}%
\end{equation}
we are lead to the \textit{autonomous }system
\[
\text{\qquad\qquad\qquad\qquad\qquad}\left.
\begin{array}
[c]{c}%
y^{\prime}=-\gamma y-\left\vert Y\right\vert ^{(2-p)/(p-1)}Y,\qquad
\qquad\qquad\quad\\
\\
Y^{\prime}=-(\gamma+N)Y+\varepsilon(\alpha y-\left\vert Y\right\vert
^{(2-p)/(p-1)}Y).
\end{array}
\right\}  \text{\qquad\qquad\qquad\qquad\qquad(\textbf{S})}%
\]
Its study is fundamental: its phase portrait allows to study all the
\textit{signed} solutions of equation (\textbf{E}$_{w})$. Equation (\ref{phi})
takes the form%
\[
\text{\qquad\qquad\quad}(p-1)y^{\prime\prime}+(N+\gamma p)y^{\prime}%
+\gamma(\gamma+N)y+\varepsilon\left\vert \gamma y+y^{\prime}\right\vert
^{2-p}\left(  y^{\prime}+(\alpha+\gamma)y\right)  =0,\text{\qquad\qquad
\quad\textbf{(E}}_{y}\text{) }%
\]
Notice that $J_{N}(r)=r^{N+\gamma}(y(\tau)-\varepsilon Y(\tau)).\medskip$

\begin{remark}
\label{scaling}Since (\textbf{S}) is autonomous, for any solution $w$ of
(\textbf{E}$_{w}$) of the problem, all the functions $w_{\xi}(r)=\xi^{-\gamma
}w(\xi r),\xi>0,$ are also solutions.
\end{remark}

\begin{notation}
In the sequel we set $\varepsilon\infty:=+\infty$ if $\varepsilon=1,$
$\varepsilon\infty:=-\infty$ if $\varepsilon=-1.$
\end{notation}

\subsection{The phase plane of system (\textbf{S})\label{stat}}

In the phase plane $(y,Y)$ we denote the four quadrants by
\[
\mathcal{Q}_{1}=\left(  0,\infty\right)  \times\left(  0,\infty\right)
,\quad\mathcal{Q}_{2}=\left(  -\infty,0\right)  \times\left(  0,\infty\right)
,\quad\mathcal{Q}_{3}=-\mathcal{Q}_{1},\quad\mathcal{Q}_{4}=-\mathcal{Q}_{2}.
\]

\begin{remark}
\label{vf} The vector field at any point $\left(  0,\xi\right)  ,\xi>0$
satisfies $y^{\prime}=-\xi^{1/(p-1)}<0,$ thus points to $\mathcal{Q}_{2};$
moreover $Y^{\prime}<0$ if $\varepsilon=1$. The field at any point $\left(
\varphi,0\right)  ,\varphi>0$ satisfies $Y^{\prime}=\varepsilon\alpha\varphi,$
thus points to $\mathcal{Q}_{1}$ if $\varepsilon\alpha>0$ and to
$\mathcal{Q}_{4}$ if $\varepsilon\alpha<0;$ moreover $y^{\prime}%
=-\gamma\varphi<0.$
\end{remark}

If $\varepsilon(\gamma+\alpha)\geqq0,$ system (\textbf{S}) has a unique
stationary point $(0,0).$ If $\varepsilon(\gamma+\alpha)<0,$ it admits three
stationary points:
\begin{equation}
(0,0),\qquad M_{\ell}=(\ell,-(\gamma\ell)^{p-1})\in\mathcal{Q}_{4},\qquad
M_{\ell}^{\prime}=-M_{\ell}\in\mathcal{Q}_{2}, \label{statio}%
\end{equation}
where $\ell$ is defined at (\ref{defl}). The point $(0,0)$ is singular because
$p>2;$ its study concern in particular the solutions $w$ with a \textit{double
zero. }When $\varepsilon(\gamma+\alpha)<0,$ the point $M_{\ell}$ is associated
to the solution $w\equiv\ell r^{\gamma}$ of equation (\textbf{E}$_{w}$) given
at (\ref{del}). \medskip\ 

\noindent\textbf{Linearization around }$M_{\ell}.$ Near the point \textbf{
}$M_{\ell},$ setting
\begin{equation}
y=\ell+\overline{y},\qquad Y=-(\gamma\ell)^{p-1}+\overline{Y}, \label{tran}%
\end{equation}
system (\textbf{S}) is equivalent in $\mathcal{Q}_{4}$ to
\begin{equation}
\overline{y}^{\prime}=-\gamma\overline{y}-\varepsilon\nu(\alpha)\overline
{Y}+\Psi(\overline{Y}),\qquad\overline{Y}^{\prime}=\varepsilon\alpha
\overline{y}-(\gamma+N+\nu(\alpha))\overline{Y}+\varepsilon\Psi(\overline{Y}),
\label{dec}%
\end{equation}
where
\begin{equation}
\nu(\alpha)=-\frac{\gamma(N+\gamma)}{(p-1)(\gamma+\alpha)},\text{ and }%
\Psi(\vartheta)=((\gamma\ell)^{p-1}-\vartheta)^{1/(p-1)}-\gamma\ell
+\frac{(\gamma\ell)^{2-p}}{p-1}\vartheta,\qquad\vartheta<(\gamma\ell)^{p-1},
\label{nu}%
\end{equation}
thus $\varepsilon\nu(\alpha)>0.$ The linearized problem is given by%
\[
\overline{y}^{\prime}=-\gamma\overline{y}-\varepsilon\nu(\alpha)\overline
{Y},\qquad\overline{Y}^{\prime}=\varepsilon\alpha\overline{y}-(\gamma
+N+\nu(\alpha))\overline{Y}.
\]
Its eigenvalues $\lambda_{1}\leqq\lambda_{2}$ are the solutions of equation%
\begin{equation}
\lambda^{2}+(2\gamma+N+\nu(\alpha))\lambda+p^{\prime}(N+\gamma)=0 \label{egv}%
\end{equation}
The discriminant $\Delta$ of the equation (\ref{egv}) is given by
\begin{equation}
\Delta=(2\gamma+N+\nu(\alpha))^{2}-4p^{\prime}(N+\gamma)=(N+\nu(\alpha
))^{2}-4\nu(\alpha)\alpha. \label{lta}%
\end{equation}
$\medskip$For $\varepsilon=1,$ $M_{\ell}$ is a \textit{sink}, and a node
point, since $\nu(\alpha)>0,$ and $\alpha<0,$ thus $\Delta>0$. For
$\varepsilon=-1,$ we have $\nu(\alpha)<0$; the nature of $M_{\ell}$ depends on
the critical value $\alpha^{\ast}$ defined at (\ref{star}); indeed
\[
\alpha=\alpha^{\ast}\Longleftrightarrow\lambda_{1}+\lambda_{2}=0.
\]
Then $M_{\ell}$ is a \textit{sink} when $\alpha>\alpha^{\ast}$ and a
\textit{source} when $\alpha<\alpha^{\ast}.$ Moreover $\alpha^{\ast}$
corresponds to a spiral point, and $M_{\ell}$ is a node point when
$\Delta\geqq0,$ that means $\alpha\leqq\alpha_{1},$ or $\gamma>N/2+\sqrt
{p^{\prime}(N+\gamma)}$ and $\alpha_{2}\leqq\alpha,$ where
\begin{equation}
\alpha_{1}=-\gamma+\frac{\gamma(N+\gamma)}{(p-1)(2\gamma+N+2(p^{\prime
}(N+\gamma))^{1/2})},\qquad\alpha_{2}=-\gamma+\frac{\gamma(N+\gamma
)}{(p-1)(2\gamma+N-2(p^{\prime}(N+\gamma))^{1/2})}. \label{aun}%
\end{equation}
When $\Delta>0$, and $\lambda_{1}<\lambda_{2},$ one can choose a basis of
eigenvectors
\begin{equation}
e_{1}=(-\varepsilon\nu(\alpha),\lambda_{1}+\gamma)\quad\text{and\quad}%
e_{2}=(\varepsilon\nu(\alpha),-\gamma-\lambda_{2}). \label{bas}%
\end{equation}

\begin{remark}
\label{sign} One verifies that $\alpha^{\ast}<-1;$ and $\alpha^{\ast
}<-(p-1)/(p-2)$ if and only if $p>N.$ Also $\alpha_{2}\leqq0,$ and $\alpha
_{2}=0\Longleftrightarrow N=p/((p-2)^{2};$ and $\alpha_{2}$ $>-p^{\prime}$
$\Longleftrightarrow$ $\gamma^{2}-7\gamma-8N<0,$ which is not always
true.\medskip
\end{remark}

As in \cite[Theorem 2.16]{Bi1} we prove that the Hopf bifurcation point is not
degenerate, which implies the existence of small cycles near $\alpha^{\ast}.$

\begin{proposition}
\label{ws}Let $\varepsilon=-1,$ and $\alpha=\alpha^{\ast}>-\gamma.$ Then
$M_{\ell}$ is a weak source. If $\alpha>\alpha^{\ast}$ and $\alpha
-\alpha^{\ast}$ is small enough, there exists a unique limit cycle in
$\mathcal{Q}_{4},$ attracting at $-\infty.$
\end{proposition}

\subsection{Other systems for positive solutions}

When $w$ has a constant sign, we define two functions associated to $(y,Y):$
\begin{equation}
\zeta(\tau)=\frac{\left\vert Y\right\vert ^{(2-p)/(p-1)}Y}{y}(\tau
)=-\frac{rw^{\prime}(r)}{w(r)},\qquad\sigma(\tau)=\frac{Y}{y}(\tau
)=-\frac{\left\vert w^{\prime}(r)\right\vert ^{p-2}w^{\prime}(r)}{rw(r)}.
\label{dze}%
\end{equation}
Thus $\zeta$ describes the behaviour of $w^{\prime}/w$ and $\sigma$ is the
slope in the phase plane $\left(  y,Y\right)  .$ They satisfy the system%
\[
\qquad\left.
\begin{array}
[c]{c}%
\zeta^{\prime}=\zeta(\zeta-\eta)+\varepsilon\left\vert \zeta y\right\vert
^{2-p}(\alpha-\zeta)/(p-1)=\zeta(\zeta-\eta+\varepsilon(\alpha-\zeta
)/(p-1)\sigma),\\
\\
\sigma^{\prime}=\varepsilon(\alpha-N)+\left(  \left\vert \sigma y\right\vert
^{(2-p)/(p-1)}\sigma-N\right)  (\sigma-\varepsilon)=\varepsilon(\alpha
-\zeta)+\left(  \zeta-N\right)  \sigma.
\end{array}
\right\}  \text{\qquad\qquad\quad(\textbf{Q})}%
\]
In particular, System (\textbf{Q}) provides a short proof of the local
existence and uniqueness of the \textit{regular} solutions: they correspond to
its stationary point $(0,\varepsilon\alpha/N)$, see Section \ref{regu}%
.\medskip\ 

Moreover, if $w$ and $w^{\prime}$ have a strict constant sign, that means in
any quadrant $\mathcal{Q}_{i},$ we can define
\begin{equation}
\psi=\frac{1}{\sigma}=\frac{y}{Y} \label{dzp}%
\end{equation}
We obtain a new system relative to $(\zeta,\psi):$%
\[
\qquad\qquad\qquad\qquad\qquad\left.
\begin{array}
[c]{c}%
\zeta^{\prime}=\zeta(\zeta-\eta+\varepsilon(\alpha-\zeta)\psi/(p-1)),\\
\\
\psi^{\prime}=\psi(N-\zeta+\varepsilon\left(  \zeta-\alpha)\psi\right)
.\qquad\,\quad
\end{array}
\right\}  \text{\qquad\qquad}\qquad\qquad\qquad\text{\qquad\quad\textbf{(P})}%
\]
We are reduced to a polynomial system, thus with no singularity. System
(\textbf{P}) gives the existence of singular solutions when $p>N,$
corresponding to its stationary point $(\eta,0),$ see Section \ref{new}%
.\medskip

We will also consider another system in any $\mathcal{Q}_{i}:$ setting
\begin{equation}
\zeta=-1/g,\qquad\sigma=-s,\qquad d\tau=gsd\nu=\left\vert Y\right\vert
^{(p-2)/(p-1)}d\nu, \label{gsnu}%
\end{equation}
we find
\[
\qquad\qquad\qquad\qquad\left.
\begin{array}
[c]{c}%
dg/d\nu=g(s(1+\eta g)+\varepsilon(1+\alpha g)/(p-1)),\\
\\
ds/d\nu=-s(\varepsilon(1+\alpha g)+(1+Ng)s).\quad\quad\quad
\end{array}
\right\}  \qquad\qquad\text{\qquad\qquad}\qquad\quad\text{\quad\textbf{(R})}%
\]
System (\textbf{R}) allows to get the existence of solutions $w$ with a hole
or a compact support, and other solutions, corresponding to its stationary
points $(0,-\varepsilon)$ and ($-1/\alpha,0)$; it provides a complete study of
the singular point $(0,0)$ of system (\textbf{S}), see Sections \ref{doub},
\ref{new}; and of the focussing solutions, see Section \ref{four}.

\begin{remark}
\label{parti}The particular solutions can be found again in the different
phase planes, where their trajectories are lines:

For $\alpha=N,$ the solutions (\ref{fut}) correspond to $Y\equiv\varepsilon
y,$ that means $\sigma\equiv\varepsilon.$

For $\alpha=\eta\neq0$ the solutions (\ref{aeta}) correspond to $\zeta
\equiv\eta.$

For $\alpha=-p^{\prime},$ the solutions (\ref{moi}) are given by
$\zeta+\varepsilon N\sigma\equiv\alpha.$

For $N=1,$ $\alpha=-(p-2)/(p-1),$ the solutions (\ref{exp}) satisfy $\alpha
g+\varepsilon s\equiv-1.$
\end{remark}

\section{Global existence\label{glob}}

\subsection{Local existence and uniqueness\label{regu}}

\begin{proposition}
\label{evi}Let $r_{1}>0$ and $a,b\in\mathbb{R}.$ If $(a,b)\neq(0,0)$, there
exists a unique solution $w$ of equation (\textbf{E}$_{w}$) in a neighborhood
$\mathcal{V}$ of $r_{1}$, such that $w$ and $\left\vert w^{\prime}\right\vert
^{p-2}w^{\prime}\in C^{1}\left(  \mathcal{V}\right)  $ and $w(r_{1})=a,$
$w^{\prime}(r_{1})=b.$ It extends on a maximal interval $I$ where
$(w(r),w^{\prime}(r))\neq(0,0).$
\end{proposition}

\begin{proof}
If $b\neq0,$ the Cauchy theorem directly applies to system (\textbf{S}). If
$b=0$ the system is a priori singular on the line $\left\{  Y=0\right\}  $
since $p>2.$ In fact it is only singular at $(0,0).$ Indeed near any point
$\left(  \xi,0\right)  $ with $\xi\neq0,$ one can take $Y$ as a variable, and
\[
\frac{dy}{dY}=F(Y,y),\qquad F(Y,y):=\frac{\gamma y+\left\vert Y\right\vert
^{(2-p)/(p-1)}Y}{(\gamma+N)Y+\varepsilon(\left\vert Y\right\vert
^{(2-p)/(p-1)}Y-\alpha y)},
\]
where $F$ is continuous in $Y$ and $C^{1}$ in $y,$ hence local existence and
uniqueness hold.
\end{proof}

\begin{notation}
For any point $P_{0}=(y_{0},Y_{0})\in\mathbb{R}^{2}\backslash\left\{
(0,0)\right\}  ,$ the unique trajectory in the phase plane $(y,Y)$ of system
(\textbf{S}) going through $P_{0}$ is denoted by $\mathcal{T}_{\left[
P_{0}\right]  }.$ By symmetry, $\mathcal{T}_{\left[  -P_{0}\right]
}=-\mathcal{T}_{\left[  P_{0}\right]  }$.
\end{notation}

Next we show the existence of regular solutions. Our proof is short, based on
phase plane portrait, and not on a fixed point method, rather delicate because
$p>2,$ see \cite{Bi0}.

\begin{theorem}
\label{exlo} For any $a\in\mathbb{R},$ $a\neq0,$ there exists a unique
solution $w=w(.,a)$ of equation (\textbf{E}$_{w}$) in an interval $\left[
0,r_{0}\right)  ,$ such that $w$ and $\left\vert w^{\prime}\right\vert
^{p-2}w^{\prime}\in C^{1}\left(  \left[  0,r_{0}\right)  \right)  $ and
\begin{equation}
w(0)=a,\qquad w^{\prime}(0)=0; \label{ini}%
\end{equation}
and then $\lim_{r\rightarrow0}\left\vert w^{\prime}\right\vert ^{p-2}%
w^{\prime}/rw=-\varepsilon\alpha/N.$ In other words in the phase plane $(y,Y)$
there exists a unique trajectory $\mathcal{T}_{r}$ such that $\lim
_{\tau\rightarrow-\infty}y=\infty,$ and $\lim_{\tau\rightarrow-\infty
}Y/y=\varepsilon\alpha/N.$
\end{theorem}

\begin{proof}
We have assumed $\alpha\neq0$ (when $\alpha=0,w\equiv a$ from (\ref{azer})).
If such a solution $w$ exists, then from (\ref{gg}) and (\ref{jpn}),
$J_{N}^{\prime}(r)=r^{N-1}(N-\alpha)a(1+o(1))$ near $0.$ Thus $J_{N}%
(r)=r^{N-1}(1-\alpha/N)a(1+o(1)),$ hence $\lim_{r\rightarrow0}\left\vert
w^{\prime}\right\vert ^{p-2}w^{\prime}/rw=-\varepsilon\alpha/N;$ in other
words, $\lim_{\tau\rightarrow-\infty}\sigma=\varepsilon\alpha/N.$ And
$\lim_{\tau\rightarrow-\infty}y=\infty,$ thus $\lim_{\tau\rightarrow-\infty
}\zeta=0,$ and $\varepsilon\alpha\zeta>0$ near $-\infty.$ Reciprocally
consider system (\textbf{Q}). The point $(0,\varepsilon\alpha/N)$ is
stationary. Setting $\sigma=\varepsilon\alpha/N+\bar{\sigma},$ the linearized
system near this point is given by
\[
\zeta^{\prime}=p^{\prime}\zeta,\qquad\bar{\sigma}^{\prime}=\varepsilon
\zeta(\alpha-N)/N-N\bar{\sigma}.
\]
One finds is a saddle point, with eigenvalues $-N$ and $p^{\prime}.$ Then
there exists a unique trajectory $\mathcal{T}_{r}^{\prime}$ in the phase-plane
$\left(  \zeta,\sigma\right)  $ starting at $-\infty$ from $(0,\varepsilon
\alpha/N)$ with the slope $\varepsilon(\alpha-N)/N(N+p^{\prime})\neq0$ and
$\varepsilon\alpha\zeta>0.$ It corresponds to a unique trajectory
$\mathcal{T}_{r}$ in the phase plane $(y,Y),$ and $\lim_{\tau\rightarrow
-\infty}y=\infty,$ since $y=\left\vert \sigma\right\vert \left\vert
\zeta\right\vert ^{1-p})^{1/(p-2)}$.\ For any solution $\left(  \zeta
,\sigma\right)  $ describing $\mathcal{T}_{r}^{\prime}$, the function
$w(r)=r^{\gamma}(\left\vert \sigma\right\vert \left\vert \zeta\right\vert
^{1-p}(\tau))^{1/(p-2)}$ satisfies $\lim_{r\rightarrow0}\left\vert w^{\prime
}\right\vert ^{p-2}w^{\prime}/rw=-\varepsilon\alpha/N.$ As a consequence,
$w^{(p-2)/(p-1)}$ has a finite nonzero limit, and $\lim_{r\rightarrow
0}w^{\prime}=0;$ thus $w$ is regular. Local existence and uniqueness follows
for any $a\neq0,$ by Remark \ref{scaling}.
\end{proof}

\begin{definition}
The trajectory $\mathcal{T}_{r}$ in the plane $(y,Y)$ and its opposite
$-\mathcal{T}_{r}$ will be called \textit{regular trajectories}. We shall say
that $y$ is regular. Observe that $\mathcal{T}_{r}$ starts in $\mathcal{Q}%
_{1}$ if $\varepsilon\alpha>0,$ and in $\mathcal{Q}_{4}$ if $\varepsilon
\alpha<0.$
\end{definition}

\begin{remark}
From Theorem \ref{exlo} and Remark \ref{scaling}, all regular solutions are
obtained from one one of them: $w(r,a)=aw(a^{-1/\gamma}r,1).$ Thus they have
the same behaviour near $\infty.$
\end{remark}

\subsection{Sign properties}

Next we give informations on the zeros of $w$ or $w^{\prime},$ by using the
monotonicity properties of the functions $y_{d},Y_{d},$ in particular $y,Y,$
and $\zeta$ and $\sigma$. At any extremal point $\tau$, they satisfy
respectively%
\begin{equation}
y_{d}^{\prime\prime}(\tau)=y_{d}(\tau)\left(  d(\eta-d)+\frac{\varepsilon
(d-\alpha)}{p-1}e^{((p-2)d+p)\tau}\left\vert dy_{d}(\tau)\right\vert
^{2-p}\right)  , \label{yde}%
\end{equation}%
\begin{equation}
Y_{d}^{\prime\prime}(\tau)=Y_{d}(\tau)\left(  (p-1)^{2}(\eta-d)(p^{\prime
}+d)+\varepsilon(d-\alpha)e^{((p-2)d+p)\tau}\left\vert Y_{d}(\tau)\right\vert
^{(2-p)/(p-1)}\right)  , \label{Yd}%
\end{equation}%
\begin{equation}
(p-1)y^{\prime\prime}(\tau)=\gamma^{2-p}y(\tau)\left(  -\gamma^{p-1}%
(N+\gamma)-\varepsilon(\gamma+\alpha)\left\vert y(\tau)\right\vert
^{2-p}\right)  =-\left\vert Y(\tau)\right\vert ^{(2-p)/(p-1)}Y^{\prime}(\tau),
\label{sqr}%
\end{equation}%
\begin{equation}
Y^{\prime\prime}(\tau)=Y(\tau)\left(  -\gamma(N+\gamma)-\varepsilon
(\gamma+\alpha)\left\vert Y(\tau)\right\vert ^{(2-p)/(p-1)}\right)
=\varepsilon\alpha y^{\prime}(\tau), \label{spr}%
\end{equation}%
\begin{equation}
(p-1)\zeta^{\prime\prime}(\tau)=-\varepsilon(p-2)((\alpha-\zeta)\left\vert
\zeta\right\vert ^{2-p}\left\vert y\right\vert ^{-p}yy^{\prime})(\tau
)=\varepsilon(p-2)((\alpha-\zeta)(\gamma+\zeta)\left\vert \zeta y\right\vert
^{2-p})(\tau), \label{kse}%
\end{equation}%
\begin{equation}
(p-1)\sigma^{\prime\prime}(\tau)=-(p-2)((\sigma-\varepsilon)\left\vert
\sigma\right\vert ^{(2-p)/(p-1)}Y\left\vert y\right\vert ^{(4-3p)/(p-1)}%
y^{\prime})(\tau)=\zeta^{\prime}(\tau)(\sigma(\tau)-\varepsilon). \label{segh}%
\end{equation}

\begin{proposition}
\label{zer} Let $w\not \equiv 0$ be any solution of (\textbf{E}$_{w}$) on an
interval $I$.\medskip

\noindent(i) If $\varepsilon=1$ and $\alpha\leqq N,$ then $w$ has at most one
simple zero; if $\alpha<N$ and $w$ is regular, it has no zero. If $\alpha=N$
it has no simple zero and a compact support. If $\alpha>N$ and $w$ is regular,
it has at least one simple zero. \medskip

\noindent(ii) If $\varepsilon=-1$ and $\alpha\geqq\min(0,\eta),$ then $w$ has
at most one simple zero. If $w\not \equiv 0$ has a double zero, then it has no
simple zero. If $\alpha>0$ and $w$ is regular, it has no zero.\medskip

\noindent(iii) If $\varepsilon=-1$ and $-p^{\prime}\leqq\alpha<\min(0,\eta),$
then $w^{\prime}$ has at most one simple zero, consequently $w$ has at most
two simple zeros, and at most one if $w$ is regular. If $\alpha<-p^{\prime},$
the regular solutions have at least two zeros.
\end{proposition}

\begin{proof}
(i) Let $\varepsilon=1.$ Consider two consecutive simple zeros $\rho_{0}<$
$\rho_{1}$ of $w,$ with $w>0$ on $\left(  \rho_{0},\rho_{1}\right)  ;$ hence
$w^{\prime}(\rho_{1})<0<w^{\prime}(\rho_{0})$. If\textbf{ } $\alpha\leqq N,$
we find from (\ref{gg}),
\[
J_{N}(\rho_{1})-J_{N}(\rho_{0})=-\rho_{1}^{N-1}\left\vert w^{\prime}(\rho
_{1})\right\vert ^{p-2}-\rho_{0}^{N-1}w^{\prime}(\rho_{0})^{p-1}=(N-\alpha)%
{\displaystyle\int\nolimits_{\rho_{0}}^{\rho_{1}}}
s^{N-1}wds,
\]
which is contradictory; thus $w$ has at most one simple zero. The
contradiction holds as soon as $\rho_{0}$ is simple, even if $\rho_{1}$ is
not. If $w$ is regular with $w(0)>0,$ and $\rho_{1}$ is a first zero, and
$\alpha<N,$
\[
J_{N}(\rho_{1})=-\rho_{1}^{N-1}\left\vert w^{\prime}(\rho_{1})\right\vert
^{p-1}=(N-\alpha)%
{\displaystyle\int\nolimits_{0}^{\rho_{1}}}
s^{N-1}wds>0,
\]
which is still impossible. If $\alpha=N,$ the (Barenblatt) solutions are given
by (\ref{fut})$.$ Next suppose $\alpha>N$ and $w$ regular. If $w>0,$ then
$J_{N}<0,$ thus $w^{-1/(p-1)}w^{\prime}+r^{1/(p-1)}<0.$ Then the function
$r\mapsto r^{p^{\prime}}+\gamma w^{(p-2)/(p-1)}$ is non increasing and we
reach a contradiction for large $r.$ Thus $w$ has a first zero $\rho_{1},$ and
$J_{N}(\rho_{1})<0,$ thus $w^{\prime}(\rho_{1})\neq0.$ \medskip

\noindent(ii) Let $\varepsilon=-1$ and $\alpha\geqq\min(\eta,0).$ Here we use
the substitution (\ref{cge}) from some $d\neq0.$ If $y_{d}$ has a maximal
point, where it is positive, and is not constant, then (\ref{yde}) holds.
Taking $d\in\left(  0,\min(\alpha,\eta\right)  )$ if $\eta>0$, $d=\eta$ if
$\eta\leqq0,$ we reach a contradiction. Hence $y_{d}$ has at most a simple
zero, and no simple zero if it has a double one.\ Suppose $w$ regular and
$\alpha>0.$ Then $w^{\prime}>0$ near $0,$ from Theorem \ref{exlo}. As long as
$w$ stays positive, any extremal point $r$ is a strict minimum, from
(\textbf{E}$_{w}$), thus in fact $w^{\prime}$ stays positive. \medskip

\noindent(iii) Let $\varepsilon=-1$ and $-p^{\prime}\leqq\alpha<\min(0,\eta).$
Suppose that $w^{\prime}$ and has two consecutive zeros $\rho_{0}<$ $\rho_{1}%
$, and one of them is simple, and use again (\ref{cge}) with $d=\alpha$. Then
the function $Y_{\alpha}$ has an extremal point $\tau$, where it is positive
and is not constant; from (\ref{Yd}),
\begin{equation}
Y_{\alpha}^{\prime\prime}(\tau)=(p-1)^{2}(\eta-\alpha)(p^{\prime}%
+\alpha)Y_{\alpha}(\tau), \label{kal}%
\end{equation}
thus $Y_{\alpha}^{\prime\prime}(\tau)\geqq0,$ which is contradictory. Next
consider the regular solutions. They satisfy $Y_{\alpha}(\tau)=e^{(\alpha
(p-1)+p)\tau}(\left\vert \alpha\right\vert a/N)(1+o(1)$ near $-\infty,$ from
Theorem \ref{exlo} and (\ref{cge}), thus $\lim_{\tau\rightarrow-\infty
}Y_{\alpha}=0.$ As above $Y_{\alpha}$ cannot have any extremal point, then
$Y_{\alpha}$ is positive and increasing. In turn $w^{\prime}<0$ from
(\ref{cge}), hence $w$ has at most one zero. \medskip
\end{proof}

\begin{proposition}
\label{fini}Let $w\not \equiv 0$ be any solution of (\textbf{E}$_{w}$) on an
interval $I$. If\ $\varepsilon=1,$ then $w$ has a finite number of isolated
zeros. If $\varepsilon=-1,$ it has a finite number of isolated zeros in any
interval $\left[  m,M\right]  \cap I$ with $0<m<M<\infty.$
\end{proposition}

\begin{proof}
Let $Z$ be the set of isolated zeros on $I$. If $w$ has two consecutive
isolated zeros $\rho_{1}<$ $\rho_{2}$, and $\tau\in\left(  e^{\rho_{1}%
},e^{\rho_{2}}\right)  $ is a maximal point of $\left\vert y_{d}\right\vert $,
from (\ref{yde}), it follows that%
\begin{equation}
\varepsilon e^{((p-2)d+p)\tau}\left\vert dy_{d}(\tau)\right\vert ^{2-p}\left(
d-\alpha\right)  \leqq(p-1)d(d-\eta). \label{tig}%
\end{equation}
That means with $\rho=e^{\tau}\in\left(  \rho_{1},\rho_{2}\right)  ,$
\begin{equation}
\varepsilon\rho^{p}\left\vert w(\rho)\right\vert ^{2-p}\left(  d-\alpha
\right)  \leqq(p-1)d^{p-1}(d-\eta). \label{tog}%
\end{equation}
First suppose $\varepsilon=1$ and fix $d>\alpha.$ Consider the energy
function
\[
E(r)=\frac{1}{p^{\prime}}\left\vert w^{\prime}\right\vert ^{p}+\frac{\alpha
}{2}w^{2}.
\]
It is nonincreasing since $E^{\prime}(r)=-(N-1)r^{-1}\left\vert w^{\prime
}\right\vert ^{p}-rw^{\prime2},$ thus bounded on $I\cap\left[  \rho_{1}%
,\infty\right)  .$ Then $w$ is bounded, $\rho_{2}$ is bounded$,$ $Z$ is a
bounded set. If $Z$ is infinite, there exists a sequence of zeros $\left(
r_{n}\right)  $ converging to some point $\overline{r}\in\left[
0,\infty\right)  ,$ and a sequence $\left(  \tau_{n}\right)  $ of maximal
points of $\left\vert y_{d}\right\vert $ converging to $\overline{\tau}%
=\ln\overline{r}$. If $\overline{r}>0,$ then $w(\overline{r})=w^{\prime
}(\overline{r})=0;$ we get a contradiction by taking $\rho=\rho_{n}%
=e^{\tau_{n}}$ in (\ref{tog}), because the left-hand side tends to $\infty.$
If $\overline{r}=0,$ fixing now $d<\eta$, there exists a sequence $\left(
\tau_{n}\right)  $ of maximal points of $\left\vert y_{d}\right\vert $
converging to $-\infty$. Then $w(\rho_{n})=O(\rho_{n}^{p/(p-2)}),$ and
$w^{\prime}(\rho_{n})=-d\rho_{n}^{-1}w(\rho_{n})=O(\rho_{n}^{2/(p-2)}),$ thus
$E(\rho_{n})=o(1).$ Since $E$ is monotone, it implies $\lim_{r\rightarrow
0}E(r)=0,$ hence $E\equiv0,$ and $w\equiv0,$ which is contradictory. Next
suppose $\varepsilon=-1$ and fix $d<\alpha.$ If $Z\cap\left[  m,M\right]  $ is
infinite, we construct a sequence converging vers some $\overline{r}>0$ and
reach a contradiction as above.\ \medskip
\end{proof}

\begin{proposition}
\label{com} Let $y$ be any non constant solution of (\textbf{E}$_{y}$), on a
maximal interval $I$ where $(y,Y)\neq(0,0),$ and $s$ be an extremity of
$I.\medskip$

\noindent(i) If $y$ has a constant sign near $s,$ then the same is true for
$Y$.$\medskip$

\noindent(ii) If $y>0$ is strictly monotone near $s$, then $Y,\zeta,\sigma$
are monotone near $s.\medskip$

\noindent(iii) If $y>0$ is not strictly monotone$\ $near $s,$ then
$s=\pm\infty,$ $\varepsilon(\gamma+\alpha)<0$ and $y$ oscillates around
$\ell.\medskip$

\noindent(iv) If $y$ is oscillating around 0 near $s$, then $\varepsilon
=-1,s=\pm\infty,\alpha<-p^{\prime};$ if $\alpha>-\gamma,$ then $\left\vert
y\right\vert >\ell$ at the extremal points.
\end{proposition}

\begin{proof}
(i) The function $w$ has at most one extremal point on $I:$ at such a point,
it satisfies $(\left\vert w^{\prime}\right\vert ^{p-2}w^{\prime})^{\prime
}=-\varepsilon\alpha w$ with $\alpha\neq0.$ From (\ref{cha}), $Y$ has a
constant sign near $s.\medskip$

\noindent(ii)\ Suppose $y$ strictly monotone near $s.$ At any extremal point
$\tau$ of $Y,$ we find $Y^{\prime\prime}(\tau)=\varepsilon\alpha y^{\prime
}(\tau)$ from (\ref{spr}). Then $y^{\prime}(\tau)\neq0,$ $Y^{\prime\prime
}(\tau)$ has a constant sign. Thus $\tau$ is unique, and $Y$ is strictly
monotone near $s$. Next consider $\zeta$. If there exists $\tau_{0}$ such that
$\zeta(\tau_{0})=\alpha,$ then $\zeta^{\prime}(\tau_{0})=\alpha(\alpha-\eta),$
from system \textbf{(Q)}. If $\alpha\neq\eta,$ then $\tau_{0}$ is unique, thus
$\alpha-\zeta$ has a constant sign near $s.$ Then $\zeta^{\prime\prime}(\tau)$
has a constant sign at any extremal point $\tau$ of $\zeta,$ from (\ref{kse}),
thus $\zeta$ is strictly monotone near $s.$ If $\alpha=\eta,$ then
$\zeta\equiv\alpha$. At last consider $\sigma$. If there exists $\tau_{0}$
such that $\sigma(\tau_{0})=\varepsilon,$ then $\sigma^{\prime}(\tau
_{0})=\varepsilon(\alpha-N)$ from System \textbf{(Q).} If $\alpha\neq N$, then
$\tau_{0}$ is unique, and $\sigma-\varepsilon$ has a constant sign near $s.$
Thus $\sigma^{\prime\prime}(\tau)$ has a constant sign at any extremal point
$\tau$ of $\sigma,$ from (\ref{segh}) and assertion (i). If $\alpha=N$, then
$\sigma\equiv\varepsilon.$\medskip

\noindent(iii) Let $y$ be positive and not strictly monotone near $s.$ There
exists a sequence $(\tau_{n})$ strictly monotone, converging to $\pm\infty,$
such that $y^{\prime}(\tau_{n})=0,$ $y^{\prime\prime}(\tau_{2n})>0>$
$y^{\prime\prime}(\tau_{2n+1}).$ Since $y(\tau_{n})=$ $\gamma^{-1}\left\vert
Y\right\vert ^{(2-p)/(p-1)}Y(\tau_{n}),$ we deduce $Y<0$ near $s,$ from (i).
From (\ref{spr}),
\begin{equation}
-\varepsilon(\gamma+\alpha)y(\tau_{2n+1})^{2-p}\leqq\gamma^{p-1}%
(N+\gamma)\leqq-\varepsilon(\gamma+\alpha))y(\tau_{2n})^{2-p}, \label{tru}%
\end{equation}
thus $\varepsilon(\gamma+\alpha)<0$ and $y(\tau_{2n})<\ell<y(\tau_{2n+1}),$
and $Y(\tau_{2n+1})<-\left(  \gamma\ell\right)  ^{p-1}<Y(\tau_{2n}).$ If $s$
is finite, then $y(s)=y^{\prime}(s)=0,$ which is impossible; thus $s=\pm
\infty.$ \medskip

\noindent(iv) If $y$ is changing sign, then $\varepsilon=-1$ and
$\alpha<-p^{\prime},$ from Propositions \ref{zer} and \ref{fini}. At any
extremal point $\tau,$
\[
(\alpha+\gamma)\left\vert y(\tau)\right\vert ^{2-p}\leqq\gamma^{p-1}%
(N+\gamma)
\]
from (\ref{sqr}); if $\alpha>-\gamma$ it means $\left\vert y\right\vert
>\ell.$
\end{proof}

\subsection{Double zeros and global existence\label{doub}}

\begin{theorem}
\label{tpsilon}For any $\overline{r}>0,$ there exists a unique solution $w$ of
(\textbf{E}$_{w}$) defined in a interval $\left[  \overline{r},\overline{r}\pm
h\right)  $ such that
\[
w>0\quad\text{on }\left(  \overline{r},\overline{r}\pm h\right)
\quad\text{and}\quad w(\overline{r})=w^{\prime}(\overline{r})=0.
\]
Moreover $\varepsilon h<0$ and
\begin{equation}
\lim_{r\rightarrow\overline{r}}\left\vert (\overline{r}-r)\right\vert
^{(p-1)/(2-p)}\overline{r}^{1/(2-p)}w(r)=\pm((p-2)/(p-1))^{(p-1)/(p-2)}.
\label{rip}%
\end{equation}
In other words in the phase plane $(y,Y)$ there exists a unique trajectory
$\mathcal{T}_{\varepsilon}$ converging to $(0,0)$ at $\varepsilon\infty.$ It
has the slope $\varepsilon$ and converges in finite time; it depends locally
continuously of $\alpha.$
\end{theorem}

\begin{proof}
Suppose that a solution $w\not \equiv 0$ exists on $\left[  \overline
{r},\overline{r}\pm h\right)  $ with $w(\overline{r})=w^{\prime}(\overline
{r})=0.$ From Propositions \ref{fini} and \ref{com}, up to a symmetry, $y>0,$
$\left\vert Y\right\vert >0$ near $\bar{\tau}=\ln\overline{r},$ and
$\lim_{\tau\rightarrow\ln\overline{r}}y=$ $0,$ and $\sigma,\zeta$ are monotone
near $\ln r$. Let $\mu$ and $\lambda$ be their limits. If $\left\vert
\mu\right\vert =\infty,$ then $\left\vert \lambda\right\vert =\infty,$ because
$\zeta=\left\vert Y\right\vert ^{(2-p)/(p-1)}\sigma$, $\left\vert
\zeta\right\vert ^{p-2}\zeta=\sigma y^{2-p};$ then $f=1/\zeta$ tends to $0;$
but
\begin{equation}
f^{\prime}=-1+\eta f+\varepsilon\frac{1-\alpha f}{(p-1)\sigma}, \label{ff}%
\end{equation}
thus $f^{\prime}$ tends to $-1,$ which is impossible. Thus $\mu$ is finite. If
$\lambda$ is finite, then $\mu=0,$ thus $\lambda=\alpha,$ from system
(\textbf{Q}), $\ln w$ is integrable at $\overline{r}$, which is not true. Then
$\lambda=\varepsilon\infty,$ hence
\[
\mu=\lim_{\tau\rightarrow\ln\overline{r}}\sigma=\varepsilon,
\]
from system (\textbf{Q}). Then $\varepsilon Y>0$ near $\bar{\tau},$ then
$\varepsilon w^{\prime}<0$ near $\overline{r},$ thus $\varepsilon h<0.$
Consider system (\textbf{R}): as $\tau$ tends to $\bar{\tau},$ $\nu$ tends to
$\pm\infty,$ and $(g,s)$ converges to the stationary point $(0,-\varepsilon
).\medskip$

Reciprocally, setting $s=-\varepsilon/\beta+h,$ the linearized system of
system (\textbf{R}) at this point is given by
\[
\frac{dg}{d\nu}=-\varepsilon\frac{p-2}{p-1}g,\qquad\frac{dh}{d\nu}%
=(\alpha-N)g+\varepsilon h.
\]
The eigenvalues are $-\varepsilon(p-2)/(p-1)$ and $\varepsilon,$ thus we find
a saddle point. There are two trajectories converging to $(0,-\varepsilon).$
The first one satisfies $g\equiv0,$ it does not correspond to a solution of
the initial problem.\ Then there exists a unique trajectory converging to
$(0,-\varepsilon),$ as $\nu$ tends to $\varepsilon\infty,$ with $g>0$ near
$\varepsilon\infty.$ It is associated to the eigenvalue $-\varepsilon
(p-2)/(p-1)$ and the eigenvector $((2p-3)/(p-1),\varepsilon(N-\alpha)).$ It
satisfies $dg/d\nu=-\varepsilon((p-2)/(p-1))g(1+o(1)),$ thus $dg/d\tau
=((p-2)/(p-1))(1+o(1)).$ Then $\tau$ has a finite limit $\bar{\tau},$ and
$\tau$ increases to $\bar{\tau}$ if $\varepsilon=1$ and decreases to
$\bar{\tau}$ if $\varepsilon=-1.$ In turn $\left\vert Y\right\vert
^{(p-2)/(p-1)}=gs$ tends to $0,$ and $s$ tends to $\varepsilon,$ thus $(y,Y)$
tends to $(0,0)$ as $\tau$ tends to $\bar{\tau}.$ Then $w$ and $w^{\prime}$
converges to $0$ at $\overline{r}=e^{\bar{\tau}}.$ And $w^{\prime}%
w^{-1/(p-1)}+(\varepsilon+o(1))r^{1/(p-1)}=0,$ which implies (\ref{rip}).

\begin{corollary}
\label{prol}Let $r_{1}>0,$ and $a,b\in\mathbb{R}$ and $w$ be any local
solution such that $w(r_{1})=a,$ $w^{\prime}(r_{1})=b.\medskip$

\noindent(i) If $(a,b)=(0,0),$ then $w$ has a unique extension by $0$ on
$\left(  r_{1},\infty\right)  $ if $\varepsilon=1,$ on $\left(  0,r_{1}%
\right)  $ if $\varepsilon=-1.\medskip$

\noindent(ii) If $(a,b)\neq(0,0)$, $w$ has a unique extension to $\left(
0,\infty\right)  .$
\end{corollary}
\end{proof}

\begin{proof}
(i) Assume $a=b=0,$ the function $w\equiv0$ is a solution. Let $w$ be any
local solution near $r_{1}$, defined in an interval $\left(  r_{1}-h_{1}%
,r_{1}+h_{1}\right)  $ with $w(r_{1})=w^{\prime}(r_{1})=0$. Suppose that there
exists $h_{2}\in\left(  0,h_{1}\right)  $ such that $w(r_{1}+\varepsilon
h_{1})\neq0.$ Let $\bar{h}=\inf\left\{  h\in\left(  0,h_{1}\right)
:w(r_{1}+\varepsilon h)\neq0\right\}  ,$ and $\bar{r}=r_{1}+\varepsilon\bar
{h},$ thus $w(\bar{r})=w^{\prime}(\bar{r})=0$, and for example $w>0$ on some
interval $(\bar{r},\bar{r}+\varepsilon k))$ with $k>0.$ This contradicts
theorem \ref{tpsilon}. Thus $w\equiv0$ on $\left(  r_{1},r_{1}+\varepsilon
h_{1}\right)  .\medskip$

\noindent(ii) From Theorems \ref{tpsilon} and \ref{exlo}, $w$ has no double
zero for $\varepsilon\left(  r-r_{1}\right)  <0,$ and has a unique extension
to a maximal interval with no double zero. From (i) it has a unique extension
to $\left(  0,\infty\right)  .$ In particular any local regular solution is
defined on $\left[  0,\infty\right)  $.
\end{proof}

\section{Asymptotic behaviour}

Next the function $y$ is supposed to be monotone, thus $w$ has a constant sign
near $0$ or $\infty,$ we can assume that $w>0.$

\begin{proposition}
\label{esup}Let $y$ be any solution of (\textbf{E}$_{y}$) strictly monotone
and positive near $s=\pm\infty$.$\medskip$

(1) Then $(\zeta,\sigma)$ has a limit $(\lambda,\mu)$ near $s,$ given by is
some of the values
\begin{align}
A_{\gamma}  &  =\left(  -\gamma,\varepsilon\frac{\alpha+\gamma}{N+\gamma
}\right)  ,\quad A_{r}=\left(  0,\varepsilon\alpha/N\right)  ,\quad A_{\alpha
}=\left(  \alpha,0\right)  ,\nonumber\\
L_{\eta}  &  =\eta\left(  1,\infty\right)  (\text{if }p\neq N),\quad
L_{+}=\left(  0,\infty\right)  (\text{if }p\geqq N),\quad L_{-}=\left(
0,-\infty\right)  (\text{if }p>N). \label{poi}%
\end{align}

(2) More precisely,$\medskip$

\noindent(i) Either $\varepsilon\left(  \gamma+\alpha\right)  <0$ and $(y,Y)$
converges to $\pm M_{\ell}.$ Then $(\lambda,\mu)=A_{\gamma}$ and
($\varepsilon=1,s=\infty)$ or ($\varepsilon=-1,s$=-$\infty$ for $\alpha
\leqq\alpha^{\ast},$ $s=\infty$ for $\alpha>\alpha^{\ast}).\medskip$

\noindent(ii) Or $(y,Y)$ converges to $(0,0).$ Then ($s=\infty$ and
$-\gamma<\alpha)$ or ($s=-\infty$ and $\alpha<-\gamma)$, or ($s=\varepsilon
\infty$ and $\alpha=-\gamma)$ and $(\lambda,\mu)=A_{\alpha}.\medskip$

\noindent(iii) Or $\lim_{\tau\rightarrow s}y=\infty.$ Then $s=-\infty.$ If
$p<N,$ then $(\lambda,\mu)=A_{r}$ or $L_{\eta}$.\ If $p=N,$ then $(\lambda
,\mu)=A_{r}$ or $L_{+}.$ If $p>N,$ then $(\lambda,\mu)=A_{r},L_{\eta},L_{+}$
or $L_{-}.$
\end{proposition}

\begin{proof}
(1) The functions $Y,\sigma,\zeta$ are also monotone, and by definition
$\zeta\sigma>0$. Thus $\zeta$ has a limit $\lambda\in\left[  -\infty
,\infty\right]  $ and $\sigma$ has a limit $\mu\in\left[  -\infty
,\infty\right]  $, and $\lambda\mu\geqq0.\medskip$

\noindent(i) $\lambda$ is finite. Indeed if $\lambda=\pm\infty,$ then
$f=1/\zeta$ tends to $0.$ From (\ref{ff}), either $\mu=\pm\infty,$ then
$f^{\prime}$ tends to $-1,$ which is imposible; or $\mu$ is finite, thus
$\mu=\varepsilon$ from system (\textbf{Q}), then $f^{\prime}$ tends to
$(2-p)/(p-1),$ which is still contradictory. $\medskip$

\noindent(ii) Either $\mu$ is finite, thus $(\lambda,\mu)$ is a stationary
point of system \textbf{(Q)}, equal to $A_{\gamma},A_{r}$ or $A_{\alpha
}.\medskip$

\noindent(iii) 0r $\mu=\pm\infty$ and $(\lambda,0)$ is a stationary point of
system (\textbf{P}). $\medskip$

$\bullet$ If $p\neq N,$ either $\lambda=\eta\neq0$ and $(\lambda,\mu)=L_{\eta
};$ or $\lambda=0$ and $(\lambda,\mu)=L_{+}$ or $L_{-}.$ In the last case
$(\zeta,\psi)$ converges to $(0,0),$ and $\zeta^{\prime}/\psi^{\prime}%
=-(\eta\zeta/N\psi)(1+o(1)),$ thus $\eta<0,$ that means $p>N.\medskip$

$\bullet$ If $p=N,$ then again $(\zeta,\psi)$ converges to $(0,0),$ thus
$\mu=$ $\pm\infty,$ and $\psi^{\prime}=N\psi(1+o(1)),$ and necessarily
$s=-\infty.$ We make the substitution (\ref{sysd}) with $d=0.$ Then
$y_{0}(\tau)=w(r)$, and $y_{0}$ satisfies
\[
y_{0}^{\prime}=-\left\vert Y_{0}\right\vert ^{(2-p)/(p-1)}Y_{0}=-\zeta
y_{0}=o(y_{0}),\qquad Y_{0}^{\prime}=\varepsilon e^{p\tau}y_{0}(\alpha
-\zeta)=\varepsilon e^{p\tau}y_{0}\alpha(1+o(1).
\]
Thus for any $\upsilon>0,$ we get $y_{0}=O(e^{-\upsilon\tau})$ and
$1/y_{0}=O(e^{\upsilon\tau}).$ Then $Y_{0}^{\prime}$ is integrable, and
$Y_{0}$ has a finite limit $\left\vert k\right\vert ^{p-2}k.$ Suppose that
$k=0.$ Then $Y_{0}=O(e^{(p-\upsilon)\tau}),$ and $y_{0}$ has a finite limit
$a\geqq0.$ If $a\neq0,$ then $Y_{0}^{\prime}=\varepsilon\alpha ae^{p\tau
}(1+o(1));$ in turn $Y_{0}=p^{-1}\varepsilon\alpha ae^{p\tau}(1+o(1)),$ and
$\psi=e^{p\tau}y_{0}/Y_{0}$ does not tend to $0.$ If $a=0,$ then
$y_{0}=O(e^{p^{\prime}\tau}),$ which contradicts the estimate of $1/y_{0}.$
Thus $k>0$ and
\begin{equation}
y_{0}=-k\tau(1+o(1),\qquad Y_{0}=k^{p-1}(1+o(1)); \label{hlm}%
\end{equation}
hence $(\lambda,\mu)=L_{+}.\medskip$

(2) Since $y$ is monotone, we encounter one of the three following cases:
$\medskip$

\noindent(i) $(y,Y)$ converges to $\pm M_{\ell}.$ Then $(\lambda
,\mu)=A_{\gamma}$ and $M_{\ell}$ is a source (or a weak source) for
$\alpha\leqq\alpha^{\ast},$ a sink for $\alpha>\alpha^{\ast}$. \medskip

\noindent(ii) $y$ tends to $0.$ Since $\lambda$ is finite, $(y,Y)$ converges
to $(0,0).$ And $\left\vert \sigma\right\vert =$ $\left\vert \zeta\right\vert
^{p-1}y^{p-2}$ tends to $0,$ thus $(\lambda,\mu)=A_{\alpha}.$ If
$-\gamma<\alpha,$ seeing that $y^{\prime}=-y(\gamma+\zeta)<0$ we find
$s=\infty.$ If $\alpha<-\gamma,$ then $s=-\infty.$ If $\alpha=-\gamma<0,$ then
$\varepsilon\left(  \gamma+\zeta\right)  >0,$ from the first equation of
\textbf{(Q),} thus $\varepsilon y^{\prime}<0,$ hence $s=\varepsilon\infty
.$\medskip

\noindent(iii) $y$ tends to $\infty.$ Either $\lambda\neq0,$ thus $\left\vert
\sigma\right\vert =$ $\left\vert \zeta\right\vert ^{p-1}y^{p-2}$ tends to
$\infty,$ and $\lambda=\eta$ from system (\textbf{Q}), thus $p\neq N,$
$(\lambda,\mu)=L_{\eta}$. Or $\lambda=0$ and $\mu$ is finite, thus
$\mu=\varepsilon\alpha/N,$ $(\lambda,\mu)=A_{r}$. Or $(\lambda,\mu)=L_{0};$
then either $p=N,$ $L_{0}=L_{\eta},$ or $p>N.$ In any case, $y^{\prime
}=-y(\gamma+\zeta)<0,$ from (\ref{rel}), hence $s=-\infty.$
\end{proof}

Next we apply these results to the functions $w:$

\begin{proposition}
\label{dbw}We keep the assumptions of Proposition \ref{esup}. Let $w$ be the
solution of (\textbf{E}$_{w}$) associated to $y$ by (\ref{cha}).$\medskip$

\noindent(i) If $(\lambda,\mu)=A_{\gamma}$ (near 0 or $\infty),$ then%
\begin{equation}
\lim r^{-\gamma}w=\ell. \label{gam}%
\end{equation}

\noindent(ii) If $(\lambda,\mu)=A_{\alpha}$ (near 0 or $\infty),$ then
\begin{equation}
\lim r^{\alpha}w=L>0\qquad\qquad\qquad\qquad\qquad\qquad\qquad\qquad
\quad\text{if }\alpha\neq-\gamma, \label{val}%
\end{equation}%
\begin{equation}
\lim r^{-\gamma}(\ln r)^{1/(p-2)}w=((p-2)\gamma^{p-1}(N+\gamma))^{-1/(p-2)}%
\qquad\text{if }\alpha=-\gamma. \label{kap}%
\end{equation}

\noindent(iii) If $p<N$ and $(\lambda,\mu)=L_{\eta},$ then
\begin{equation}
\lim_{r\rightarrow0}r^{\eta}w=c>0. \label{vet}%
\end{equation}

\noindent(iv) If $p>N$ and $(\lambda,\mu)=L_{\eta},$ then%
\begin{equation}
\lim_{r\rightarrow0}r^{-\left\vert \eta\right\vert }w=c>0. \label{viu}%
\end{equation}

\noindent(v) If $p=N$ and $(\lambda,\mu)=L_{+},$ then
\begin{equation}
\lim_{r\rightarrow0}\left\vert \ln r\right\vert ^{-1}w=k>0\text{,}\qquad
\lim_{r\rightarrow0}rw^{\prime}=-k\qquad\text{if }p=N. \label{vpk}%
\end{equation}

\noindent(vi) If $p>N$ and $(\lambda,\mu)=L_{+}$, or $L_{-},$ then
\begin{equation}
\lim_{r\rightarrow0}w=a>0,\qquad\lim_{r\rightarrow0}(-r^{(N-1)/(p-1)}%
w^{\prime})=c>0, \label{via+}%
\end{equation}

or
\begin{equation}
\lim_{r\rightarrow0}w=a>0,\qquad\lim_{r\rightarrow0}(-r^{(N-1)/(p-1)}%
w^{\prime})=c<0. \label{via-}%
\end{equation}

\end{proposition}

\begin{proof}
(i) This follows directly from (\ref{cha}).\medskip

\noindent(ii) From (\ref{dze}), $rw^{\prime}(r)=-\alpha w(r)(1+o(1).$ We are
lead to three cases.\medskip

$\bullet$ Either $-\gamma<\alpha,$ and $s=\infty$. For any $\upsilon>0,$ we
find $w=O(r^{-\alpha+\upsilon})$ and $1/w=O(r^{\alpha+\upsilon})$ near
$\infty$ and $w^{\prime}=O(r^{-\alpha-1+\upsilon})$.\ Then $J_{\alpha}%
^{\prime}(r)=O(r^{\alpha(2-p)-p-1+\upsilon}),$ hence $J_{\alpha}^{\prime}$ is
integrable, $J_{\alpha}$ has a limit $L.$ And $\lim r^{\alpha}w=L,$ seeing
that $J_{\alpha}(r)=r^{\alpha}w(1+o(1)).$ If $L=0,$ then $r^{\alpha
}w=O(r^{\alpha(2-p)-p+\upsilon}),$ which contradicts the estimate of
$1/w=O(r^{\alpha+\upsilon})$ for $\upsilon$ small enough. Thus $L>0$.\medskip

$\bullet$ Or $\alpha<-\gamma$ and $s=-\infty$. For any $\upsilon>0,$ we find
$w=O(r^{-\alpha-\upsilon})$ and $1/w=O(r^{\alpha+\upsilon})$ near $0$ and
$w^{\prime}=O(r^{-\alpha-1-\upsilon})$.\ Then $J_{\alpha}^{\prime
}(r)=O(r^{\alpha(2-p)-p-1-\upsilon}),$ and $J_{\alpha}^{\prime}$ is still
integrable, $J_{\alpha}$ has a limit $L,$ and $\lim r^{\alpha}w=L.$ If $L=0,$
then $r^{\alpha}w=O(r^{\alpha(2-p)-p-\upsilon}),$ which contradicts the
estimate of $1/w$. Thus again $L>0$.\medskip

$\bullet$ Or $\alpha=-\gamma$ and $s=\varepsilon\infty.\ $Then $Y=-\gamma
^{p-1}y^{p-1}(1+o(1)),$ and $\mu=0,$ thus $y-\varepsilon Y=y(1+o(1)).\ $From
System (\textbf{S}),
\[
\left(  y-\varepsilon Y\right)  ^{\prime}=\varepsilon(N+\gamma)Y=-\varepsilon
(N+\gamma)\gamma^{p-1}\left(  y-\varepsilon Y\right)  ^{p-1}(1+o(1)).
\]
Then $y=(N+\gamma)\gamma^{p-1}(p-2)\left\vert \tau\right\vert )^{-1/(p-2)}%
(1+o(1)),$ which is equivalent to (\ref{kap}).\medskip

\noindent(iii) From (\ref{dze}), we get $rw^{\prime}(r)=-\eta w(r)(1+o(1).$ We
use (\ref{cge}) with $d=\eta,$ thus $y_{\eta}=r^{\eta}w.$ We find $y_{\eta
}=O(e^{-\upsilon\tau}),$ $1/y_{\eta}=O(e^{-\upsilon\tau}),$ in turn $Y_{\eta
}=O(e^{-\upsilon\tau})$. From (\ref{sysd}), $Y_{\eta}^{\prime}%
=O(e^{(p+(p-2)\eta-\upsilon)\tau}),$ thus $Y_{\eta}^{\prime}$ is integrable,
hence $Y_{\eta}$ has a finite limit. Now $(e^{-\eta\tau}y_{\eta})^{\prime
}=-e^{-\eta\tau}Y_{\eta}^{1/(p-1)},$ and $\eta>0,$ thus $y_{\eta}$ has a limit
$c.$ If $c=0,$ then $Y_{\eta}=O(e^{(p+(p-2)\eta-\upsilon)\tau}),$ $y_{\eta
}=O(e^{((p+(p-2)\eta)/(p-1)-\upsilon)\tau}),$ which contradicts $1/y_{\eta
}=O(e^{-\upsilon\tau})$ for $\upsilon$ small enough. Then (\ref{vet})
holds.\medskip

\noindent(iv) As above, $Y_{\eta}$ has a finite limit. In turn $r^{-\left\vert
\eta\right\vert +1}w^{\prime}=\left\vert Y_{\eta}\right\vert ^{(2-p)/(p-1)}%
Y_{\eta}$ has a limit $c\left\vert \eta\right\vert $ and $w$ has a limit
$a\geqq0.\ $From (\ref{dze}), $rw^{\prime}=\left\vert \eta\right\vert
w(1+o(1),$ hence $a=0.$ Then $c\geqq0;$ if $b=0,$ then $Y<0,$ the function
$v=-e^{(\gamma+N)\tau}Y>0$ tends to $0$ and
\[
v^{\prime}=-e^{(\gamma+N)\tau}\varepsilon(\alpha-\eta)y(1+o(1))=-\varepsilon
(\alpha-\eta)\left\vert \eta\right\vert e^{-(\gamma+N)(p-2)/(p-1)\tau
}v^{1/(p-1)};
\]
we reach again a contradiction.Thus $a=0$ and $c>0,$ and (\ref{viu})
holds.\medskip

\noindent(v) Assertion (\ref{vpk}) follows from (\ref{hlm}).\medskip

\noindent(vi) Here $rw^{\prime}=o(w),$ thus $w+\left\vert w^{\prime
}\right\vert =O(r^{-k})$ for any $k>0.$ Then $J_{N}^{\prime}$ is integrable,
$J_{N}$ has a limit at $0,$ and $\lim_{r\rightarrow0}$ $r^{N}w=0.\ $Thus
$\lim_{r\rightarrow0}r^{(N-1)/(p-1)}w^{\prime}=-c\in\mathbb{R},$
$\lim_{r\rightarrow0}J_{N}=-\varepsilon\left\vert c\right\vert ^{p-2}c,$
$\lim_{r\rightarrow0}w=a\geq0.$ If $c=0,$ then $J_{N}(r)=%
{\displaystyle\int\limits_{0}^{r}}
J_{N}^{\prime}(s)ds,$ implying that $\lim_{r\rightarrow0}w^{\prime}=0.$ Either
$a>0$ and then $w$ is regular, then $\lim_{\tau\rightarrow-\infty}%
\sigma=\varepsilon;$ or $a=0,$ then $w^{\prime}>0$ and $(w^{\prime}%
)^{p-1}=O(rw);$ in both cases we get a contradiction. Thus $c\neq0.$ If $a=0,$
we find $\lim_{\tau\rightarrow-\infty}\zeta=\eta,$ which is not true, hence
$a>0.$ In any case (\ref{via+}) or (\ref{via-}) holds.\medskip
\end{proof}

Now we study the cases where $y$ is not monotone, and eventually changing sign.

\begin{proposition}
\label{nmo} Suppose $\varepsilon=-1.$ Let $w\not \equiv 0$ be any solution of
(\textbf{E}$_{w}$).\medskip\ 

\noindent(i) If $\alpha\leqq-\gamma,$ then $w$ is oscillating near $0$ at
$\infty.$\medskip

\noindent(ii) If $\alpha<0,$ then $y$ and $Y$ are bounded at $\infty.$
\end{proposition}

\begin{proof}
(i) Suppose by contradiction that $w\geqq0$ for large $r,$ then $y\geqq0$ for
large $\tau.$ If $y>0$ near $\infty,$ then from Proposition \ref{com}, either
$y$ is constant, which is impossible since $(0,0)$ is the unique stationary
point; or $y$ is strictly monotone, which contradicts Proposition \ref{esup}.
Then there exists a sequence ($\tau_{n})$ tending to $\infty$ such that
$y(\tau_{n})=y^{\prime}(\tau_{n})=0;$ from Theorem \ref{prol}, $y\equiv0$ on
$\left(  -\infty,\tau_{n}\right)  ,$ thus $y\equiv0$.\medskip

\noindent(ii) Consider the function
\[
\tau\mapsto R(\tau)=\frac{y^{2}}{2}+\frac{\left\vert Y\right\vert ^{p^{\prime
}}}{p^{\prime}\left\vert \alpha\right\vert };
\]
it satisfies%
\[
R^{\prime}\left(  \tau\right)  =-\gamma y^{2}+\frac{1}{\left\vert
\alpha\right\vert }\left\vert Y\right\vert ^{2/(p-1)}-\frac{N+\gamma
}{\left\vert \alpha\right\vert }\left\vert Y\right\vert ^{p^{\prime}}.
\]
From the Young inequality,%
\[
\left\vert \alpha\right\vert (R^{\prime}\left(  \tau\right)  +\gamma
R(\tau))=\left\vert Y\right\vert ^{2/(p-1)}-(N+\frac{1}{p-2})\left\vert
Y\right\vert ^{p^{\prime}}\leqq(\frac{2}{Np+\gamma})^{(p-2)/2}\leqq1
\]
thus $R(\tau)$ is bounded for large $\tau,$ at least by $1/\left\vert
\alpha\right\vert \gamma.$
\end{proof}

\begin{proof}
\begin{proposition}
\label{cyc} (i) Assume $\varepsilon=1,$ or $\varepsilon=-1,$ $\alpha
\not \in \left(  \alpha_{2},\alpha_{1}\right)  .$ Then for any trajectory of
system (\textbf{S})\textbf{ }in $\mathcal{Q}_{4}$ near $\pm\infty,$ $y$ is
strictly monotone near $\pm\infty$.\medskip

\noindent(ii) Assume $\varepsilon=1,$ and $\alpha\leqq\alpha^{\ast}$ or
$-p^{\prime}\leqq\alpha.$ Then system (\textbf{S})\textbf{ }admits no cycle in
$\mathcal{Q}_{4}$ (or $\mathcal{Q}_{2}).$
\end{proposition}
\end{proof}

\begin{proof}
(i) In any case $M_{\ell}$ is a node point. Following \cite[Theorem 2.24]%
{Bi1}, we use the linearization defined by (\ref{tran}). Consider the line $L$
given by the equation $A\overline{y}+\overline{Y}=0$, where $A$ is a real
parameter. The points of $L$ are in $\mathcal{Q}_{4}$ whenever $\overline
{Y}<(\gamma\ell)^{p-1}$ and $-\ell<\overline{y}.$ We get
\[
A\overline{y}^{\prime}+\overline{Y}^{\prime}=\left(  \varepsilon\nu
(\alpha)A^{2}+(N+\nu(\alpha))A+\varepsilon\alpha\right)  \overline
{y}+(A+\varepsilon)\Psi(\overline{Y}).
\]
From (\ref{lta}), apart from the case $\varepsilon=1,\alpha=N,$ we can find an
$A$ such that
\[
\varepsilon\nu(\alpha)A^{2}+(N+\nu(\alpha))A+\varepsilon\alpha=0,
\]
and $A+\varepsilon\neq0$. Moreover $\Psi(\overline{Y})\leqq0$ on
$L\cap\mathcal{Q}_{4}.$ Indeed $(p-1)\Psi^{\prime}(t)=-((\gamma\ell
)^{p-1}-t)^{(2-p)/(p-1)}+(\gamma\ell)^{2-p},$ thus $\Psi$ has a maximum $0$ on
$\left(  -\infty,(\delta\ell)^{p-1}\right)  $ at point $0$. Then the
orientation of the vector field does not change along $L\cap\mathcal{Q}_{4}.$
In particular $y$ cannot oscillate around $\ell,$ thus $y$ is monotone, from
Proposition \ref{com}. If $\varepsilon=1,\alpha=N,$ then $Y\equiv y\in\left(
\ell,\infty\right)  $ defines the trajectory $\mathcal{T}_{r}$, corresponding
to the solutions given by (\ref{fut}) with $K>0$. No solution $y$ can
oscillate around $\ell,$ since the trajectory cannot meet $\mathcal{T}_{r}.$
\medskip

\noindent(ii) Suppose that there exists a cycle in $\mathcal{Q}_{4}$.\medskip

$\bullet$ Assume $\alpha\leqq\alpha^{\ast}.$ Here $M_{\ell}$ is a source, or a
weak source, from Proposition \ref{ws}. Any trajectory starting from $M_{\ell
}$ at $-\infty$ has a limit cycle in $\mathcal{Q}_{1},$ which is attracting at
$\infty.$ Writing System (\textbf{S}) under the form $y^{\prime}%
=f_{1}(y,Y),Y^{\prime}=f_{2}(y,Y),$ the mean value of the Floquet integral on
the period $\left[  0,\mathcal{P}\right]  $ is given by
\begin{equation}
I=%
{\displaystyle\oint}
(\frac{\partial f_{1}}{\partial y}(y,Y)+\frac{\partial f_{2}}{\partial
Y}(y,Y))d\tau=%
{\displaystyle\oint}
(\frac{\left\vert Y\right\vert ^{(2-p)/(p-1)}}{p-1}-2\gamma-N)d\tau.
\label{hhh}%
\end{equation}
Such a cycle is not unstable, thus $I\leqq0.$ Now
\[%
{\displaystyle\oint}
(\alpha y^{\prime}-\gamma Y^{\prime})d\tau=0=(\alpha+\gamma)%
{\displaystyle\oint}
\left\vert Y\right\vert ^{1/(p-1)}d\tau-\gamma(\gamma+N)%
{\displaystyle\oint}
\left\vert Y\right\vert d\tau.
\]
From the Jensen and H\"{o}lder inequalities, since $1/(p-1)<1,$
\[
\gamma(\gamma+N)(%
{\displaystyle\oint}
\left\vert Y\right\vert ^{1/(p-1)}d\tau)^{p-2}\leqq\alpha+\gamma,
\]%
\[
1\leqq\left(
{\displaystyle\oint}
\left\vert Y\right\vert ^{(2-p)/(p-1)})d\tau\right)  \left(
{\displaystyle\oint}
\left\vert Y\right\vert ^{1/(p-1)}d\tau\right)  ^{p-2}\leqq\frac
{(p-1)(2\gamma+N)}{\gamma(\gamma+N)}(\alpha+\gamma),
\]
then $\alpha^{\ast}<\alpha,$ which is contradictory.\medskip

$\bullet$ Assume $-p^{\prime}\leqq\alpha<0.$ Consider the functions
$y_{\alpha}=e^{(\alpha+\gamma)\tau}y$ and $Y_{\alpha}=e^{(\alpha
+\gamma)(p-1)\tau}Y$ defined by (\ref{cge}) with $d=\alpha.$ They vary
respectively from $0$ to $\infty$ and from $0$ to $-\infty.$ They have no
extremal point. Indeed at such a point, from (\ref{yde}) and (\ref{Yd})
$y_{\alpha}^{\prime\prime}$ or $Y_{\alpha}^{\prime\prime}$ have a strict
constant sign for $\alpha\neq\eta,p^{\prime},$ which is contradictory. If
$\alpha=\eta$ or $p^{\prime},$ from uniqueness $y_{\alpha}$ or $Y_{\alpha}$ is
constant, thus $y$ or $Y$ is monotone, which is impossible. In any case
$y_{\alpha}^{\prime}>0>Y_{\alpha}^{\prime}$ on $\left(  -\infty,\infty\right)
.$ Next, from (\ref{phi}) and (\ref{phu}),
\begin{equation}
\frac{y_{\alpha}^{\prime\prime}}{y_{\alpha}^{\prime}}+\eta-2\alpha-\frac
{1}{p-1}Y^{(2-p)/(p-1)}=\alpha(\eta-\alpha)\frac{y_{\alpha}}{y_{\alpha
}^{\prime}}, \label{aa}%
\end{equation}%
\begin{equation}
\frac{Y_{\alpha}^{\prime\prime}}{Y_{\alpha}^{\prime}}+(p-1)(\eta
-2\alpha-p^{\prime})-\frac{1}{p-1}Y^{(2-p)/(p-1)}=(p-1)^{2}(\eta
-\alpha)(p^{\prime}+\alpha)\frac{Y_{\alpha}}{Y_{\alpha}^{\prime}}. \label{bb}%
\end{equation}
Let us integrate on the period $\mathcal{P}.$ If $\eta\leqq\alpha<0,$ then
$\eta-N-2(\alpha+\gamma)\geqq0$ from (\ref{aa}), which is contradictory. If
$-p^{\prime}\leqq\alpha<\eta,$ then $-2(\alpha+p^{\prime}+\gamma)>0$ from
(\ref{bb}), still contradictory.
\end{proof}

\section{New local existence results\label{new}}

At Proposition \ref{esup} we gave all the \textit{possible} behaviours of the
positive solutions near $\pm\infty.$ Next we prove their existence, and
uniqueness or multiplicity. The case $p>N$ is very delicate.

\begin{theorem}
\label{teta} (i) Suppose $p<N.$ In the phase plane $(y,Y)$ of system
(\textbf{S}) there exist an infinity of trajectories $\mathcal{T}_{\eta}$ such
that $\lim_{\tau\rightarrow-\infty}(\zeta,\sigma)=L_{\eta};$ the corresponding
$w$ satisfy (\ref{vet}).\medskip

\noindent(ii) Suppose $p>N.$ There exist a unique trajectory $\mathcal{T}_{u}$
such that $\lim_{\tau\rightarrow-\infty}(\zeta,\sigma)=L_{\eta}$; in other
words for any $c\neq0,$ there exists a unique solution $w$ of equation
(\textbf{E}$_{w}$) such that (\ref{viu}) holds.
\end{theorem}

\begin{proof}
Suppose that such a trajectory exists in the plane $(y,Y)$. In the phase plane
$\left(  \zeta,\psi\right)  $ of System (\textbf{P}), $\zeta$ and $\psi$ keep
a strict constant sign, because the two axes $\zeta=0$ and $\psi=0$ contain
particular trajectories, and $\left(  \zeta,\psi\right)  $ converges to
$(\eta,0)$ at $-\infty.$ Reciprocally, setting $\zeta=\eta+\bar{\zeta},$ the
linearized problem at point $(\eta,0)$
\[
\bar{\zeta}^{\prime}=\eta\bar{\zeta}+\eta(\alpha-\eta)\varepsilon
\psi/(p-1),\qquad\psi^{\prime}=(N-\eta)\psi,
\]
admits the eigenvalues $\eta$ and $N-\eta.$ The trajectories linked to the
eigenvalue $\eta$ are tangent to the line $\psi=0.$\medskip

\noindent(i) Case $p<N.$ Then $\eta>0,$ and $(\eta,0)$ is a source. In the
plane $\left(  \zeta,\psi\right)  $ there exist an infinity of trajectories,
starting from this point at $-\infty,$ such that $\psi>0,$ and $\lim
_{\tau\rightarrow-\infty}\zeta=\eta,$ thus $\zeta>0.$ In the phase plane
$(y,Y),$ setting $y=(\psi\left\vert \zeta\right\vert ^{p-2}\zeta)^{2-p}$ and
$Y=y/\psi,$ they correspond to an infinity of trajectories in the plane
$(y,Y)$ such that $\lim_{\tau\rightarrow-\infty}(\zeta,\sigma)=L_{\eta},$ and
(\ref{vet}) holds from Proposition (\ref{dbw}).\medskip

\noindent(ii) Case $p>N.$ Then $\eta<0,$ and $(\eta,0)$ is a saddle point. In
the plane $\left(  \zeta,\psi\right)  ,$ there exists a unique trajectory
starting from $(\eta,0),$ tangentially to the vector $\left(  \eta(\alpha
-\eta)\varepsilon/(p-1),N-\eta\right)  ,$ with $\psi<0;$ it defines a unique
trajectory $\mathcal{T}_{u}$ in the plane $(y,Y)$, and (\ref{viu}) holds. From
Remark \ref{scaling}, we get a solution for any $c\neq0.$
\end{proof}

\begin{theorem}
\label{tzero}(i) Suppose $p=N.$ In the phase plane $(y,Y),$ there exists an
infinity of trajectories $\mathcal{T}_{+}$ such that $\lim_{\tau
\rightarrow-\infty}(\zeta,\sigma)=L_{+};$ then $w$ satisfies (\ref{vpk}%
).\medskip

\noindent(ii) Suppose $p>N.$ Then there exist an infinity of trajectories
$\mathcal{T}_{+}$ (resp. $\mathcal{T}_{-})$ such $\lim_{\tau\rightarrow
-\infty}(\zeta,\sigma)=L_{+}$ (resp. $L_{-});$ then the corresponding
solutions $w$ of (\textbf{E}$_{w}$) satisfy (\ref{via+}) (resp. (\ref{via-}).

More precisely for any $k>0$ (for $p=N)$ or any $a>0$ and $c\neq0$ (for $p>N)$
there exists a unique function $w$ satisfying those conditions.
\end{theorem}

\begin{proof}
If $\lim_{\tau\rightarrow-\infty}(\zeta,\sigma)=L_{\pm},$ then $\lim
_{\tau\rightarrow-\infty}(\zeta,\psi)=(0,0),$ with $\zeta\psi>0$ in case of
$L_{+},$ $\zeta\psi<0$ in case of $L_{-.}.$ The linearization of System
(\textbf{P}) near $(0,0)$ is given by
\[
\zeta^{\prime}=\left\vert \eta\right\vert \zeta,\qquad\psi^{\prime}=N\psi.
\]

\noindent(i) Case $p=N.$ The phase plane study is delicate because $0$ is a
center, thus we use a fixed method. Suppose that such a trajectory exists, and
consider the substitution (\ref{cge}) with $d=0.$ From (\ref{hlm}), there
exists $k>0$ such that $\zeta=\left\vert Y_{0}\right\vert ^{(2-p)/(p-1)}%
/y_{0}=-\tau^{-1}(1+o(1))>0$, and $\psi=-k^{2-p}\tau e^{N\tau}(1+o(1))>0.$
Then $\zeta^{\prime}=\tau^{-2}(1+o(1))$ from System (\textbf{P}). The
function
\[
V=\psi e^{-N/\zeta}\zeta
\]
satisfies $\lim_{\tau\rightarrow-\infty}V=k^{2-p},$ and
\[
V^{\prime}=\frac{Ve^{N/\zeta}}{(N-1)\zeta^{2}}(\varepsilon\left(  \alpha
-\zeta\right)  (N-(N-2)\zeta)V+2N(N-1)\zeta^{2}e^{-N/\zeta}).
\]
Thus $\varepsilon\alpha(V-k^{2-p})>0$ near $-\infty.$ Moreover $\lim
_{\tau\rightarrow-\infty}\zeta^{\prime}/V^{\prime}=0,$ so that $\zeta$ can be
considered as a function of $V$ near $k^{2-p},$ with $\lim_{V\rightarrow
k^{2-p}}\zeta=0$ and
\[
\frac{d\zeta}{dV}=K(V,\zeta),\quad\quad K(V,\zeta):=\frac{\zeta^{2}}{V}%
\frac{\varepsilon\left(  \alpha-\zeta\right)  V+(N-1)\zeta^{2}e^{-N/\zeta}%
}{\varepsilon\left(  \alpha-\zeta\right)  (N-(N-2)\zeta)V+2N(N-1)\zeta
^{2}e^{-N/\zeta}}.
\]
Reciprocally, extending the function $\zeta^{2}e^{-N/\zeta}$ by $0$ for
$\zeta\leqq0,$ the function $K$ is of class $C^{1}$ near $(k^{2-p},0).\ $ For
any $k>0,$ there exists a unique local solution $V\mapsto\zeta(V)$ on a
interval $\mathcal{V}$ where $\varepsilon\alpha(V-k^{2-p})>0,$ such that
$\zeta(k^{2-p})=0.$ And $d\zeta/dV=(\zeta^{2}/Nk^{2-p})(1+o(1))$ near $0,$
thus $\zeta>0.$ In the plane $(\zeta,\psi)$, taking one point $P$ on the curve
$\mathcal{C}=\left\{  (\zeta(V),V\zeta(V)e^{N/\zeta(V)}):v\in\mathcal{V}%
\right\}  ,$ there exists a unique solution of System (\textbf{P}) issued from
$P$ at time $0.$ Its trajectory is on $\mathcal{C},$ thus it converges to
$(0,0),$ with $\zeta,\psi>0.$ It corresponds to a unique trajectory
$\mathcal{T}_{+}$ in the plane $(y,Y),$ and $(\zeta,\sigma)$ converges to
$L_{+},$ as $\tau$ tends to $-\infty,$ from Proposition \ref{esup}.\ The
corresponding functions $w$ satisfy (\ref{vpk}) from Proposition
(\ref{dbw}).\medskip

\noindent(ii) Case $p>N.$ Here $(0,0)\ $is a source for System (\textbf{P}).
The lines $\zeta=0$ and $\psi=0$ contain trajectories. There exists an
infinity of trajectories converging to $(0,0),$ with $\zeta\psi\neq0;$
moreover, if $N\geqq2,$ then $\left\vert \eta\right\vert <N,$ thus $\lim
_{\tau\rightarrow-\infty}(\psi/\zeta)=0.$ Our claim is more precise. Given
$a>0$ and $c\neq0,$ we look for a solution $w$ of (\textbf{E}$_{w})$ such that
$\lim_{r\rightarrow0}w=a,$ $\lim_{r\rightarrow0}r^{\eta+1}w^{\prime}=-c.$ By
scaling we can assume $a=1.$ If $w_{1}$ is a such a solution, then $\zeta$ and
$\psi$ have the sign of $c$ near $0,$ and $\zeta\left(  \tau\right)
=ce^{\left\vert \eta\right\vert \tau}(1+o(1))$ and $\left\vert c\right\vert
^{p-2}c\psi\left(  \tau\right)  =e^{N\tau}(1+o(1)).$ The function
\[
v=c(\left\vert c\right\vert ^{p-2}c\psi)^{1/\kappa}/\zeta,\qquad\text{with
}\kappa=N/\left\vert \eta\right\vert >1,
\]
satisfies $\lim_{\tau\rightarrow-\infty}v=1,$ and can be expressed locally as
a function of $\zeta,$ and
\[
\frac{dv}{d\zeta}=H(\zeta,v),\qquad H(\zeta,v):=-\frac{v}{\kappa}%
\frac{(p-1)(\kappa+1)+\varepsilon(\kappa-p+1)\left\vert c\right\vert
^{1-p-\kappa}(\zeta-\alpha)\left\vert \zeta\right\vert ^{\kappa-1}v^{\kappa}%
}{(p-1)(\zeta-\eta)+\varepsilon\left\vert c\right\vert ^{1-p-\kappa}%
(\alpha-\zeta)\left\vert \zeta\right\vert ^{\kappa-1}\zeta v^{\kappa}}.
\]
Reciprocally, there exists a unique solution $\zeta\mapsto v(\zeta)$ of this
equation on a small interval $\left[  0,hc\right)  ,$ with $h>0,$ such that
$v(0)=1.$ Indeed $H$ is locally continuous in $\xi$ and $C^{1}$ in $v.$ Taking
one point $P$ on the curve $\mathcal{C}^{\prime}=\left\{  (\zeta,\left\vert
c\right\vert ^{1-p-\kappa}\left\vert \zeta\right\vert ^{\kappa-1}\zeta
v(\zeta)):\zeta\in\left[  0,hc\right)  \right\}  ,$ there exists a unique
solution of System (\textbf{P}) issued from $P$ at time $0.$ Its trajectory is
on $\mathcal{C}^{\prime},$ thus converges to $(0,0)$ with $\zeta\psi>0.$ It
corresponds to a solution $(y,Y)$ of System (\textbf{S}), such that
$(\zeta,\sigma)$ converges to $L_{+},$ as $\tau$ tends to $-\infty,$ from
Proposition \ref{esup}.\ The corresponding function, called $w_{2},$ satisfies
$\lim_{r\rightarrow0}r^{\eta+1}w_{2}^{\gamma^{-1}\left\vert \eta\right\vert
-1}w_{2}^{\prime}=-c;$ thus $w_{2}$ has a limit $a_{2}$, and $\lim
_{r\rightarrow0}r^{\eta-1}w_{2}^{\prime}=a_{2}^{1-s}b.$ Moreover $a_{2}\neq0,$
because $a_{2}=0$ implies that $r^{-\gamma}w_{2}$ has a nonzero limit, thus
$(\zeta,\sigma)$ converges to $A_{\gamma}.$ The function $w(r)=a_{2}^{-1}%
w_{2}(a_{2}^{1/\gamma}r)$ satisfies $\lim_{r\rightarrow0}w=1,$ and
$\lim_{r\rightarrow0}r^{\eta-1}w^{\prime}=-c,$ and the proof is done.

\begin{theorem}
\label{talpha} (i) In the phase plane $(y,Y),$ for any $\alpha\neq0$ there
exists at least a trajectory $\mathcal{T}_{\alpha}$ converging to $(0,0)$ with
$y>0,$ and $\lim(\zeta,\sigma)=A_{\alpha}.$ The convergence holds at $\infty$
if $-\gamma<\alpha,$ or $-\infty$ if $\alpha<-\gamma$, or $\varepsilon\infty$
if $\alpha=-\gamma.\medskip$

\noindent(ii) If $\varepsilon(\gamma+\alpha)<0,$ $\mathcal{T}_{\alpha}$ is
unique, it is the unique trajectory converging to $(0,0)$ at $-\varepsilon
\infty$ with $y>0,$ and it depends locally continuously of $\alpha.$
\end{theorem}
\end{proof}

\begin{proof}
(i) Suppose that such a trajectory exists. Then $\tau$ tends to $\infty$ if
$-\gamma<\alpha,$ or $-\infty$ if $\alpha<-\gamma$, or $\varepsilon\infty$ if
$\alpha=-\gamma,$ from Proposition \ref{esup}. Consider System (\textbf{R}),
where $g,s$ and $\nu$ are defined by (\ref{gsnu}). Then $(g,s)$ converges to
$(-1/\alpha,0)$, with $gs>0,$ and $\nu$ tends to the same limits as $\tau$,
since $Y$ converges to $0.$ Reciprocally, in the plane $(g,s),$ let us show
the existence of a trajectory converging to $(-1/\alpha,0),$ different from
the line $s=0.$ Setting $g=-1/\alpha+\bar{g},$ the linearized system at this
point is
\[
\frac{d\bar{g}}{d\nu}=-\frac{\varepsilon}{p-1}\bar{g}+\frac{\eta-\alpha
}{\alpha^{2}}s,\qquad\frac{ds}{d\nu}=0,
\]
thus we find a center: the eigenvalues are $0$ and $\lambda=\varepsilon
/(p-1)$. Since the system is polynomial, it is known that System (\textbf{R})
admits a trajectory, depending locally continuously of $\alpha,$ such that
$sg>0,$ and tangent to the eigenvector $((p-1)(\eta-\alpha),\varepsilon
\alpha^{2})$. It satisfies $ds/d\nu=(p-2)(\alpha+\gamma)s^{2}(1+o(1)).$ Then
$ds/d\tau=-(p-2)\alpha(\alpha+\gamma)s(1+o(1)),$ thus $\tau$ tends to
$\pm\infty.$ And $\left\vert y\right\vert ^{p-2}=\left\vert s\right\vert $
$\left\vert g\right\vert ^{1/(p-1)},$ then $y$ tends to $0,$ $(y,Y)$ converges
to $(0,0),$ and $\lim(\zeta,\sigma)=A_{\alpha}.\medskip$

\noindent(ii) Suppose $\varepsilon(\gamma+\alpha)<0$. Consider two
trajectories $\mathcal{T}_{1},\mathcal{T}_{2}$ in the plane $(y,Y),$
converging to $(0,0)$ at $-\varepsilon\infty,$ with $y>0.$ They are different
from $\mathcal{T}_{\varepsilon}$ which converges at $\varepsilon\infty,$ thus
$\lim(\zeta_{i},\sigma_{i})=(\alpha,0)$ from Proposition \ref{esup}. Then
$\zeta_{1},\zeta_{2}$ can locally be expressed as a function of $y,$ and
\[
y\frac{d(\zeta_{1}-\zeta_{2})^{2}}{dy}=2(F(\zeta_{1},y)-F(\zeta_{2},y))\left(
\zeta_{1}-\zeta_{2}\right)
\]
near $0,$ where
\[
F(\zeta,y)=\frac{1}{\gamma+\zeta}(-\zeta(\zeta-\eta)+\frac{\varepsilon}%
{p-1}\left\vert \zeta y\right\vert ^{2-p}(\zeta-\alpha)).
\]
Then $(\zeta_{1}-\zeta_{2})^{2}$ is nonincreasing, seeing that $\partial
F/\partial\zeta(\zeta,y)=-((p-1)\varepsilon(\gamma+\alpha))^{-1}\left\vert
\alpha y\right\vert ^{2-p}(1+o(1)).$ Hence $\zeta_{1}\equiv\zeta_{2}$ near
$0,$ and $\mathcal{T}_{1}\equiv\mathcal{T}_{2}.$
\end{proof}

\section{The case $\varepsilon=1,$ $-\gamma\leqq\alpha$\label{one}}

In that Section and in Sections \ref{two}, \ref{three} and \ref{four} we
describe the solutions of (\textbf{E}$_{w}$). When we give a
\textit{uniqueness} result, we mean that $w$ is unique, \textit{up to a
scaling,} from Remark \ref{scaling}.

\begin{theorem}
\label{pin}Assume $\varepsilon=1,$ $-\gamma\leqq\alpha$ ($\alpha\neq
0)$.\medskip\ 

Any solution $w$ of (\textbf{E}$_{w}$) has a finite number of simple zeros,
and satisfies (\ref{val}) or (\ref{kap}) near $\infty$ or has a compact
support. Either $w$ is regular, or $\left\vert w\right\vert $ satisfies
(\ref{vet}),(\ref{vpk}), (\ref{viu}),(\ref{via+}) or (\ref{via-}) near $0$,
and there exist solutions of each type.\medskip

(1) Case $\alpha<N.$ All regular solutions have a strict constant sign, and
satisfy (\ref{val}) or (\ref{kap}) near $\infty.$ Moreover there exist (and
exhaustively, up to a symmetry)\medskip

\noindent(i) a unique nonnegative solution with (\ref{vet})or (\ref{vpk}) or
(\ref{via+})) near $0$, and compact support;

\noindent(ii) positive solutions with the same behaviour at $0$ and
(\ref{val}) or (\ref{kap}) near $\infty;$

\noindent(iii) solutions with one simple zero, and $\left\vert w\right\vert $
has the same behaviour at $0$ and $\infty;$

\noindent(iv) for $p>N,$ a unique positive solution with (\ref{viu}) near $0$,
and (\ref{val}) or (\ref{kap}) near $\infty$;

\noindent(v) for $p>N,$ positive solutions with (\ref{via-}) near $0$, and
(\ref{val}) or (\ref{kap}) near $\infty.$\medskip

(2) Case $\alpha=N.$ Then the regular (Barenblatt) solutions have a constant
sign with compact support. If $p\leqq N,$ all the other solutions are of type
(iii). If $p>N,$ there exist also solutions of type (iv) and (v).\medskip

(3) Case $\alpha>N.\medskip$

Either the regular solutions have $m$ simple zeros and satisfy satisfies
(\ref{val}) near $\infty.$ Then there exist

\noindent(vi) a unique solution with $m$ simple zeros, $\left\vert
w\right\vert $ satisfies (\ref{vet}), (\ref{vpk}) or(\ref{via+}) near $0,$
with compact support;

\noindent(vii) solutions with $m+1$ simple zeros, $\left\vert w\right\vert $
satisfies (\ref{vet}), (\ref{vpk}) or (\ref{via+}) near $0,$ and (\ref{val})
or (\ref{kap}) near $\infty;$

\noindent(viii) for $p>N,$ solutions with $m$ simple zeros, $\left\vert
w\right\vert $ satisfies (\ref{via+}),(\ref{viu}) or (\ref{via-}) near $0,$
and (\ref{val}) or (\ref{kap}) near $\infty.\medskip$

Or the regular solutions have $m$ simple zeros and a compact support. Then the
other solutions are of type (vii) or (viii).
\end{theorem}

\[%
\begin{array}
[c]{cc}%
\raisebox{-0pt}{\includegraphics[
natheight=187.289200pt,
natwidth=187.289200pt,
height=189.096pt,
width=189.096pt
]%
{../../../AMarie/MesArticles/107050selfp+2/pin1.bmp}%
}%
&
\raisebox{-0pt}{\includegraphics[
natheight=187.289200pt,
natwidth=187.289200pt,
height=189.096pt,
width=189.096pt
]%
{../../../AMarie/MesArticles/107050selfp+2/pin2.bmp}%
}%
\\
\text{th \ref{pin},fig1: }\varepsilon=1,N=2,p=3,\alpha=-2 & \text{th
\ref{pin},fig2: }\varepsilon=1,N=2,p=3,\alpha=1
\end{array}
\]%
\[%
\begin{array}
[c]{cc}%
\raisebox{-0pt}{\includegraphics[
natheight=187.289200pt,
natwidth=187.289200pt,
height=189.096pt,
width=189.096pt
]%
{../../../AMarie/MesArticles/107050selfp+2/pin3.bmp}%
}%
&
\raisebox{-0pt}{\includegraphics[
natheight=187.289200pt,
natwidth=187.289200pt,
height=189.096pt,
width=189.096pt
]%
{../../../AMarie/MesArticles/107050selfp+2/pin4.bmp}%
}%
\\
\text{th \ref{pin},fig3: }\varepsilon=1,N=2,p=3,\alpha=2 & \text{th
\ref{pin},fig4: }\varepsilon=1,N=2,p=3,\alpha=50
\end{array}
\]

\begin{proof}
All the solutions $w$ have a finite number of simple zeros, from Proposition
\ref{fini} and Theorem \ref{tpsilon}. Either they have a compact support. Or
$y$ has a strict constant sign and is monotone near $\infty$, and converge to
$(0,0)$ at $\infty,$ and (\ref{val}) or (\ref{kap}) holds$,$ from Propositions
\ref{com}, \ref{esup}. \medskip

In the phase plane $(y,Y),$ system (\textbf{S}) admits only one stationary
point $(0,0)$. The trajectory $\mathcal{T}_{r}$ starts in $\mathcal{Q}_{4}$
when $\alpha<0,$ in $\mathcal{Q}_{1}$ when $\alpha>0,$ and $\lim
_{\tau\rightarrow-\infty}y=\infty,$ with an asymptotical direction of slope
$\alpha/N$. From Propositions \ref{esup} and \ref{dbw} all the nonregular
solutions $\pm w$ satisfy (\ref{vet}), (\ref{vpk}), (\ref{viu}), (\ref{via+})
or (\ref{via-}) near $-\infty.$ The existence of solutions of any kind is
proved at Theorems \ref{teta} and \ref{tzero}. When $p\leqq N,$ they
correspond to trajectories $\pm\mathcal{T}_{\eta}$ such that $\mathcal{T}%
_{\eta}$ starts in $\mathcal{Q}_{1}$ with an infinite slope, in any case above
$\mathcal{T}_{r}.$ When $p>N,$ there is a unique trajectory $\mathcal{T}_{u}$
satisfying (\ref{viu}), starting in $\mathcal{Q}_{4},$ under $\mathcal{T}%
_{r};$ the trajectories $\mathcal{T}_{+}$ start from $\mathcal{Q}_{1},$ above
$\mathcal{T}_{r};$ the trajectories $\mathcal{T}_{-}$ start in $\mathcal{Q}%
_{4}$ under $\mathcal{T}_{r}.$ From Theorem \ref{tpsilon}, there exists a
unique trajectory $\mathcal{T}_{\varepsilon}$ converging to $(0,0)$ in
$\mathcal{Q}_{1}$ at $\infty,$ with the slope $1.$\medskip

(1) Case $\alpha<N.$ From Proposition \ref{zer}, all the solutions $w$ have at
most one simple zero.

The regular solutions stay positive, and $\mathcal{T}_{r}$ stays in its
quadrant, $\mathcal{Q}_{4}$ or $\mathcal{Q}_{1},$ from Remark \ref{vf} (see
figures 1 and 2). Then $\mathcal{T}_{\varepsilon}$ stays in $\mathcal{Q}_{1}$,
because it cannot meet $\mathcal{T}_{r}$ for $\alpha>0,$ or the line $\left\{
Y=0\right\}  $ for $\alpha<0,$ from Remark \ref{vf}; and the corresponding $w$
is of type (i).

Consider any trajectory $\mathcal{T}_{\left[  P\right]  }$ with $P\in
\mathcal{Q}_{1}$ above $\mathcal{T}_{\varepsilon}.$ It cannot stay in
$\mathcal{Q}_{1}$ because it does not meet $\mathcal{T}_{\varepsilon}$ and
converges to $(0,0)$ with a slope $0$. Thus it enters $\mathcal{Q}_{2}$ from
Remark \ref{vf}. Then $y$ has a unique zero, and $\mathcal{T}_{\left[
P\right]  }$ stays in $\mathcal{Q}_{1}$ before $P,$ and in $\mathcal{Q}%
_{2}\cup\mathcal{Q}_{3}$ after $P.$ Since $\mathcal{T}_{\left[  P\right]  }$
cannot meet $\pm\mathcal{T}_{\varepsilon},$ and $\lim_{\tau\rightarrow\infty
}\zeta=\alpha,$ $\mathcal{T}_{\left[  P\right]  }$ ends up in $\mathcal{Q}%
_{3}$ if $\alpha>0,$ in $\mathcal{Q}_{2}$ if $\alpha<0.$ It has the same
behaviour as $\mathcal{T}_{\varepsilon}$ at $-\infty,$ and $w$ is of type (iii).

Next consider $\mathcal{T}_{\left[  P\right]  }$ for any $P\in\mathcal{Q}%
_{1}\cup$ $\mathcal{Q}_{4}$ between $\mathcal{T}_{\varepsilon}$ and
$\mathcal{T}_{r}.$ Then $y$ stays positive, and $\mathcal{T}_{\left[
P\right]  }$ necessarily starts from $\mathcal{Q}_{1},$ and $w$ is of type (ii).

At least take any $P\in$ $\mathcal{Q}_{1}\cup\mathcal{Q}_{4}$ under
$\mathcal{T}_{r}.$ If $p\leqq N,$ $\mathcal{T}_{\left[  P\right]  }$ starts
from $\mathcal{Q}_{3}$ and $y$ has a unique zero, and $-w$ is of type (iii).
If $p>N,$ either $-w$ is of type (iii), or $\mathcal{T}_{\left[  P\right]  }$
stays in $\mathcal{Q}_{4}.$ From Theorems \ref{teta}, \ref{tzero}, either
$\mathcal{T}_{\left[  P\right]  }$ coincides with $\mathcal{T}_{u},$ and $w$
is of type (iv), or with one of the trajectories $\mathcal{T}_{-},$ thus $w$
is of type (v).\medskip

(2) Case $\alpha=N.$ All the solutions are given by (\ref{fat}), which is
equivalent to $J_{N}\equiv C,$ where $J_{N}$ is defined by (\ref{gg}). For
$C=0,$ the regular (Barenblatt) solutions, given by (\ref{fut}), are
nonnegative, with a compact support. In other words the trajectory
$\mathcal{T}_{\varepsilon}$ given by Theorem \ref{talpha} coincides with
$\mathcal{T}_{r},$ it is given by $y\equiv Y,$ $y>0$ (see figure 3). The only
change in the phase plane is the nonexistence of solutions of type
(ii).\medskip

(3) Case $\alpha>N.$

The regular solutions have a number $m\geqq1$ of simple zeros, from
Proposition \ref{zer} (see figure 4). As above, $\mathcal{T}_{r}$ starts from
$\mathcal{Q}_{1}$ with a finite slope $\alpha/N.$

Either $\mathcal{T}_{r}\neq\mathcal{T}_{\varepsilon}.$ Then the regular
solutions satisfy $\lim_{r\rightarrow\infty}r^{\alpha}w=$ $L\neq0.$ Since
$\mathcal{T}_{\varepsilon}$ cannot meet $\mathcal{T}_{r},$ $\mathcal{T}%
_{\varepsilon}$ also cuts the line $\left\{  y=0\right\}  $ at $m$ points, and
the corresponding $w$ is of type (vi). For any $P\in\mathcal{Q}_{1}$ above
$\mathcal{T}_{r},$ the trajectory $\mathcal{T}_{\left[  P\right]  }$ cuts the
line $\left\{  y=0\right\}  $ at $m+1$ points and $w$ is of type (vii). If
$p>N,$ there exist trajectories starting from $\mathcal{Q}_{1}$ between
$\mathcal{T}_{\varepsilon}$ and $\mathcal{T}_{r}$, with (\ref{via+}), such
that $w$ has $m$ simple zeros, and trajectories with (\ref{viu}) or
(\ref{via-}), $m$ zeros, and $\lim_{r\rightarrow\infty}r^{\alpha}w=$ $L\neq0.$

Or $\mathcal{T}_{r}=\mathcal{T}_{\varepsilon},$ the regular solutions have a
compact support, and we only find solutions of type (vii), (viii).
\end{proof}

\begin{remark}
In the case $\alpha=\eta<0$, the solutions (iv) are given by (\ref{aeta}). In
the case $N=1,$ $\alpha=-(p-1)/(p-2),$ the solutions of types (i) and (v) are
given by (\ref{exp}).
\end{remark}

\begin{remark}
We conjecture that there exists an increasing sequence $\left(  \bar{\alpha
}_{m}\right)  ,$ with $\bar{\alpha}_{0}=N$ such that the regular solutions $w$
have $m$ simple zeros for $\alpha\in\left(  \bar{\alpha}_{m-1},\bar{\alpha
}_{m}\right)  ,$ with $\lim_{r\rightarrow\infty}r^{\alpha}w=$ $L\neq0,$ and
$m$ simple zeros and a compact support for $\alpha=\bar{\alpha}_{m}$, in which
case $\mathcal{T}_{r}=\mathcal{T}_{\varepsilon}.$
\end{remark}

\section{The case $\varepsilon=-1,\alpha\leqq-\gamma$\label{two}}

\begin{theorem}
\label{osc}Assume $\varepsilon=-1,$ $\alpha\leqq-\gamma.$ Then all the
solutions $w$ of (\textbf{E}$_{w}$), among them the regular ones, are
ocillating near $\infty$ and $r^{-\gamma}w$ is asymptotically periodic in $\ln
r$. There exist$\medskip$

\noindent(i) solutions such that $r^{-\gamma}w$ is periodic in $\ln r;$

\noindent(ii) a unique solution with a hole;

\noindent(iii) flat solutions $w$ with (\ref{val}) or (\ref{kap}) near $0;$

\noindent(iv) solutions with (\ref{vet}) or(\ref{vpk}) or (\ref{via+}) or also
(\ref{via-}) near $0$;

\noindent(v) for $p>N,$ a unique solution with (\ref{viu}) near $0$.
\end{theorem}

\[%
\begin{array}
[c]{c}%
\raisebox{-0pt}{\includegraphics[
natheight=187.289200pt,
natwidth=187.289200pt,
height=189.096pt,
width=189.096pt
]%
{../../../AMarie/MesArticles/107050selfp+2/osc.bmp}%
}%
\\
\text{th \ref{osc},fig5: }\varepsilon=-1,N=1,p=3,\alpha=-4
\end{array}
\]

\begin{proof}
Here again, $(0,0)$ is the unique stationary point in the plane $(y,Y)$. Any
solution $y$ of (\textbf{E}$_{y}$) oscillates near $\infty,$ and $(y,Y)$ is
bounded from Proposition \ref{nmo}. From the strong form of the
Poincar\'{e}-Bendixon theorem, see \cite[p.239]{HuWe}, all the trajectories
have a limit cycle or are periodic. In particular $\mathcal{T}_{r}$ starts in
$\mathcal{Q}_{1}$, since $\varepsilon\alpha>0,$ with the asymptotical
direction $\varepsilon\alpha/N$. and it has a limit cycle $\mathcal{O}.$ There
exists a periodic trajectory of orbit $\mathcal{O},$ thus $w$ is of type (i)
(see figure 5).$\medskip$

From Theorem \ref{tzero} there exists a unique trajectory $\mathcal{T}%
_{\varepsilon}$ starting from $(0,0)$ with the slope $-1,$ $y>0;$ it has a
limit cycle $\mathcal{O}_{\varepsilon}\subset\mathcal{O},$ and $w$ is of type
(ii). For any $P$ in the bounded domain delimitated by $\mathcal{O}%
_{\varepsilon},$ not located on $\mathcal{T}_{\varepsilon},$ the trajectory
$\mathcal{T}_{\left[  P\right]  }$ does not meet $\mathcal{T}_{\varepsilon},$
and admits $\mathcal{O}_{\varepsilon}$ as limit cycle; near $-\infty,$ $y$ has
a constant sign, is monotone and converges to $(0,0)$ from Propositions
\ref{com} and \ref{esup}, and $\lim_{\tau\rightarrow-\infty}\zeta=\alpha$.
This show again the existence of such trajectories, proved at Theorem
\ref{teta}, and there is an infinity of them; and $w$ is if type
(iii)$.\medskip$

From Theorems \ref{teta} and \ref{tzero}, there exist trajectories starting
from infinity, with $\mathcal{O}$ as limit cycle, and $w$ is of type (iv) or
(v). If $\mathcal{O=O}_{\varepsilon},$ all the solutions are described.
\end{proof}

\section{Case $\varepsilon=1,\alpha<-\gamma.$ \label{three}}

\begin{theorem}
\label{mel}Assume $\varepsilon=1,$ $\alpha<-\gamma.$ Then $w\equiv\pm\ell
r^{\gamma}$ is a solution of (\textbf{E}$_{w}$). All regular solutions have a
strict constant sign, and satisfy (\ref{gam}) near $\infty.$ Moreover there
exist (exhaustively, up to a symmetry)$\medskip$

\noindent(i) a unique positive flat solution with (\ref{val}) near $0$ and
(\ref{gam}) near $\infty$;

\noindent(ii) a unique nonnegative solution with (\ref{vet}) or (\ref{vpk}) or
(\ref{via+}) near $0$, and compact support;

\noindent(iii) positive solutions with the same behaviour near 0 and
(\ref{gam}) near $\infty$;

\noindent(iv) solutions with one zero and the same behaviour near 0, and
$\left\vert w\right\vert $ satisfies (\ref{gam}) near $\infty;$

\noindent(v) for $p>N,$ positive solutions with (\ref{viu}) near 0 and
(\ref{gam}) near $\infty;$

\noindent(vi) for $p>N,$ positive solutions with (\ref{via-}) near 0 and
(\ref{gam}) near $\infty.$
\end{theorem}

\[%
\begin{array}
[c]{c}%
\raisebox{-0pt}{\includegraphics[
natheight=187.289200pt,
natwidth=187.289200pt,
height=189.096pt,
width=189.096pt
]%
{../../../AMarie/MesArticles/107050selfp+2/mel.bmp}%
}%
\\
\text{th \ref{mel}, fig6: }\varepsilon=1,N=2,p=3,\alpha=-6
\end{array}
\]

\begin{proof}
Here system (\textbf{S}) admits three stationary points in the plane $(y,Y)$,
given at (\ref{statio}), thus $w\equiv\pm\ell r^{\gamma}$ is a solution; and
$M_{\ell}$ is a sink (see figure 6). Any solution $y$ of (\textbf{E}$_{y}$)
has at most one zero, and is strictly monotone near $\pm\infty,$ from
Propositions \ref{zer} and \ref{com}.$\medskip$

From Theorems \ref{tpsilon} and \ref{talpha}, there exists a unique trajectory
$\mathcal{T}_{\varepsilon}$ converging to $(0,0)$ in $\mathcal{Q}_{1}$ at
$\infty,$ and a unique trajectory $\mathcal{T}_{\alpha}$ converging to $(0,0)$
in $\mathcal{Q}_{4}$ at $-\infty.$ The trajectory $\mathcal{T}_{r}$ starts in
$\mathcal{Q}_{4}$ with the asymptotical direction $-\left\vert \alpha
\right\vert /N$. From Remark \ref{vf}, $\mathcal{Q}_{4}$ is positively
invariant, and $\mathcal{Q}_{1}$ negatively invariant.\ Then $\mathcal{T}%
_{\varepsilon}$ stays in $\mathcal{Q}_{1},$ and $\mathcal{T}_{\alpha}$ and
$\mathcal{T}_{r}$ in $\mathcal{Q}_{4}$. From Proposition \ref{esup}, all the
trajectories, apart from $\pm\mathcal{T}_{\varepsilon},$ converge to $\pm
M_{\ell}$ at $\infty.$ Then $\mathcal{T}_{r}$ converges to $M_{\ell},$ and $w$
satisfies (\ref{gam}) near $\infty.$ And $\mathcal{T}_{\alpha}$ also converges
to $M_{\ell},$ and $w$ is of type (i).$\medskip$

From Propositions \ref{esup}, Theorems \ref{teta} and \ref{tzero}, all the
nonregular solutions which are positive near $-\infty$ satisfy (\ref{vet}),
(\ref{vpk}), (\ref{via+}), (\ref{via-}) or (\ref{viu})$,$ and there exist such
solutions. For $p<N$ (resp. $p=N),$ they correspond to trajectories
$\mathcal{T}_{\eta}$ (resp. $\mathcal{T}_{+})$ starting in $\mathcal{Q}_{1}$.
For $p>N,$ there is a unique trajectory $\mathcal{T}_{u}$ satisfying
(\ref{viu}), starting in $\mathcal{Q}_{4}$ under $\mathcal{T}_{r};$ and the
trajectories $\mathcal{T}_{+}$ satisfying (\ref{via+}) start from
$\mathcal{Q}_{1}$; the trajectories $\mathcal{T}_{-}$ satisfying (\ref{via-})
and the unique trajectory $\mathcal{T}_{u}$ satisfying (\ref{viu}) start from
$\mathcal{Q}_{4}$, under $\mathcal{T}_{r}.$ Since $\mathcal{T}_{\varepsilon}$
stays in $\mathcal{Q}_{1},$ it defines solutions $w$ of type (ii).$\medskip$

Consider the basis of eigenvectors $\left(  e_{1},e_{2}\right)  $ defined at
(\ref{bas}), where $\nu(\alpha)>0,$ associated to the eigenvalues $\lambda
_{1}<\lambda_{2}.$ One verifies that $\lambda_{1}<-\gamma<\lambda_{2};$ thus
$e_{1}$ points towards $\mathcal{Q}_{3}$ and $e_{2}$ points towards
$\mathcal{Q}_{4}.$ There exist unique trajectories $\mathcal{T}_{e_{1}}$ and
$\mathcal{T}_{-e_{1}}$ converging to $M_{\ell},$ tangentially to $e_{1}$ and
$-e_{1}$. All the other trajectories converging to $M_{\ell}$ at $\infty$ are
tangent to $\pm e_{2}.$ Let
\[
\mathcal{M}=\left\{  \left\vert Y\right\vert ^{(2-p)/(p-1)}Y=-\gamma
y\right\}  ,\qquad\mathcal{N}=\left\{  (N+\gamma)Y+\varepsilon\left\vert
Y\right\vert ^{(2-p)/(p-1)}Y=\varepsilon\alpha y\right\}
\]
\ be the sets of extremal points of $y$ and $Y.$ $\medskip$

The trajectory $\mathcal{T}_{r}$ starts above the curves $\mathcal{M}$ and
$\mathcal{N}$, thus $y^{\prime}<0$ and $Y^{\prime}>0$ near $-\infty.$ And
$\mathcal{T}_{r}$ converges to $M_{\ell}$ at $\infty,$ tangentially to
$e_{2}.$ Indeed if $\mathcal{T}_{r}=\mathcal{T}_{e_{1}},$ then $y$ has a
minimal point such that $y<\ell$ and $Y<-\left(  \gamma\ell\right)  ^{p-1},$
then $(y,Y)$ cannot be on $\mathcal{M}$. If $\mathcal{T}_{r}=\mathcal{T}%
_{-e_{1}},$ then $Y$ has a maximal point such that $y>\ell$ and $Y<-\left(
\gamma\ell\right)  ^{p-1},$ then also $(y,Y)$ cannot be on $\mathcal{N}$.
Finally $\mathcal{T}_{r}$ cannot end up tangentially to $-e_{2},$ it would
intersect $\mathcal{T}_{e_{1}}$ or $\mathcal{T}_{-e_{1}}.$ $\medskip$

The trajectory $\mathcal{T}_{\alpha}$ converge to $M_{\ell}$ tangentially to
$-e_{2}.$ Indeed if $\mathcal{T}_{\alpha}=\mathcal{T}_{e_{1}},$ then $Y$ has a
maximal point such that $y<\ell$ and $Y<-\left(  \gamma\ell\right)  ^{p-1};$
if $\mathcal{T}_{\alpha}=\mathcal{T}_{-e_{1}},$ then $y$ has a maximal point
such that $y>\ell$ and $Y>-\left(  \gamma\ell\right)  ^{p-1}.$ In any case we
reach a contradiction. Moreover $\mathcal{T}_{e_{1}}$ does not stay in
$\mathcal{Q}_{4}:$ $y$ would have a minimal point such that $y<\ell$ and
$Y<-\left(  \gamma\ell\right)  ^{p-1},$ which is impossible; thus
$\mathcal{T}_{e_{1}}$ starts in $\mathcal{Q}_{3},$ and enters $\mathcal{Q}%
_{4}$ at some point $\left(  \xi_{1},0\right)  $ with $\xi_{1}<0.$ And $-w$ is
of type (iv).$\medskip$

Any trajectory $\mathcal{T}_{\left[  P\right]  },$ with $P$ in the domain of
$\mathcal{Q}_{1}\cup\mathcal{Q}_{4}$ delimitated by $\mathcal{T}%
_{r},\mathcal{T}_{\alpha}$ and $\mathcal{T}_{\varepsilon},$ comes from
$\mathcal{Q}_{1},$ and converges to $M_{\ell}$ in $\mathcal{Q}_{4}$, in
particular $\mathcal{T}_{-e_{1}};$ the corresponding $w$ are of type
(iii).$\medskip$

Any trajectory $\mathcal{T}_{\left[  P\right]  },$ with $P$ in the domain of
$\mathcal{Q}_{3}\cup\mathcal{Q}_{4}$ delimitated by $\mathcal{T}_{e_{1}%
},\mathcal{T}_{\alpha}$ and $-\mathcal{T}_{\varepsilon},$ goes from
$\mathcal{Q}_{3}$ to $\mathcal{Q}_{4}$, and $\mathcal{T}_{\left[  P\right]  }$
converges to $M_{\ell}$ at $\infty,$ and $-w$ is of type (iv). For any
$\xi<\xi_{1},$ the trajectory $\mathcal{T}_{\left[  (0,\xi)\right]  }$ is of
the same type. If $p\leqq N,$ any trajectory in the domain under
$\mathcal{T}_{r},$ and $\mathcal{T}_{e_{1}}$ is of the same type. $\medskip$

If $p>N,$ moreover in this domain there exists a the unique trajectory
$\mathcal{T}_{u}$ and trajectories of the type $\mathcal{T}_{-}$ corresponding
to solutions $w$ of type (v) and (vi), from Theorems \ref{teta} and
\ref{tzero}. Up to a symmetry, all the solutions are described, and all of
them do exist.
\end{proof}

\section{Case $\varepsilon=-1,-\gamma<\alpha$\label{four}}

Here again System (S) admits the three stationary points (\ref{statio}), thus
$w\equiv\pm\ell r^{\gamma}$ is a solution of (\textbf{E}$_{w})$. The behaviour
is very rich: it depends on the position of $\alpha$ with respect to
$\alpha^{\ast}$ defined at (\ref{star}), and $0,$ $-p^{\prime},$ and $\eta$
(in case $p>N),$ and also $\alpha_{1},\alpha_{2}$ defined at (\ref{aun}). We
start from some general remarks$.$

\begin{remark}
\label{ploc1}(i) There exists a unique trajectory $\mathcal{T}_{\varepsilon}$
starting from $(0,0)$ in $\mathcal{Q}_{4}$ with the slope $-1$, from Theorem
\ref{tpsilon}.

\noindent(ii) There exists a unique trajectory $\mathcal{T}_{\alpha}$
converging to $(0,0)$ at $\infty,$ in $\mathcal{Q}_{1}$ if $\alpha>0,$ in
$\mathcal{Q}_{4}$ if $\alpha<0,$ with a slope $0$ at $(0,0),$ and $\lim
_{\tau\rightarrow\infty}\zeta=\alpha,$ from Theorem \ref{talpha}.

\noindent(iii) From Remark \ref{vf}, if $\alpha>0,$ $\mathcal{Q}_{4}$ is
positively invariant and $\mathcal{Q}_{1}$ negatively invariant. If
$\alpha<0,$ at any point $(0,\xi),\xi<0,$ the vector field points to
$\mathcal{Q}_{4},$ and at any point $(\varphi,0),\varphi>0,$ it points to
$\mathcal{Q}_{1}.$ Thus if $\mathcal{T}_{\varepsilon}$ does not stay in
$\mathcal{Q}_{1},$ then $\mathcal{T}_{\alpha}$ stays in the bounded domain
delimitated by $\mathcal{Q}_{4}\cap\mathcal{T}_{\varepsilon}$.\ If
$\mathcal{T}_{\alpha}$ does not stay in $\mathcal{Q}_{4},$ then $\mathcal{T}%
_{\varepsilon}$ stays in the bounded domain delimitated by $\mathcal{Q}%
_{4}\cap\mathcal{T}_{\alpha}$. If $\mathcal{T}_{\varepsilon}$ is homoclinic,
in other words $\mathcal{T}_{\varepsilon}=\mathcal{T}_{\alpha},$ it stays in
$\mathcal{Q}_{4}$.
\end{remark}

\begin{remark}
\label{ploc2} From Propositions \ref{esup}, Theorems \ref{teta} and
\ref{tzero}, all the nonregular solutions positive near $-\infty$ satisfy
(\ref{vet}) for $p<N$, (\ref{vpk}) for $p=N,$ corresponding to trajectories
$\mathcal{T}_{\eta},\mathcal{T}_{+}$ starting from$\mathcal{Q}_{1};$ and
(\ref{via+}), (\ref{via-}) or (\ref{viu}) for $p>N,$ corresponding to
trajectories $\mathcal{T}_{+}$ starting from $\mathcal{Q}_{1}$, and
$\mathcal{T}_{-},\mathcal{T}_{u}$ starting from $\mathcal{Q}_{4}.$
\end{remark}

\begin{remark}
\label{ploc3} Any trajectory $\mathcal{T}$ is bounded near $\infty$ from
Proposition \ref{nmo}. From the strong form of the Poincar\'{e}-Bendixon
theorem, any trajectory $\mathcal{T}$ bounded at $\pm\infty$ converges to
$(0,0)$ or $\pm M_{\ell},$ or its limit set $\Gamma_{\pm}$ at $\pm\infty$ is a
cycle, or it is homoclinic, namely $\mathcal{T}_{\varepsilon}$ $=\mathcal{T}%
_{\alpha}.$ If there exists a limit cycle surrounding $(0,0),$ it also
surrounds the points $\pm M_{\ell}$, from Proposition \ref{com}.
\end{remark}

The simplest case is $\alpha>0.$

\begin{theorem}
\label{int} Assume $\varepsilon=-1,$ $\alpha>0.$\medskip\ 

Then $w\equiv\ell r^{\gamma}$ is a solution $w$ of (\textbf{E}$_{w}$). All
regular solutions have a strict constant sign; and satisfy (\ref{gam}) near
$\infty.$ There exist (exhaustively, up to a symmetry)\medskip

\noindent(i) a unique nonnegative solution with a hole, and (\ref{gam}) near
$\infty$;

\noindent(ii) a unique positive solution with (\ref{vet}), or (\ref{vpk}) or
(\ref{via+})$,$ and (\ref{val}) near $\infty$;

\noindent(iii) positive solutions with the same behaviour near $0$, and
(\ref{gam}) near $\infty;$

\noindent(iv) solutions with one zero, the same behaviour near $0,$ and
$\left\vert w\right\vert $ satisfies (\ref{gam}) near $\infty;$

\noindent(v) for $p>N,$ a unique positive solution with (\ref{viu}) near $0,$
and (\ref{gam}) near $\infty$;

\noindent(vi) for $p>N,$ positive solutions with (\ref{via-}) near $0,$ and
(\ref{gam}) near $\infty.$
\end{theorem}

\[%
\begin{array}
[c]{cc}%
\raisebox{-0pt}{\includegraphics[
natheight=187.289200pt,
natwidth=187.289200pt,
height=189.096pt,
width=189.096pt
]%
{../../../AMarie/MesArticles/107050selfp+2/int1.bmp}%
}%
&
\raisebox{-0pt}{\includegraphics[
natheight=187.289200pt,
natwidth=187.289200pt,
height=189.096pt,
width=189.096pt
]%
{../../../AMarie/MesArticles/107050selfp+2/int2.bmp}%
}%
\\
\text{th \ref{int}, fig7: }\varepsilon=-1,N=1,p=3,\alpha=0.7 & \text{th
\ref{int}, fig8: }\varepsilon=-1,N=1,p=3,\alpha=1
\end{array}
\]

\begin{proof}
Any solution $y$ of (\textbf{E}$_{y}$) has at most one zero, and $y$ is
strictly monotone near $\infty,$ from Propositions \ref{zer} and \ref{cyc}.
The point $M_{\ell}$ is a sink and a node point, since $\alpha>0\geqq
\alpha_{2}$ (see figure 7). Consider the basis eigenvectors $(e_{1},$
$e_{2}),$ defined at (\ref{bas}), where $\nu(\alpha)<0,$ associated to the
eigenvalues $\lambda_{1}<\lambda_{2}<0.$ One verifies that $\lambda
_{1}<-\gamma<\lambda_{2},$ thus $e_{1}$ points towards $\mathcal{Q}_{3}$ and
$e_{2}$ points towards $\mathcal{Q}_{4}.$ There exist unique trajectories
$\mathcal{T}_{e_{1}}$ and $\mathcal{T}_{-e_{1}}$ tangent to $e_{1}$ and
$-e_{1}$ at $\infty.$ All the other trajectories which converge to $M_{\ell}$
end up tangentially to $\pm e_{1}.$ \medskip

The trajectory $\mathcal{T}_{\alpha}$ stays in $\mathcal{Q}_{1}$ from Remark
\ref{ploc1}; near $-\infty$ it is of type $\mathcal{T}_{\eta}$ for $p<N,$ and
$\mathcal{T}_{+}$ for $p\geqq N;$ it defines the solution of type (ii). Since
$\mathcal{T}_{\alpha}$ is the unique trajectory converging to $(0,0)$ at
$\infty,$ all the trajectories, apart from $\pm\mathcal{T}_{\alpha},$ converge
to $\pm M_{\ell}$ at $\infty,$ from Propositions \ref{com} and \ref{esup}%
.\medskip

The trajectories $\mathcal{T}_{r}$ and $\mathcal{T}_{\varepsilon}$ start in
$\mathcal{Q}_{4}$, and stay in it from Remark \ref{ploc1}, and both converge
to $M_{\ell}$ at $\infty,$ then $w$ satisfies (\ref{gam}); and $\mathcal{T}%
_{r}$ starts with the asymptotical direction $-\alpha/N$. And $\mathcal{T}%
_{\varepsilon}$ defines the solution of type (i). \medskip

As in the proof of Theorem \ref{mel}, $\mathcal{T}_{r}$ ends up tangentially
to $e_{2},$ and $\mathcal{T}_{\varepsilon}$ tangentially to $-e_{2}.$ Moreover
$\mathcal{T}_{e_{1}}$ does not stay in $\mathcal{Q}_{4},$ it starts in
$\mathcal{Q}_{3},$ and converges to $M_{\ell}$ in $\mathcal{Q}_{4},$ and $-w$
is of type (iv). Any trajectory $\mathcal{T}_{\left[  P\right]  },$ with $P$
in the domain of $\mathcal{Q}_{4}$ between $\mathcal{T}_{e_{1}},$
$\mathcal{T}_{\varepsilon},$ starts from $\mathcal{Q}_{3}$, enters
$\mathcal{Q}_{4}$ at some point $(0,\xi),\xi>\xi_{1},$ and has the same type
as $\mathcal{T}_{e_{1}}$. Any trajectory $\mathcal{T}_{\left[  (0,\xi)\right]
}$ with $\xi<\xi_{1}$ is of the same type.\medskip

Any trajectory $\mathcal{T}_{\left[  P\right]  },$ with $P$ in the domain of
$\mathcal{Q}_{1}\cup\mathcal{Q}_{4}$ above $\mathcal{T}_{r}\cup\mathcal{T}%
_{\varepsilon},$ starts from $\mathcal{Q}_{1},$ and converges to $M_{\ell}$ in
$\mathcal{Q}_{4}$, in particular $\mathcal{T}_{-e_{1}};$ the corresponding $w$
are of type (iii). If $p\leqq N,$ all the solutions are described. If $p>N$,
moreover there exist trajectories staying in $\mathcal{Q}_{4}:$ $\mathcal{T}%
_{u}$ and the $\mathcal{T}_{-}$, starting under $\mathcal{T}_{r},$
corresponding to types (v) and (vi).
\end{proof}

\begin{remark}
For $\alpha=N,$ $\mathcal{T}_{r}$ and $\mathcal{T}_{\varepsilon}$ are given by
(\ref{fut}), respectively with $K>0$ and $K<0.\ $ The trajectory
$\mathcal{T}_{\varepsilon}$ describes the portion $0M_{\ell}$ of the line
$\left\{  Y=-y\right\}  ,$ and $\mathcal{T}_{r}$ the complementary half-line
in $\mathcal{Q}_{4}$ (see figure 8).
\end{remark}

Next we assume $-p^{\prime}\leqq\alpha<0.$ The case $p>N$ is delicate: indeed
the special value $\alpha=\eta$ is involved, because $\eta<0.$

\begin{theorem}
\label{pom}Assume $\varepsilon=-1,p\leqq N,$ and $-p^{\prime}\leqq\alpha<0.$
Then $w\equiv\ell r^{\gamma}$ is a solution $w$ of (\textbf{E}$_{w}%
$).$\medskip$

There exist a unique nonnegative solution with a hole, satisfying (\ref{gam})
at $\infty.\medskip$

(1) If $\alpha\neq-p^{\prime},$ all regular solutions have one zero, and
$\left\vert w\right\vert $ satisfies (\ref{gam}) near $\infty.$ There exist
(exhaustively, up to a symmetry)$\medskip$

\noindent$\bullet$ for $p\leqq N,\medskip$

(i) a unique solution with one zero, with (\ref{vet}) or (\ref{vpk}) near
$0,$and (\ref{val}) near $\infty$;

(ii) solutions with one zero, with (\ref{vet}) or (\ref{vpk}) near $0$, and
$\left\vert w\right\vert $ satisfies (\ref{gam}) near $\infty;$

(iii) solutions with two zeros, with (\ref{vet}) or (\ref{vpk}) near $0,$ and
(\ref{gam}) near $\infty;\medskip$

\noindent$\bullet$ for $p>N,$ $\eta<\alpha,\medskip$

(iv) a unique positive solution, with (\ref{via-}) near $0,$ and (\ref{val})
near $\infty$;

(v) a unique positive solution, with (\ref{viu}) near $0,$ and (\ref{gam})
near $\infty$;

(vi) positive solutions, with (\ref{via-}) near $0,$ and (\ref{gam}) near
$\infty;$

(vii) solutions with one zero with (\ref{via-}) or (\ref{via+}) near $0,$ and
(\ref{gam}) near $\infty;\medskip$

\noindent$\bullet$ for $p>N,\alpha<\eta,\medskip$

(viii) a unique solution with one zero, with (\ref{via+}) near $0,$ and
(\ref{val}) near $\infty$;

(ix) a unique solution with one zero, with (\ref{viu}) near $0,$ and
$\left\vert w\right\vert $ satisfies (\ref{gam}) near $\infty$;

(x) solutions with one zero, with (\ref{via+}) or (\ref{via+}) near $0,$ and
$\left\vert w\right\vert $ satisfies (\ref{gam}) near $\infty;$

(xi) solutions with two zeros, with (\ref{via+}) near $0,$ and (\ref{gam})
near $\infty.\medskip$

\noindent$\bullet$ for $p>N,\alpha=\eta,$ solutions of the form
$w=cr^{\left\vert \eta\right\vert }$ $(c>0)$. The other solutions are of type
(vii).$\medskip$

(2) If $\alpha=-p^{\prime},$ all regular solutions have one zero and satisfy
(\ref{val}) near $\infty.$ The solutions without hole are of types (ii), (iii)
for $p\leqq N,$ (ix), (x), (xi) for $p>N.$
\end{theorem}

\[%
\begin{array}
[c]{cc}%
\raisebox{-0pt}{\includegraphics[
natheight=187.289200pt,
natwidth=187.289200pt,
height=189.096pt,
width=189.096pt
]%
{../../../AMarie/MesArticles/107050selfp+2/pom1.bmp}%
}%
&
\raisebox{-0pt}{\includegraphics[
natheight=187.289200pt,
natwidth=187.289200pt,
height=189.096pt,
width=189.096pt
]%
{../../../AMarie/MesArticles/107050selfp+2/pom2.bmp}%
}%
\\
\text{th \ref{pom}, fig9: }\varepsilon=-1,N=1,p=3,\alpha=-0.7 & \text{th
\ref{pom}, fig10: }\varepsilon=-1,N=1,p=3,\alpha=-1.49
\end{array}
\]%
\[%
\begin{array}
[c]{c}%
\raisebox{-0pt}{\includegraphics[
natheight=187.289200pt,
natwidth=187.289200pt,
height=189.096pt,
width=189.096pt
]%
{../../../AMarie/MesArticles/107050selfp+2/pom3.bmp}%
}%
\\
\text{th \ref{pom}, fig11: }\varepsilon=-1,N=1,p=3,\alpha=-3/2
\end{array}
\]

\begin{proof}
Here again $M_{\ell}$ is a sink; but it is a node point only if $\alpha
\geqq\alpha_{2.}$. The phase plane $(y,Y)$ does not contain any cycle, from
Proposition \ref{cyc}. From Proposition \ref{zer}, any solution $y$ has at
most two zeros, and $Y$ at most one.$\medskip$

The unique trajectory $\mathcal{T}_{\alpha}$ ends up in $\mathcal{Q}_{4}$ with
the slope $0.$ From the uniqueness of $\mathcal{T}_{\alpha}$ and
$\mathcal{T}_{\varepsilon},$ all the trajectories, apart from $\pm
\mathcal{T}_{\alpha},$ converge to $\pm M_{\ell}$ at $\infty,$ from
Proposition \ref{esup} and Remark \ref{ploc3}. Since $\varepsilon\alpha>0,$
the trajectory $\mathcal{T}_{r}$ starts in $\mathcal{Q}_{1},$ and $y$ has at
most one zero. Then $\mathcal{T}_{r}$ converges to $-M_{\ell}$ in
$\mathcal{Q}_{2}$, or $\mathcal{T}_{r}=-\mathcal{T}_{\alpha}$.$\medskip$

The trajectory $\mathcal{T}_{\varepsilon}$ starts in $\mathcal{Q}_{4}$ with
the slope $-1,$ satisfies $y\geqq0$ from Proposition \ref{zer}. If
$\mathcal{T}_{\varepsilon}$ converge to $(0,0),$ then $\mathcal{T}%
_{\varepsilon}=\mathcal{T}_{\alpha},$ thus it is homoclinic. Then $M_{\ell}$
is in the bounded component defined by $\mathcal{T}_{\varepsilon}$, and
$\mathcal{T}_{\varepsilon}$ meets $\mathcal{T}_{r},$ which is impossible.
Hence $\mathcal{T}_{\varepsilon}$ converges to $M_{\ell}$ in $\mathcal{Q}%
_{4},$ and $w$ is nonnegative with a hole and satisfies (\ref{gam}) near
$\infty.\medskip$

If $\alpha\neq-p^{\prime},$ we claim that $\mathcal{T}_{r}\neq-\mathcal{T}%
_{\alpha}.$ Indeed suppose $\mathcal{T}_{r}=-\mathcal{T}_{\alpha}.$ Consider
the functions $y_{\alpha},Y_{\alpha},$ defined by (\ref{cge}) with $d=\alpha.$
Then $Y_{\alpha}$ stays positive, and $Y_{\alpha}=O(e^{(\alpha\left(
p-1)+p\right)  \tau})$ at $\infty$, thus
\[
\lim_{\tau\rightarrow\infty}Y_{\alpha}=0,\quad\lim_{\tau\rightarrow\infty
}Y_{\alpha}=c>0,\quad\lim_{\tau\rightarrow-\infty}y_{\alpha}=\infty,\quad
\lim_{\tau\rightarrow\infty}y_{\alpha}=L<0.
\]
Moreover $y_{\alpha},Y_{\alpha}$ have no extremal point: at such a point, from
(\ref{yde}), (\ref{Yd}) the second derivatives have a strict constant sign;
then $Y_{\alpha}^{\prime}>0>y_{\alpha}^{\prime}.$ If $\alpha<\eta$ (in
particular if $p\leqq N),$ from (\ref{bb}), near $\infty,$
\[
(p-1)Y_{\alpha}^{\prime\prime}/Y_{\alpha}^{\prime}\geqq\left\vert Y\right\vert
^{(2-p)/(p-1)}(1+o(1)),
\]
thus $Y_{\alpha}^{\prime\prime}>0$ near $\infty,$ which is contradictory; if
$\alpha>\eta$, from (\ref{aa})
\[
(p-1)y_{\alpha}^{\prime\prime}/y_{\alpha}^{\prime}\geqq\left\vert Y\right\vert
^{(2-p)/(p-1)}(1+o(1)),
\]
thus $y_{\alpha}^{\prime\prime}<0$ near $\infty,$ still contradictory. If
$\alpha=\eta,$ $\mathcal{T}_{\alpha}=\mathcal{T}_{u}$ from (\ref{aeta}), thus
again $\mathcal{T}_{r}\neq-\mathcal{T}_{\alpha}.$\medskip

If $p>N$ and $\alpha\neq\eta,$ we claim that $\mathcal{T}_{\alpha}%
\neq\mathcal{T}_{u}.$ Indeed suppose $\mathcal{T}_{\alpha}=\mathcal{T}_{u}.$
This trajectory stays $\mathcal{Q}_{4},$ the function $\zeta$ stays negative,
and $\lim_{\tau\rightarrow-\infty}\zeta=\eta,$ $\lim_{\tau\rightarrow\infty}$
$\zeta=\alpha.$ If $\zeta$ has an extremal point $\vartheta,$ then
$\vartheta\in\left(  \alpha,\eta\right)  $ from System (\textbf{Q}), and
$\zeta^{\prime\prime}$ has a constant sign, the sign of $\alpha-\zeta$; it is
impossible. Thus $\zeta$ is monotone; then $(\alpha-\eta)\zeta^{\prime}>0,$
which contradicts System (\textbf{Q}).$\medskip$

(1) Case\textbf{ }$\alpha\neq-p^{\prime}.$ Since $\mathcal{T}_{r}%
\neq-\mathcal{T}_{\alpha},$ $\mathcal{T}_{r}$ converges to $-M_{\ell},$ and
$y$ has one zero, and $\left\vert w\right\vert $ satisfies (\ref{gam}).
$\medskip$

$\bullet$ Case $p\leqq N$. All the other trajectories start in $\mathcal{Q}%
_{3}$ or $\mathcal{Q}_{1},$ from Remarks \ref{ploc1} and \ref{ploc2}. For any
$\varphi>0,$ the trajectory $\mathcal{T}_{\left[  (\varphi,0)\right]  }$ goes
from $\mathcal{Q}_{4}$ into $\mathcal{Q}_{1}$, and converges to $-M_{\ell}$ in
$\mathcal{Q}_{2},$ since it cannot meet $\mathcal{T}_{r}$ and $-\mathcal{T}%
_{\varepsilon}$; thus $y$ has two zeros, and $w$ is of type (iii). The
trajectory $\mathcal{T}_{\alpha}$ cannot meet $\mathcal{T}_{\left[
(\varphi,0)\right]  },$ thus $y$ has one zero, and it has the same behaviour
at $-\infty$, and $w$ is of type (i). All the trajectories $\mathcal{T}%
_{\left[  P\right]  }$ with $P$ in the interior domain of $\mathcal{Q}_{1}$
delimitated by $-\mathcal{T}_{\varepsilon}$ and $\mathcal{T}_{r}$ start from
$\mathcal{Q}_{1}$ and converge to $-M_{\ell},$ $y$ has precisely one zero, and
has the same behaviour at $-\infty$, and $w$ is of type (ii).\medskip

$\bullet$ Case $p>N,$ $\eta<\alpha$ (see figure 9). Any solution $y$ has at
most one simple zero. The trajectory $\mathcal{T}_{\alpha}$ stays in
$\mathcal{Q}_{4}$. Indeed if it started in $\mathcal{Q}_{3},$ then for any
trajectory $\mathcal{T}_{\left[  (0,\xi)\right]  }$ with $(0,\xi)$ above
$-\mathcal{T}_{\alpha}$, the function $y$ would have two zeros. Since
$\mathcal{T}_{\alpha}\neq\mathcal{T}_{u},$ we have $\mathcal{T}_{\alpha}%
\in\mathcal{T}_{-},$ and $w$ is of type (iv). The trajectory $\mathcal{T}_{u}$
necessarily stays in $\mathcal{Q}_{4}$ and converges to $M_{\ell},$ and $w$ is
of type (v). The trajectories $\mathcal{T}_{\left[  P\right]  },$ with $P$ in
the domain delimitated by $\mathcal{T}_{u},\mathcal{T}_{\alpha}$ and
$\mathcal{T}_{\varepsilon},$ are of type $\mathcal{T}_{-}$ and converge in
$\mathcal{Q}_{4}$ to $M_{\ell},$ and $w$ is of type (vi). The trajectories
$\mathcal{T}_{\left[  P\right]  },$ with $P$ in the domain delimitated by
$\mathcal{T}_{r},\mathcal{T}_{\alpha}$ and $-\mathcal{T}_{\varepsilon},$ are
of type $\mathcal{T}_{-}$, and converge to $-M_{\ell},$ and $y$ has one zero.
The trajectories $\mathcal{T}_{\left[  P\right]  },$ with $P$ in the domain
delimitated by $\mathcal{T}_{r}$ and $-\mathcal{T}_{u},$ are of type
$\mathcal{T}_{+,}$converge to $-M_{\ell},$ and $y$ has one zero. Both define
solutions $w$ of type (vii).\medskip

$\bullet$ Case $p>N,\alpha<\eta$ (see figure 10). We have seen that
$\mathcal{T}_{r}\neq-\mathcal{T}_{\alpha}.$ If $\mathcal{T}_{\alpha}%
\in\mathcal{T}_{+},$ then $\zeta$ decreases from $0$ to $\alpha,$ which
contradicts System (\textbf{Q}) at $\infty.$ Then $\mathcal{T}_{\alpha}$ does
not stay in $\mathcal{Q}_{4},$ it starts in $\mathcal{Q}_{3}$ and
$-\mathcal{T}_{\alpha}\in\mathcal{T}_{-},$ hence $y$ has a zero, and $w$ is of
type (viii). Then $\mathcal{T}_{u}$ and the trajectories $\mathcal{T}_{-}$
converge to $-M_{\ell},$ and $y$ has one zero. The trajectories $\mathcal{T}%
_{\left[  P\right]  },$ with $P$ in the domain delimitated by $\mathcal{T}%
_{r},-\mathcal{T}_{\alpha}$ and $-\mathcal{T}_{\varepsilon},$ are of type
$\mathcal{T}_{+}$ and converge to $-M_{\ell},$ $y$ has one zero. They
correspond to $w$ is of type (ix) or (x). The trajectories $\mathcal{T}%
_{\left[  P\right]  },$ with $P$ in $\mathcal{Q}_{4}$ above $\mathcal{T}_{r},$
cut the line $\left\{  y=0\right\}  $ twice, and converge to $M_{\ell},$ and
$w$ is of type (xi).\medskip

$\bullet$ Case $p>N,\alpha=\eta.$ Then $\mathcal{T}_{\alpha}=\mathcal{T}_{u}$,
the functions $w=cr^{-\eta}$ $(c>0)$ are particular solutions$.$ The phase
plane study is the same, and gives only solutions of type (vii).\medskip

(2) Case\textbf{ } $\alpha=-p^{\prime}$ (see figure 11). Here $\mathcal{T}%
_{r}=-\mathcal{T}_{\alpha},$ since the regular solutions are given by
(\ref{moi})$.$ Thus there exist no more solutions of type (ii) or
(viii).\medskip
\end{proof}

Next we study the behaviour of all the solutions when $\alpha<-p^{\prime}.$ In
particular we prove the existence and uniqueness of an $\alpha_{c}$ for which
there exists an homoclinic trajectory. Thus we find again some results
obtained in \cite{GiVa}, with new detailed proofs. We also improve the bounds
for $\alpha_{c},$ in particular $\alpha^{\ast}<\alpha_{c}.$

\begin{lemma}
\label{alphap}Let
\[
\alpha_{p}:=-(p-1)/(p-2).
\]
If $N=1,$ for $\alpha=\alpha_{p},$ then there exists an homoclinic trajectory
in the phase plane $\left(  y,Y\right)  .$ If $N\geqq2,$ for $\alpha
=\alpha_{p},$ there is no homoclinic trajectory, moreover $\mathcal{T}%
_{\alpha}$ converges to $M_{\ell}$ at $-\infty$ or has a limit cycle in
$\mathcal{Q}_{4}$.
\end{lemma}

\begin{proof}
In the case $N=1,$ $\alpha=\alpha_{p},$ the explicit solutions (\ref{exp})
define an homoclinic trajectory in the phase plane $(y,Y)$, namely
$\mathcal{T}_{\varepsilon}=\mathcal{T}_{\alpha}.$ In the phase plane $(g,s)$
of System \textbf{(R}), from Remark \ref{parti}, they correspond to the line
$s\equiv1+\alpha g,$ joining the stationary points $(0,1)$ and ($-1/\alpha
,0).\medskip$

Next assume $N\geqq2$ and consider the trajectory $\mathcal{T}_{\alpha}$ in
the plane $(y,Y).$ In the plane $(g,s)$ of System (\textbf{R}), the
corresponding trajectory $\mathcal{T}_{\alpha}^{\prime}$ ends up at
($-1/\alpha,0),$ as $\nu$ tends to $\infty$ from (\ref{gsnu}), with the slope
$-k_{p}.$ If $\mathcal{T}_{\alpha}$ is homoclinic, then $\mathcal{T}_{\alpha
}^{\prime}$ converges to $(0,1)$ as $\nu$ tends to $-\infty.$ Consider the
segment
\[
T=\left\{  (g,-k(g+1/\alpha_{p}):g\in\left[  0,1/\left\vert \alpha
_{p}\right\vert \right]  \right\}  ,\text{ \quad with\quad\ }k=p^{\prime
}\alpha_{p}^{2}/(N+2/(p-2))>k_{p}.
\]
Its extremity $(0,k/\left\vert \alpha_{p}\right\vert )$ is strictly under
$(0,1)$. The domain $\mathcal{R}$ delimitated by the axes, which are
particular orbits, and $T,$ is negatively invariant: indeed, at any point of
$T,$ we find
\[
k\frac{dg}{d\nu}+\frac{ds}{d\nu}=(N-1)p^{\prime}ks(g-\frac{1}{\gamma})^{2}.
\]
The trajectory $\mathcal{T}_{\alpha}^{\prime}$ ends up in $\mathcal{R},$ thus
it stays in it, hence $\mathcal{T}_{\alpha}^{\prime}$ cannot join $(0,1).$ In
the phase plane $(y,Y),$ $\mathcal{T}_{\alpha}$ is not homoclinic, and
$\mathcal{T}_{\alpha}$ stays in $\mathcal{Q}_{4},$ and Remark \ref{ploc3} applies.
\end{proof}

\begin{remark}
Notice that $\alpha^{\ast}\leqq\alpha_{p}\Leftrightarrow N\leqq p.$
\end{remark}

\begin{theorem}
\label{clin}Assume $\varepsilon=-1,$ and $\alpha<-p^{\prime}.$ There exists a
unique $\alpha_{c}<0$ such that there exists an homoclinic trajectory in the
plane $\left(  y,Y\right)  ;$ in other words $\mathcal{T}_{\varepsilon
}=\mathcal{T}_{\alpha}.$ If $N=1,$ then $\alpha_{c}=\alpha_{p}.$ If $N\geqq2,$
then
\begin{equation}
max(\alpha^{\ast},\alpha_{p})<\alpha_{c}<\min(\alpha_{2},-p^{\prime}).
\label{fic}%
\end{equation}

\end{theorem}

\begin{proof}
In order to prove the existence of an homoclinic orbit for System
(\textbf{S}), we could consider a Poincar\'{e} application as in \cite{Bi1},
but it does not give uniqueness. Thus we consider the system \textbf{(R}%
$_{\beta}$) obtained from \textbf{(R}) by setting $s=\beta S$:
\[
\left.
\begin{array}
[c]{c}%
\qquad\quad\frac{dg}{d\nu}=gF(g,S),\qquad\qquad F(g,S):=\beta S(1+\eta
g)-\frac{1}{p-1}(1+\alpha g),\\
\\
\qquad\quad\frac{dS}{d\nu}=SG(g,S),\qquad\qquad G(g,S):=1+\alpha
g-\beta(1+Ng)S.\qquad
\end{array}
\right\}  \qquad\qquad\quad\text{\textbf{(R}}_{\beta}\text{)}%
\]
Its stationary points are
\[
(0,0),\qquad A^{\prime}=(1/\left\vert \alpha\right\vert ,0),\qquad B^{\prime
}=(0,1/\beta),\qquad M^{\prime}=(1/\gamma,1/(N+\gamma)(p-2)),
\]
where $M^{\prime}$ corresponds to $M_{\ell}.$ The existence of homoclinic
trajectory for System (\textbf{S}) resumes to the existence of a trajectory
for System \textbf{(R}$_{\beta}$) in the plane $(g,S)$, starting from
$B^{\prime}$ and ending at $A^{\prime}.\medskip$

(i) \textit{Existence}. We can assume that $\alpha\in\left(  \alpha_{1}%
,\min(\alpha_{2},-p^{\prime})\right)  $, from Proposition \ref{cyc}. In the
plane $(g,S),$ consider the trajectories $\mathcal{T}_{\varepsilon}^{\prime}$
and $\mathcal{T}_{\alpha}^{\prime}$ corresponding to $\mathcal{T}%
_{\varepsilon}\cap$ $\mathcal{Q}_{4}$ and $\mathcal{T}_{\alpha}\cap
\mathcal{Q}_{4}$ in the plane $(y,Y).$ Then $\mathcal{T}_{\varepsilon}%
^{\prime}$ starts from $B^{\prime}$ and $\mathcal{T}_{\alpha}^{\prime}$ ends
up at $A^{\prime}.$ From Remark \ref{ploc1}, for any $\alpha\in\left(
\alpha_{1},\alpha_{2}\right)  ,$ with $\alpha\leqq-p^{\prime},$ we have three
possibilities:$\medskip$

$\bullet$ $\mathcal{T}_{\varepsilon}^{\prime}$ is converging to $M^{\prime}$
as $\nu$ tends to $\infty$ and turns around this point, since $\alpha$ is a
spiral point, or it has a limit cycle in $\mathcal{Q}_{1}$ around $M^{\prime
}.$ And $\mathcal{T}_{\alpha}^{\prime}$ admits the line $g=0$ as an asymptote
as $\nu$ tends to $-\infty,$ which means that $\mathcal{T}_{\alpha}$ does not
stay in $\mathcal{Q}_{4}$ in the plane $(y,Y).$ Then $\mathcal{T}%
_{\varepsilon}^{\prime}$ meats the line
\[
L:=\left\{  g=1/\gamma\right\}
\]
at a first point $(1/\gamma,S_{0}(\alpha)).$ And $\mathcal{T}_{\alpha}%
^{\prime}$ meats $L$ at a last point $(1/\gamma,S_{1}(\alpha)),$ such that
$S_{0}(\alpha)-S_{1}(\alpha)<0$;\medskip\ 

$\bullet$ $\mathcal{T}_{\alpha}^{\prime}$ is converging to $M^{\prime}$ at
$-\infty$ or it has a limit cycle in $\mathcal{Q}_{1}$ around $M^{\prime}.$
And $\mathcal{T}_{\varepsilon}^{\prime}$ admits the line $S=0$ as an asymptote
at $\infty,$ which means that $\mathcal{T}_{\varepsilon}$ does not stay in
$\mathcal{Q}_{4}$. Then with the same notations, $S_{0}(\alpha)-S_{1}%
(\alpha)>0.\medskip$

$\bullet$ $\mathcal{T}_{\varepsilon}^{\prime}=\mathcal{T}_{\alpha}^{\prime},$
equivalently $S_{0}(\alpha)-S_{1}(\alpha)=0.\medskip$

The function $\alpha\mapsto\varphi(\alpha)=S_{0}(\alpha)-S_{1}(\alpha)$ is
continuous, from Theorems \ref{tpsilon} and \ref{talpha}. If $-p^{\prime
}<\alpha_{2},$ then $\varphi(-p^{\prime})$ is well defined and $\varphi
(-p^{\prime})<0;$ indeed $\mathcal{T}_{\alpha}=-\mathcal{T}_{r},$ thus
$\mathcal{T}_{\alpha}$ does not stay in $\mathcal{Q}_{4}$ from Theorem
\ref{pom}. If $\alpha_{2}\leqq-p^{\prime},$ in the plane $(y,Y),$ the
trajectory $\mathcal{T}_{\alpha_{2}}$ leaves $\mathcal{Q}_{4},$ from
Proposition \ref{cyc}, because $\alpha_{2}$ is a sink, and transversally from
Remark \ref{ploc1}. The same happens for $\mathcal{T}_{\alpha_{2-\upsilon}}$
for $\upsilon>0$ small enough, by continuity, thus $\varphi(\alpha
_{2}-\upsilon)<0.$ From Lemma \ref{alphap}, $\varphi(\alpha_{p})>0$ if
$N\geqq2,$ and $\varphi(\alpha_{p})=0$ if $N=1.$ In any case there exists at
least an $\alpha_{c}$ satisfying (\ref{fic}), such that $\varphi(\alpha
_{c})=0.\medskip$

(ii) \textit{Uniqueness}. First observe that $1+\eta g>0;$ indeed
$1+\eta/\left\vert \alpha\right\vert >(p^{\prime}+\eta)/\left\vert
\alpha\right\vert >0.$ Now
\[
(p-1)F+G=p\beta S(1/\gamma-g)=(p-2)\beta S(1-\gamma g),
\]
hence the curves $\left\{  F=0\right\}  $ and $\left\{  G=0\right\}  $
intersect at $M^{\prime}$ and $A^{\prime},$ $\left\{  G=0\right\}  $ contains
$B^{\prime}$ and is above $\left\{  F=0\right\}  $ for $g\in$ $\left(
0,1/\gamma\right)  $ and under it for $g\in\left(  1/\gamma,1/\left\vert
\alpha\right\vert \right)  .$ Moreover $\mathcal{T}_{\varepsilon}^{\prime}$
has a negative slope at $B^{\prime},$ thus $F>0>G$ near $0$ from
(\textbf{R}$_{\beta}$). And $\mathcal{T}_{\varepsilon}^{\prime}$ cannot meet
$\left\{  G=0\right\}  $ for $\left(  0,1/\gamma\right)  ,$ because on this
curve the vector field is $(gF,0)$ and $F>0.$ Thus $\mathcal{T}_{\varepsilon
}^{\prime}$ satisfies $F>0>G$ on $\left(  0,1/\gamma\right)  $. In the same
way $\mathcal{T}_{\alpha}^{\prime}$ has a negative slope $-\theta\alpha
^{2}/(p-1)(\eta+\left\vert \alpha\right\vert )<0$ at $1/\left\vert
\alpha\right\vert ,$ thus $F>0>G$ near $1/\left\vert \alpha\right\vert .$ And
$\mathcal{T}_{\alpha}^{\prime}$ cannot meet $\left\{  F=0\right\}  ,$ because
the vector field on this curve is $(0,SG)$ and $G<0$. Thus $\mathcal{T}%
_{\alpha}^{\prime}$ satisfies $F>0>G$ on $\left(  1/\gamma,1/\left\vert
\alpha\right\vert \right)  .$ \medskip\ 

Let $\alpha<\bar{\alpha}$. Then $\mathcal{T}_{\varepsilon}^{\prime}$ is above
$\mathcal{\bar{T}}_{\varepsilon}^{\prime}$ near $g=0$, and $\mathcal{T}%
_{\alpha}^{\prime}$ is at the left of $\mathcal{T}_{\bar{\alpha}}^{\prime}$
near $S=0.$ We show that $\varphi(\alpha)>\varphi(\bar{\alpha}).$ First
suppose that $\mathcal{T}_{\varepsilon}^{\prime}$ and $\mathcal{\bar{T}%
}_{\varepsilon}^{\prime}$ (or $\mathcal{T}_{\alpha}^{\prime}$ and
$\mathcal{\bar{T}}_{\bar{\alpha}}^{\prime})$ intersect at a first point
$P_{1}$ (or a last point) such $g\neq$ $1/\gamma.$ Then at this point
\begin{equation}
\frac{1}{p-1}\frac{g}{S}\frac{dS}{dg}+1=\frac{(p-2)(1-\gamma g)S}%
{(p-1)S(1+\eta g)-\beta^{-1}(1+\alpha g)}=\frac{(p-2)(1-\gamma g)S}%
{h_{S}(g)-\beta^{-1}(1-\gamma g)} \label{ttt}%
\end{equation}
with $h_{S}(g)=(p-1)S(1+\eta g)-g/(p-2).$ Thus the denominator, which is
positive, is increasing in $\alpha$ on $\left(  0,1/\gamma\right)  ,$
decreasing on $\left(  1/\gamma,1/\left\vert \alpha\right\vert \right)  ;$ in
any case $dS/dg>d\overline{S}/dg$ at $P_{1},$ which is contradictory. Next
suppose that there is an intersection on $L.$ At such a point $P_{1}%
=(1/\gamma,S_{1})=(1/\gamma,\overline{S}_{1})$ the derivatives are equal from
(\ref{ttt}), and $P_{1}$ is above $M^{\prime}$, because $F>0.$ At any points
$(g,S(g))\in\mathcal{T}_{\varepsilon}^{\prime}$ (or $\mathcal{T}_{\alpha
}^{\prime}),$ $(g,\overline{S}(g))\in\mathcal{\bar{T}}_{\varepsilon}^{\prime}$
(or $\mathcal{\bar{T}}_{\bar{\alpha}}^{\prime}),$ setting $g=1/\gamma+u,$
\[
\Phi(u)=(\frac{1}{p-1}\frac{g}{S}\frac{dS}{dg}+1)\frac{1}{(p-2)S}%
=-\frac{\gamma}{h_{S}(1/\gamma)}u+\frac{1}{h_{S}^{2}(1/\gamma)}(\frac{\gamma
}{\beta}+h_{S}^{\prime}(1/\gamma))u^{2}(1+o(1)),
\]%
\[
\bar{\Phi}(u)=(\frac{1}{p-1}\frac{g}{\overline{S}}\frac{d\overline{S}}%
{dg}+1)\frac{1}{(p-2)\overline{S}}=-\frac{\gamma}{h_{\overline{S}}(1/\gamma
)}u+\frac{1}{h_{\overline{S}}^{2}(1/\gamma)}(\frac{\gamma}{\beta}%
+h_{\overline{S}}^{\prime}(1/\gamma))u^{2}(1+o(1)),
\]
And $h_{S}(1/\gamma)=h_{\overline{S}}(1/\gamma)>0$, and $h_{S}^{\prime
}(1/\gamma)=h_{\overline{S}}^{\prime}(1/\gamma),$ then%
\[
(\Phi-\bar{\Phi})(u)=\frac{\gamma u^{2}(1/\beta-1/\bar{\beta})}{h(1/\gamma
)}(1+o(1)).
\]
This implies $d^{2}(S-\overline{S})/dg^{2}=0$ and $d^{3}(S-\overline
{S})/dg^{3}=2S_{1}\gamma^{2}(p-1)(p-2)(1/\beta-1/\bar{\beta})>0,$ which is a
contradiction. Then $\mathcal{T}_{\varepsilon}^{\prime}$ and $\mathcal{\bar
{T}}_{\varepsilon}^{\prime}$ cannot intersect on this line, similarly for
$\mathcal{T}_{\alpha}^{\prime}$ and $\mathcal{\bar{T}}_{\bar{\alpha}}^{\prime
}.$ Hence $\varphi(\alpha)>\varphi(\bar{\alpha}),$ which proves the
uniqueness. \medskip

As a consequence, for $\alpha<\alpha_{c},$ $\varphi(\alpha)>0,$ in the plane
$(y,Y),$ $\mathcal{T}_{\varepsilon}$ does not stay in $\mathcal{Q}_{4};$ for
$\alpha>\alpha_{c},$ $\varphi(\alpha)<0,$ $\mathcal{T}_{\alpha}$ does not stay
in $\mathcal{Q}_{4}.$ From Lemma \ref{alphap}, it follows that $\alpha
_{p}<\alpha_{c}$ if $N\geqq2.$ Moreover $\alpha^{\ast}<\alpha_{c}.$ Indeed
$\alpha^{\ast}$ is a weak source from Proposition \ref{ws}, thus for
$\alpha>\alpha^{\ast}$ small enough, there exists a unique cycle $\mathcal{O}$
around $M_{\ell},$ which is unstable. For such an $\alpha,$ $\mathcal{T}%
_{\varepsilon}$ cannot stay in $\mathcal{Q}_{4}:$ it would have $\mathcal{O}$
as a limit cycle at $\infty$, which contradicts the unstability.\medskip
\end{proof}

Next we discuss according to the position of $\alpha$ with respect to
$\alpha^{\ast}$ and $\alpha_{c}.$

\begin{theorem}
\label{sou}Assume $\varepsilon=-1,$ and $\alpha\leqq\alpha^{\ast}.$
Then$\medskip$

\noindent(i) there exist a unique flat positive solution $w$ of (\textbf{E}%
$_{w}$) with (\ref{gam}) near $0,$ and (\ref{val}) near $\infty;$

\noindent(ii) All the other solutions are oscillating at $\infty,$ among them
the regular ones, and $r^{-\gamma}w$ is asymptotically periodic in $\ln r$.
There exist solutions with a hole, also with (\ref{gam}), (\ref{vet}) or
(\ref{via+}) or (\ref{via+}) or (\ref{viu}) near $0.$ There exist solutions
such that $r^{-\gamma}w$ is periodic in $\ln r.$
\end{theorem}

\[%
\begin{array}
[c]{cc}%
\raisebox{-0pt}{\includegraphics[
natheight=196.946396pt,
natwidth=196.946396pt,
height=198.691pt,
width=198.691pt
]%
{../../../AMarie/MesArticles/107050selfp+2/sou1.bmp}%
}%
&
\raisebox{-0pt}{\includegraphics[
natheight=187.289200pt,
natwidth=187.289200pt,
height=189.096pt,
width=189.096pt
]%
{../../../AMarie/MesArticles/107050selfp+2/sou2.bmp}%
}%
\\
\text{th \ref{sou},fig 12: }\varepsilon=-1,N=1,p=3,\alpha=-2.53 & \text{th
\ref{sou}, fig 13: }\varepsilon=-1,N=1,p=3,\alpha=-2.2
\end{array}
\]

\begin{proof}
Here $\alpha<\alpha_{c},$ from Theorem \ref{clin}, and the trajectory
$\mathcal{T}_{\alpha}$ stays in $\mathcal{Q}_{4}$. From Proposition \ref{cyc},
it converges at $-\infty$ to $M_{\ell},$ and $w$ is of type (i). \medskip

The trajectory $\mathcal{T}_{\varepsilon}$ leaves $\mathcal{Q}_{4}$, and
cannot converge either to $(0,0)$ since $\mathcal{T}_{\varepsilon}%
\neq\mathcal{T}_{\alpha}$, or to $\pm M_{\ell},$ because this point is a
source, or a weak source. Recall that $M_{\ell}$ is a node point for
$\alpha\leqq\alpha_{1}$ (see figure 12,, where $\alpha_{1}\cong-2.50$), or a
spiral point (see figure 13). And $\mathcal{T}_{\varepsilon}$ is bounded at
$\infty$ from Proposition \ref{nmo}. Then it has a limit cycle $\mathcal{O}%
_{\varepsilon}$ surrounding $(0,0)$ from Proposition \ref{cyc}, and $\pm
M_{\ell}$ from Remark \ref{ploc3}. Thus $w$ is oscillating around $0$ near
$\infty$, $r^{-\gamma}w$ is asymptotically periodic in $\ln r.$ \medskip

The solutions $w$ corresponding to $\mathcal{O}_{\varepsilon}$ are oscillating
and $r^{-\gamma}w$ is periodic in $\ln r.$ Any trajectory $\mathcal{T}%
_{\left[  P\right]  }$ with $P$ in the interior domain delimitated by
$\mathcal{O}_{\varepsilon}$ converges to $M_{\ell}$ at $-\infty$ and has the
same limit cycle at $\infty.$ The trajectory $\mathcal{T}_{r}$ starts in
$\mathcal{Q}_{1},$ with $\lim_{\tau\rightarrow-\infty}y=\infty$ and cannot
converge to any stationary point at $\infty$. It is bounded, thus has a limit
cycle $\mathcal{O}_{r}$ surrounding $\mathcal{O}_{0}$. For any $P\not \in
\mathcal{T}_{r}$ in the exterior domain to $\mathcal{O}_{r}$, the trajectory
$\mathcal{T}_{\left[  P\right]  }$ admits $\mathcal{O}_{r}$ as a limit cycle
at $\infty,$ and $y$ is necessarily monotone at $-\infty,$ thus (\ref{vet}) or
(\ref{via+}) or (\ref{via+}) or (\ref{viu}) near $0;$ all those solutions
exist. The question of the uniqueness of the cycle ($\mathcal{O}_{r}$
$=\mathcal{O}_{\varepsilon})$ is open.\medskip

\begin{theorem}
\label{orb}Let $\alpha_{c}$ be defined by Theorem \ref{clin}.\medskip\ 

(1) Let $\alpha^{\ast}<\alpha<\alpha_{c}.$ Then all regular solutions $w$ of
(\textbf{E}$_{w}$) are oscillating around $0$ near $\infty,$ and $r^{-\gamma
}w$ is asymptotically periodic in $\ln r$. There exist$\medskip$

\noindent(i) \textbf{positive} solutions, such that $r^{-\gamma}w$ is periodic
in $\ln r;$

\noindent(ii) a unique positive solution such that $r^{-\gamma}w$ is
asymptotically periodic in $\ln r$ near $0$, with (\ref{val}) near $\infty$;

\noindent(iii) positive solutions such that $r^{-\gamma}w$ is asymptotically
periodic in $\ln r$ near $0$, with (\ref{gam}) near $\infty;$

\noindent(iv) solutions oscillating around $0$ such that $r^{-\gamma}w$ is
periodic in $\ln r;$

\noindent(v) solutions with a hole, oscillating near $\infty,$ such that
$r^{-\gamma}w$ is asymptotically periodic in $\ln r;$

\noindent(vi) solutions satisfying (\ref{vet}) or (\ref{via+}) or (\ref{via+})
or (\ref{viu}) near $0$, oscillating around $0$ near $\infty,$ such that
$r^{-\gamma}w$ is asymptotically periodic in $\ln r;$

\noindent(vii) solutions positive near $0,$ oscillating near $\infty,$ such
that $r^{-\gamma}w$ is asymptotically periodic in $\ln r$ near $0$ and
$\infty.\medskip$

(2) Let $\alpha=\alpha_{c}.\medskip$

\noindent(viii) There exist a \textbf{unique nonnegative} solution
\textbf{with a hole}, with (\ref{val}) near $\infty$.

\noindent The regular solutions are as above. There exist solutions of types
(iv), (vi), and

\noindent(ix) positive solutions such that $r^{-\gamma}w$ is bounded from
above near $0$, with (\ref{gam}) near $\infty$.
\end{theorem}
\end{proof}

\[%
\begin{array}
[c]{cc}%
\raisebox{-0pt}{\includegraphics[
natheight=187.289200pt,
natwidth=187.289200pt,
height=189.096pt,
width=189.096pt
]%
{../../../AMarie/MesArticles/107050selfp+2/orb1.bmp}%
}%
&
\raisebox{-0pt}{\includegraphics[
natheight=187.289200pt,
natwidth=187.289200pt,
height=189.096pt,
width=189.096pt
]%
{../../../AMarie/MesArticles/107050selfp+2/orb2.bmp}%
}%
\\
\text{th \ref{orb},fig 14: }\varepsilon=-1,N=1,p=3,\alpha=-2.1 & \text{th
\ref{orb},fig 15: }\varepsilon=-1,N=1,p=3,\alpha=-2
\end{array}
\]

\begin{proof}
(1) Let $\alpha^{\ast}<\alpha<\alpha_{c}$ (see figure 14). Then $\mathcal{T}%
_{\alpha}$ stays in $\mathcal{Q}_{4},$ but cannot converge neither to
$M_{\ell}$ which is a sink, nor to $(0,0)$ since $\mathcal{T}_{\alpha}%
\neq\mathcal{T}_{\varepsilon}.$ It has a limit cycle $\mathcal{O}_{\alpha}$ in
$\mathcal{Q}_{4}$ at $-\infty,$ surrounding $M_{\ell},$ and $w$ is of type
(ii). The orbit $\mathcal{O}_{\alpha}$ corresponds to solutions of type (i).
There exist positive solutions converging to $M_{\ell}$ at $\infty,$ with a
limit cycle $\mathcal{O}_{\ell}$ at $-\infty$ surrounded by $\mathcal{O}%
_{\alpha},$ and $w$ is of type (iii). This cycle is unique ($\mathcal{O}%
_{\ell}=$ $\mathcal{O}_{\alpha})$ for $\alpha-\alpha^{\ast}$ small enough,
from Proposition \ref{ws}. The trajectory $\mathcal{T}_{\varepsilon}$ still
cannot stay in $\mathcal{Q}_{4}.$ As in the case $\alpha\leqq\alpha^{\ast},$
$\mathcal{T}_{\varepsilon}$ has a limit cycle $\mathcal{O}_{\varepsilon}$
surrounding the three stationary points, $w$ is of type (v), and
$\mathcal{T}_{r}$ is oscillating around $0,$ and there exist solutions of type
(vi). Any trajectory $\mathcal{T}_{\left[  P\right]  }$ with $P\not \in
\mathcal{T}_{\varepsilon}$ in $\mathcal{Q}_{4}$ in the domain delimitated by
$\mathcal{O}_{\alpha}$ and $\mathcal{O}_{\varepsilon}$ admits $\mathcal{O}%
_{\alpha}$ as a limit cycle at $-\infty$ and $\mathcal{O}_{\varepsilon}$ at
$\infty$, and $w$ is of type (vii).\medskip

(2) Let $\alpha=\alpha_{c}$ (see figure 15). The homoclinic trajectory
$\mathcal{T}_{\varepsilon}=\mathcal{T}_{\alpha}$ corresponds to the solution
$w$ of type (viii). The trajectory $\mathcal{T}_{r}$ has a limit cycle
$\mathcal{O}_{r}$ surrounding the three points. Thus there exist solutions of
types (iv) or (vi). Any trajectory ending up at $M_{\ell}$ at $\infty$ is
bounded, contained in the domain delimitated by $\mathcal{T}_{\varepsilon},$
and its limit set at $-\infty$ is the homoclinic trajectory $\mathcal{T}%
_{\varepsilon},$ or a cycle around $M_{\ell}$, and $w$ is of type
(ix).\medskip

\begin{theorem}
\label{ent}Assume $\varepsilon=-1,$ and $\alpha_{c}<\alpha<-p^{\prime
}.\medskip$

There exist a unique nonnegative solution $w$ of (\textbf{E}$_{w}$) with a
hole, with $r^{-\gamma}w$ bounded from above and below at $\infty.$ The
regular solutions have at least two zeros.$\medskip$

(1) Either there exist oscillating solutions such that $r^{-\gamma}w$ is
periodic in $\ln r.$ Then the regular solutions have an infinity of zeros, and
$r^{-\gamma}w$ is asymptotically periodic in $\ln r.$ There exist$\medskip$

\noindent(i) solutions satisfying (\ref{vet}) or (\ref{via+}) or (\ref{via+})
or (\ref{viu}) near $0$, oscillating near $\infty,$ such that $r^{-\gamma}w$
is asymptotically periodic in $\ln r;$

\noindent(ii) a unique solution oscillating near 0, such that $r^{-\gamma}w$
is asymptotically periodic in $\ln r,$ and with (\ref{val}) near $\infty;$

\noindent(iii) solutions positive near $0,$ with $r^{-\gamma}w$ bounded, and
oscillating near $\infty,$ such that $r^{-\gamma}w$ is asymptotically periodic
in $\ln r.\medskip$

(2) Or all the solutions have a finite number of zeros, and at least two. Two
cases may occur:$\medskip$

\noindent$\bullet$ Either regular solutions have $m$ zeros and $r^{-\gamma}w$
bounded from above and below at $\infty.$ Then there exist

(iv) solutions with $m$ zeros, with (\ref{vet}) or (\ref{via+}), with
(\ref{val}) near $\infty;$

(v) solutions with $m$ zeros with (\ref{vet}) or (\ref{via+}) and $r^{-\gamma
}w$ bounded from above and below at $\infty;$

(vi) solutions with $m+1$ zeros with (\ref{vet}) or (\ref{via+}) and
$r^{-\gamma}w$ bounded from above and below at $\infty;$

(vii) (for $p>N$) a unique solution with $m$ zeros,with (\ref{viu}) or
(\ref{via-}) and $r^{-\gamma}w$ bounded from above and below at $\infty
.\medskip$

\noindent$\bullet$ Or regular solutions have $m$ zeros and (\ref{val}) holds
near $\infty.$ Then there exist solutions of type (vi) or (vii).
\end{theorem}
\end{proof}

\[%
\begin{array}
[c]{cc}%
\raisebox{-0pt}{\includegraphics[
natheight=187.289200pt,
natwidth=187.289200pt,
height=189.096pt,
width=189.096pt
]%
{../../../AMarie/MesArticles/107050selfp+2/ent1.bmp}%
}%
&
\raisebox{-0pt}{\includegraphics[
natheight=187.289200pt,
natwidth=187.289200pt,
height=189.096pt,
width=189.096pt
]%
{../../../AMarie/MesArticles/107050selfp+2/ent2.bmp}%
}%
\\
\text{th \ref{ent},fig 16: }\varepsilon=-1,N=1,p=3,\alpha=-1.98 & \text{th
\ref{ent}, fig 17: }\varepsilon=-1,N=1,p=3,\alpha=-1.90
\end{array}
\]

\begin{proof}
Here $\mathcal{T}_{\varepsilon}$ stays in $\mathcal{Q}_{4},$ converges to
$M_{\ell}$ or has a limit cycle around $M_{\ell},$ thus $w$ has a hole and
$r^{-\gamma}w$ bounded from above and below at $\infty.$ If $\alpha\geqq
\alpha_{2},$ there is no cycle in $\mathcal{Q}_{4}$, from Proposition
\ref{cyc}, thus $\mathcal{T}_{\varepsilon}$ converges to $M_{\ell}.\medskip$

(1) Either there exists a cycle surrounding $(0,0)$ and $\pm M_{\ell}$, thus
solutions $w$ oscillating around $0,$ such that $r^{-\gamma}w$ is periodic in
$\ln r$. Then $\mathcal{T}_{r}$ has such a limit cycle $\mathcal{O}_{r}$, and
$w$ is oscillating around $0.$ The trajectory $\mathcal{T}_{\alpha}$ has a
limit cycle at $-\infty$ of the same type $\mathcal{O}_{\alpha}\subset
\mathcal{O}_{r}$, and $w$ is of type (ii). For any $P\not \in \mathcal{T}%
_{\varepsilon}$ in the interior domain in $\mathcal{O}_{\alpha},$
$\mathcal{T}_{\left[  P\right]  }$ admits $\mathcal{O}_{\alpha}$ as a limit
cycle at $-\infty$ and converges to $M_{\ell}$ at $\infty,$ or has a limit
cycle in $\mathcal{Q}_{4}$; and $w$ is of type (iii). For any $P\not \in
\mathcal{T}_{r}$, in the domain exterior to $\mathcal{O}_{r},\mathcal{T}%
_{\left[  P\right]  }$ has $\mathcal{O}_{\alpha}$ as limit cycle at $\infty,$
and w is of type (i).\medskip

(2) Or no such cycle exists. Then any trajectory converges at $\infty,$ any
trajectory, apart from $\pm\mathcal{T}_{\alpha},$ converges to $\pm M_{\ell}$
or has a limit cycle in $\mathcal{Q}_{1}$. All the trajectories end up in
$\mathcal{Q}_{2}$ or $\mathcal{Q}_{4}.$ Since $\mathcal{T}_{r}$ starts in
$\mathcal{Q}_{1},$ $y$ has at least one zero. Suppose that it is unique. Then
$\mathcal{T}_{r}$ converges to $-M_{\ell},$ thus $Y$ stays positive. Consider
the function $Y_{\alpha}=e^{(\alpha+\gamma)(p-1)\tau}Y$ defined by (\ref{cge})
with $d=\alpha.$ From Theorem \ref{exlo}, $Y_{\alpha}=(a\left\vert
\alpha\right\vert /N)e^{(\alpha(p-1)+p)\tau}(1+o(1))$ near $-\infty;$ thus
$Y_{\alpha}$ tends to $\infty,$ since $\alpha<p^{\prime}.$ And $Y_{\alpha
}=(\gamma\ell)^{p-1}e^{(\alpha+\gamma)(p-1)\tau}$ near $\infty,$ thus also
$Y_{\alpha}$ tends to $\infty;$ then it has a minimum point $\tau,$ and from
(\ref{phu}), $Y_{\alpha}^{\prime\prime}(\tau)=(p-1)^{2}(\eta-\alpha
)(p^{\prime}+\alpha)Y_{\alpha}<0,$ which is contradictory. Thus $y$ has a
number $m\geqq2$ of zeros$.$ \medskip

Either $\mathcal{T}_{r}\neq\mathcal{T}_{\alpha}.$ Since the slope of
$\mathcal{T}_{\alpha}$ near $-\infty$ is infinite and the slope of
$\mathcal{T}_{r}$ is finite, $\mathcal{T}_{\alpha}$ cuts the line $\left\{
y=0\right\}  $ at $m$ points, starts from $\mathcal{Q}_{1},$ and $w$ is of
type (iv). For any $P$ in the domain of $\mathcal{Q}_{1}$ between
$\mathcal{T}_{r}$ and $\mathcal{T}_{\alpha}$, $\mathcal{T}_{\left[  P\right]
}$ cuts $\left\{  y=0\right\}  $ at $m+1$ points, and $w$ is of type (v). For
any $P$ in the domain of $\mathcal{Q}_{1}$ above $\mathcal{T}_{r},$
$\mathcal{T}_{\left[  P\right]  }$ cuts the line $\left\{  y=0\right\}  $ at
$m+1$ points, and $w$ is of type (vi). If $p>N,$ the trajectories
$\mathcal{T}_{-}$ and $\mathcal{T}_{u}$ cut the line $\left\{  y=0\right\}  $
at $m$ points, and $w$ is of type (vii). \medskip

Or $\mathcal{T}_{r}=\mathcal{T}_{\alpha},$ and then we find only trajectories
with $w$ of type (vi) or (vii).\medskip
\end{proof}

\begin{remark}
Consider the regular solutions in the range $\alpha_{c}<\alpha<-p^{\prime}.$
We conjecture that there exists a decreasing sequence $\left(  \bar{\alpha
}_{n}\right)  ,$ with $\bar{\alpha}_{0}=-p^{\prime}$ and $\alpha_{c}%
<\bar{\alpha}_{n}$ such that for $\alpha\in\left(  \bar{\alpha}_{m}%
,\bar{\alpha}_{m-1}\right)  ,$ $y$ has $m$ zeros and converges to $\pm
M_{\ell};$ and for $\alpha=\bar{\alpha}_{m},$ $y$ has $m+1$ zeros and
converges to $(0,0),$ thus $\mathcal{T}_{r}=\mathcal{T}_{\alpha}.$ We presume
that $\left(  \bar{\alpha}_{m}\right)  $ has a limit $\bar{\alpha}>\alpha
_{c}.$ And for $\alpha<\bar{\alpha},$ $y$ has an infinity of zeros, in other
words there exists a cycle $\mathcal{O}_{r}$ surrounding $\left\{  0\right\}
$ and $\pm M_{\ell}$.$\medskip$

Numerically, for $\alpha=\alpha_{c},$ the cycle $\mathcal{O}_{r}$ seems to be
the unique cycle surrounding the three points. But for $\alpha>\alpha_{c}$ and
$\alpha-\alpha_{c}$ small enough, there exist \textbf{two different cycles}
$\mathcal{O}_{\alpha}\subset\mathcal{O}_{r}$ (see figure 15). As $\alpha$
increases, we observe the coalescence of those cycles; they disappear after
some value $\bar{\alpha}$ (see figure 16).
\end{remark}

\end{document}